\documentclass{elsarticle}

\makeatletter
\def\ps@pprintTitle{%
 \let\@oddhead\@empty
 \let\@evenhead\@empty
 \def\@oddfoot{}%
 \let\@evenfoot\@oddfoot}
\makeatother

\usepackage[utf8]{inputenc}
\usepackage{graphicx}%
\usepackage{multirow}%
\usepackage{amsmath,amssymb,amsfonts}%
\usepackage{amsthm}%
\usepackage{mathrsfs}%
\usepackage[title]{appendix}%
\usepackage{xcolor}%
\usepackage{textcomp}%
\usepackage{manyfoot}%
\usepackage{booktabs}%
\usepackage{algorithm}%
\usepackage{algorithmicx}%
\usepackage{algpseudocode}%
\usepackage{listings}%
\usepackage{verbatim}
\usepackage[hang]{subfigure}
\usepackage{bm}
\usepackage{lmodern}
\usepackage{anyfontsize}
\usepackage{todonotes}
\usepackage{nameref}
\usepackage{mathtools}
\usepackage{appendix}

\newtheorem{remark}{Remark}

\newcommand{\plbr}[1]{ \left( #1 \right) }
\newcommand{\sqbr}[1]{ \left[ #1 \right] }
\newcommand{ \eq }[1]{Eq.~(\ref{eq:#1})}
\newcommand{ \eqs }[2]{Eqs.~(\ref{eq:#1}) and (\ref{eq:#2})}

\newcommand{ \fig }[1]{Figure~\ref{fig:#1}}

\newcommand{ \figto }[2]{Figures~\ref{fig:#1} to \ref{fig:#2}}
\newcommand{ \tab }[1]{Table~\ref{tab:#1}}

\newcommand{ \alg }[1]{Algorithm~\ref{alg:#1}}

\newcommand{\cubr}[1]{ \left\{ #1 \right\} }
\newcommand{\shbr}[1]{ \left| #1 \right| }

\newcommand{ \de }[0]{\, \mathrm{d}}

\newcommand{\curl}[0]{\text{curl}} 
\newcommand{\diverg}[0]{\text{div}} 
\newcommand{\norm}[1]{ \left\lVert #1 \right\rVert } 
\newcommand{\curlcross}[0]{{\mathrm{\times}}} 

\newcommand{\Vspaced}[0]{\mathcal{V}_h} 
\newcommand{\VspacedenhV}[1]{\mathcal{V}^{V#1}_h} 
\newcommand{\Pspace}[0]{\mathcal{P}} 
\newcommand{\Pspaceunc}[0]{\mathcal{P}} 
\newcommand{\Qpolyspace}[0]{\mathcal{Q}} 
\newcommand{\kord}[0]{p} 
\newcommand{\laug}[0]{\ell_E} 
\newcommand{\kordaug}[0]{p+\ell_E} 
\newcommand{\spacedim}[0]{N_d \big|_E} 
\newcommand{\spacedimaug}[0]{N_{d+\ell_E} \big|_E} 
\newcommand{\qpoly}[0]{q} 
\newcommand{\qpolyvec}[0]{\bm{q}} 
\newcommand{\qpolyvecunc}[0]{\bm{q}} 
\newcommand{\globdom}[0]{\mathcal{T}_h} 

\newcommand{\bila}[0]{a} 
\newcommand{\bilad}[0]{a_h} 
\newcommand{\force}[0]{f} 
\newcommand{\forcevec}[0]{\bm{f}} 
\newcommand{\forced}[0]{f_h} 

\newcommand{\unku}[0]{u} 
\newcommand{\testu}[0]{\delta u} 
\newcommand{\unkud}[0]{u_h} 
\newcommand{\testud}[0]{\delta u_h} 

\newcommand{\surom}[0]{\Omega} 
\newcommand{\boundc}[0]{\Gamma} 
\newcommand{\el}[0]{E} 
\newcommand{\diammax}[0]{h} 
\newcommand{\boundel}[0]{\Gamma^\el} 
\newcommand{\bound}[0]{\Gamma^E_e} 
\newcommand{\boundcurv}[0]{\Tilde{\Gamma}^E_e} 
\newcommand{\normedge}[0]{\bm{\hat{n}}_{\Gamma_e}} 
\newcommand{\normedgecurv}[0]{\bm{\hat{n}}_{\Tilde{\Gamma}_e}} 
\newcommand{\normedgecurvtilde}[0]{\Tilde{\bm{\hat{n}}}_{\Tilde{\Gamma}_e}} 

\newcommand{\nv}[0]{n_v} 
\newcommand{\nedges}[0]{n_e} 
\newcommand{\nedgescurv}[0]{\tilde{n}_e} 
\newcommand{\Ne}[0]{N_e} 
\newcommand{\diam}[0]{h_E} 

\newcommand{\trialfcn}[0]{\psi} 

\newcommand{\Pinabla}[0]{\Pi^{\nabla}_\kord} 
\newcommand{\Piok}[0]{\Pi^{0}_\kord} 
\newcommand{\Pioktwo}[0]{\Pi^{0}_{\kord-2}} 
\newcommand{\Piokone}[0]{\Pi^{0}_{\kord-1}} 
\newcommand{\Dm}[0]{\bm{D}} 
\newcommand{\PinablaS}[1]{\Pi^{\mathrm{V#1}}_{\kordaug}} 
\newcommand{\PinablaStwo}[1]{\Pi^{\mathrm{V#1}}_{\kordaug-1}} 
\newcommand{\PiA}[0]{\Pi^{\mathcal{A}}_\kord} 

\newcommand{\KE}[0]{\bm{K}^E} 
\newcommand{\KEc}[0]{\bm{K}_c^E} 
\newcommand{\KEs}[0]{\bm{S}^E} 
\newcommand{\KEscal}[0]{S^E} 

\newcommand{\Acoeff}[0]{\mathcal{A}} 

\newcommand{\disp}[0]{\bm{u}} 
\newcommand{\dispd}[0]{\bm{u}_h} 
\newcommand{\testdisp}[0]{\delta\bm{u}} 
\newcommand{\strain}[0]{\bm{\varepsilon}} 
\newcommand{\cost}[0]{\mathbb{C}} 

\newcommand{\Vspacedvec}[0]{\bm{\mathcal{V}}_h} 

\newcommand{\PiK}[0]{\bm{\Pi}^{\bm{\varepsilon}}_\kord} 
\newcommand{\PiKVC}[0]{\bm{\Pi}^{\mathbb{C}}_\kord} 
\newcommand{\PiokK}[0]{\bm{\Pi}_\kord^{0}} 

\newcommand{\trialfcnvec}[0]{\bm{\psi}} 

\newcommand{\vel}[0]{\bm{u}} 
\newcommand{\testvel}[0]{\delta\bm{ u}} 
\newcommand{\pres}[0]{\rho} 
\newcommand{\testpres}[0]{\delta \rho} 
\newcommand{\veld}[0]{\bm{u}_h} 
\newcommand{\presd}[0]{\rho_h} 

\newcommand{\Qspaced}[0]{Q_h} 
\newcommand{\Gspace}[0]{\mathcal{G}} 
\newcommand{\gpolyperp}[0]{\bm{g}^{\perp}} 

\newcommand{\bilb}[0]{b} 

\newcommand{\Pinablavec}[0]{\bm{\Pi}^{\bm{\nabla}}_\kord} 

\newcommand\bilVC{\mathfrak{a}_\mathsf{vc}^\el}
\newcommand{\PiZvec}[0]{\bm{\Pi}^{\bm{0}}_\kord}
\newcommand\fc{\begin{color}{blue}}
\newcommand\cf{\end{color}}

\begin{document}

\title{Benchmarking stabilized and self-stabilized $\kord$-virtual element methods with variable coefficients}

\author[polimi]{Paola Pia Foligno}
\ead{paolapia.foligno@polimi.it}
\author[kaust,pavia]{Daniele Boffi}
\ead{daniele.boffi@kaust.edu.sa}
\author[kaust]{Fabio Credali\corref{cor1}}
\ead{fabio.credali@kaust.edu.sa}
\author[polimi]{Riccardo Vescovini}
\ead{riccardo.vescovini@polimi.it}

\cortext[cor1]{Corresponding author}

\address[polimi]{Dipartimento di Scienze e Tecnologie Aerospaziali, Politecnico di Milano, 20156 Milan, Italy}
\address[kaust]{CEMSE Division, King Abdullah University of Science and Technology, 23955 Thuwal, Saudi Arabia}
\address[pavia]{Dipartimento di Matematica ``F. Casorati", Universit\`a degli Studi di Pavia, 27100 Pavia, Italy}

\begin{abstract}
Standard Virtual Element Methods (VEM) are based on polynomial projections and require a stabilization term to evaluate the contribution of the non-polynomial component of the discrete space.
However, the stabilization term is not uniquely defined by the underlying variational formulation and is typically introduced in an ad hoc manner, potentially affecting the numerical response.
Stabilization-free and self-stabilized formulations have been proposed to overcome this issue, although their theoretical analysis is still less mature.

This paper provides an in-depth numerical investigation into different stabilized and self-stabilized formulations for the $\kord$-version of VEM.
The results show that self-stabilized and stabilization-free formulations achieve optimal accuracy while suffering from worse conditioning.
Moreover, a new projection operator, which explicitly accounts for variable coefficients, is introduced within the framework of standard virtual element spaces.
Numerical results show that this new approach is more robust than the existing ones for large values of $\kord$. 

\end{abstract}

\maketitle

\section{Introduction} \label{sec:Introduction}

Current research studies focus on the design of novel aerospace structures with improved efficiency. An example is given by variable stiffness panels, which are obtained by combining a variable stiffness skin, in which the fibers' orientation varies across the domain, with stringers running along curvilinear paths.
These configurations are modeled by partial differential equations with non-uniform elastic properties over the panel domain. Numerical approaches for the mentioned application have been proposed in the framework of finite elements \cite{zhao2019prestressed} and Ritz method \cite{vescovini2020semi}. However, in a finite element context, node sharing between plate and stiffeners requires careful meshing or ad-hoc formulations \cite{zhao2017vibration}. On the contrary, a Ritz approach simplifies the enforcement of the compatibility between plate and stringers, but is generally not versatile in dealing with complex geometries, unless specific procedures are employed~\cite{vescovini2023ritz}.

In this landscape, the Virtual Element Method (VEM)~\cite{daveiga2013basic} is a promising tool since it supports the use of meshes made of general polytopal elements, with possibly curved edges~\cite{daveiga2019virtual2}. Thus, meshes conforming with the stringers can be easily generated and combined with high-order discretizations.

Despite VEM has already been applied to several problems in computational mechanics, such as plane stress elasticity~\cite{daveiga2013virtual,artioli2017arbitrary,artioli2017stress,mengolini2019engineering,daltri2020error}, buckling problems~\cite{mora2020virtual,meng2020linear,adak2023c0},  Kirchhoff--Love~\cite{brezzi2013virtual} and Reissner--Mindlin plate models~\cite{daveiga2019virtual,daltri2022first}, to the best of the authors' knowledge its application to variable stiffness panels remains limited. A contribution can be found in \cite{reddy2019virtual}, in which a piecewise constant approximation of the variable coefficient is considered. For elasticity problems, spatially varying material coefficients have been the subject of other investigations in a VEM framework based on the Hellinger–Reissner variational principle. 
In the two-dimensional setting of \cite{artioli-var}, a local energy projection operator is defined, leading to elementwise approximation spaces that are no longer polynomial in the presence of variable elasticity coefficients. A similar approach is proposed in \cite{visinoni}. The formulation is here three-dimensional and an analogous projection-based approach is used for polyhedral meshes. In this context, the present paper serves as a preliminary investigation to benchmark virtual elements in their $\kord$-version~\cite{daveiga2016basic} by taking into account the challenges posed by variable stiffness panels. 

VEM can be interpreted as a generalization of finite elements on polytopal meshes. However, virtual elements give rise to a nonconforming method as the trial and test functions are solutions of local PDEs (defined in each mesh element) that are never explicitly solved when constructing the discrete formulation. Hence the name ``Virtual". Since only a polynomial projection of the trial functions is computable through the degrees of freedom, a stabilization term is introduced to handle the residual non-polynomial component and ensure the stability of the discrete problem. However, the stabilization term has an arbitrary nature and does not generally conform to the physics of the problem. In the particular case of variable coefficients, this term is scaled by a constant approximation and the polynomial projection is built either independently of the coefficient~\cite{daveiga2016virtual} or by taking into account a low-order polynomial approximation thereof~\cite{daveiga2014virtual}.

The stabilization involves the introduction of a user-defined scaling parameter, which reduces the generality of the method, potentially leading to over-stabilization and degraded solution quality: this situation is verified for anisotropic problems~\cite{berrone2022comparison} and eigenvalue problems~\cite{boffi2020approximation,alzaben2025stabilization}, where an improper choice of the stability parameters may trigger the onset of spurious modes. The relation between the stabilization parameter and the performance of VEM has been discussed in~\cite{cangiani2015hourglass,russo2022quantitative}. Despite the many stabilization forms available in the literature~\cite{daveiga2013basic,daveiga2017stability,daveiga2018exponential,mascotto2018ill,bertoluzza2022stabilization}, most of them are suitable for low-order virtual elements, while they tend to be suboptimal as the order increases. 

An approach to overcome the issues caused by the stabilization term is given by \textit{stabilization-free} virtual elements \cite{berrone2025lowest,berrone2022comparison,berrone2025stabilization}. These elements are obtained by projecting the virtual functions onto larger polynomial spaces. So, the resulting discrete problem is automatically stable. A different strategy refers to the so-called \textit{self-stabilized} elements, which are obtained by enriching the approximation space with extra internal degrees of freedom~\cite{lamperti2023hu}. 
The term \textit{self-stabilized VEM} will refer to both methods throughout this paper. 
Although not addressed here, the work in \cite{credali2024reduced} presents a reduced basis approach to cheaply reconstruct the lowest-order VEM basis functions and avoid the need of the stabilization term.

The present paper is divided into two main parts. First, the performance of several stabilized and self-stabilized formulations is systematically assessed within the framework of $\kord$-VEM. For this purpose, Laplace, linear elasticity, and Stokes problems are considered. The aim is to provide insights on the accuracy, stability, conditioning, and sensitivity to stabilization parameters on structured and distorted meshes, even with curved edges. To the best of the authors' knowledge, no previous investigations can be found on this topic with respect to $\kord$-VEM.

The second element of novelty of this paper regards a new approach for dealing with problems governed by variable coefficients in high-order approximations. 
Within the framework of well-known virtual element spaces, a new polynomial projection is introduced so that the features of the coefficient are taken into account. This formulation, denoted as $\mathrm{VC}$-VEM, can be easily embedded into stabilized and self-stabilized settings, leading to improved accuracy.

The work is organized as follows. The geometrical framework is presented in Section~\ref{sec:Domain_general}. The virtual element discretization of the Laplace equation by means of stabilized and self-stabilized methods is recalled in Section~\ref{sec:Laplace_SelfStab}, while the linear elasticity problem and the Stokes problem are discussed in Sections~\ref{sec:LinearElasticity_SelfStab} and \ref{sec:Stokes_SelfStab}, respectively. Section~\ref{sec:NumericalResults_SelfStab} is devoted to the benchmark of $\kord$-VEM across a wide range of numerical tests. Then, the $\mathrm{VC}$-VEM is presented in Section~\ref{sec:VC} for second-order elliptic and linear elasticity problems. The numerical results are presented in Section~\ref{sec:NumericalResults_VC}, where the robustness of the proposed element is shown. Conclusions are drawn in Section~\ref{sec:Conclusions}.

\section{Geometrical framework} \label{sec:Domain_general}

The domain is denoted by $\surom$ and $\cubr{\boundc_i}_{i=1,\dots,\Ne}$ is a finite set of smooth curves identifying the boundary $\boundc := \partial \surom$, where $\Ne$ is the number of edges. 
By assuming that the boundary $\boundc$ is Lipschitz and that for $m \geq 0$ each curve $\boundc_i$ is of class $\mathcal{C}^{m+1}$ \cite{daveiga2019virtual2}, one can define the regular and invertible $\mathcal{C}^{m+1}$ parametrization of $\boundc_i$ with respect to $I_i = \sqbr{a_i,b_i}$ as $\gamma_i:I_i \rightarrow \boundc_i$. 
As $\gamma_i$ is an invertible mapping, for a generic function $g$:
\begin{equation}
    \begin{aligned}
        & g = \Tilde{g} \circ \gamma_i && \forall \Tilde{g} \in \boundc_i,   \\
        & \Tilde{g} = g \circ \gamma_i^{-1} && \forall  g \in  I_i.   
    \end{aligned}
    \label{eq:mapping_gamma}
\end{equation}

In this work, curved edges are parametrized using B\'ezier curves. A generic B\'ezier curve is expressed as:
\begin{equation}
\begin{aligned}
\bm{r}\plbr{t} = \sum_{i=0}^n \binom{n}{i} \bm{P}_i \plbr{1-t}^{n-i} t^i = \sum_{i=0}^n \bm{P}_i \, B_{i,n}\plbr{t}, \quad t \in [0 \;  1],
\end{aligned}
\label{eq:bezier}
\end{equation}
where $n$ is the polynomial order, $\bm{P}_i$ are the control points coordinates, and $B_{i,n}\plbr{t}$ is a Bernstein polynomial defined as:
\begin{equation}
B_{i,n}\plbr{t} = \binom{n}{i} \plbr{1-t}^{n-i}.
    \label{eq:bernstein}
\end{equation}
The B\'ezier curve passes through the first and last control points and is enclosed by the convex hull of the control points.

The domain $\surom$ is partitioned by a finite set $\globdom$ of non-overlapping star-shaped polygons with possibly curved edges $\el \in \globdom$, as in the standard VEM theory \cite{daveiga2013basic}. The symbol $\diam$ denotes the diameter of $\el$ and $\diammax=\max_{\el \in \globdom} \diam$. The boundary of $\el$ is denoted by $\boundel = \partial \el$. More precisely, $\boundel$ is the union of straight segments $\bound$ with $e=1,\dots,\nedges$ and curved edges $\boundcurv$ with $e=1,\dots,\nedgescurv$, satisfying the regularity assumptions mentioned above.

\section{Laplace problem} \label{sec:Laplace_SelfStab}

In this section, the stabilized and self-stabilized formulations of the Virtual Element Method are presented for the Laplace problem, whose weak formulation reads:
\begin{equation}
    \begin{cases}
        \text{find } \unku \in H_0^1 \plbr{\surom} \text{ such that:}\\
        \bila \plbr{\testu,\unku} = \plbr{\testu,\force}_{\surom} \quad \forall \testu \in H_0^1 \plbr{\surom},
    \end{cases}
\label{eq:laplace_goveqs}
\end{equation}
where $\bila \plbr{\unku,\testu}$ represents the internal energy and $\plbr{\testu,\force}_{\surom}$ represents the work of external sources:
\begin{equation}
    \begin{aligned}
        \bila \plbr{\testu,\unku} = \int_{\surom} \nabla \testu^T \nabla \unku \, \de \surom,\qquad \plbr{\testu,\force}_{\surom} = \int_{\surom} \testu \, \force \, \de \surom.
    \end{aligned}
\label{eq:laplace_bilabilf}
\end{equation}

Throughout the paper it is used the compact notation $\plbr{\cdot,\cdot}_D$ to denote the scalar product in the space $L^2 \plbr{D}$.

A discrete space $\Vspaced \subset H_0^1 \plbr{\surom}$ is then introduced, and the discrete counterpart of the bilinear form $\bila \plbr{\cdot,\cdot}$ is constructed as the sum of local contributions over $\globdom$:
\begin{equation}
    \bilad \plbr{\testud,\unkud} = \sum_{\el \in \globdom} \bilad^\el\plbr{\testud,\unkud} \quad \forall \testud,\unkud \in \Vspaced.
    \label{eq:discrete_bilinearform}
\end{equation}
By performing the same operation for the right hand side, one obtains the discrete problem:
\begin{equation}
    \begin{cases}
        \text{find } \unkud \in \Vspaced \text{ such that:}\\
        \bilad\plbr{\testud,\unkud} = \plbr{\testud, \forced}_\surom \quad \forall \testud \in \Vspaced.
    \end{cases}
\label{eq:discrete_problem}
\end{equation}

Several definitions of $\Vspaced$ and of the related discrete problems are recalled in the next subsections.

\subsection{Stabilized VEM} \label{subsec:Laplace_stabilization}

The local virtual element space of order $\kord \geq 1$ is defined on a generic polygon $\el \in \globdom$, for $\kord \leq m$, as:
\begin{equation}
    \begin{aligned}
    \Vspaced \plbr{\el} = \left\{ \right. &\unkud \in H^1 \plbr{\el} \cap C^0 \plbr{\el}  :\Delta \unkud \big|_\el \in \Pspace_{\kord} \plbr{\el},\\  & \unkud \big|_{\bound} \in \Pspace_\kord \plbr{\bound} \forall  e=1,\dots,\nedges, \\
    &\left. \unkud \big|_{\boundcurv} \in \Tilde{\Pspace}_\kord \plbr{\boundcurv} \forall e=1,\dots,\nedgescurv, \right.\\
    &\plbr{\unkud,\qpoly}_\el = \plbr{\Pinabla \unkud,\qpoly}_\el \; \forall \qpoly \in \Pspace_\kord \plbr{\el} \setminus \Pspace_{\kord-2} \plbr{\el}\left. \right\}.
    \end{aligned}
    \label{eq:VEMlocalspace_enhanced}
\end{equation} 

When $\nedgescurv=0$, i.e. $\el$ has only straight edges, the inclusion ${\Pspace_{\kord} \plbr{\el} \subset \Vspaced \plbr{\el}}$ guarantees the accuracy of the method. On the other hand, for $\nedgescurv>0$, ${\Pspace_{\kord} \plbr{\el} \not\subset \Vspaced \plbr{\el}}$ and only ${\Pspace_{0} \plbr{\el} \subset \Vspaced \plbr{\el}}$ holds: in this case, the accuracy is guaranteed by the fact that curved edges are a feature of the domain and do not depend on~$h$, see~\cite{daveiga2019virtual2} for more details.
The operator $\Pinabla$ maps the discrete functions from the VEM space to the polynomial one and will be defined later. Conventionally, $\Pspace_{-1} \plbr{\el} = \cubr{0}$.

Moreover, the polynomial space $\Tilde{\Pspace}_\kord \plbr{\boundcurv}$ on the curved edges is defined by using \eq{mapping_gamma} \cite{daveiga2019virtual2}:
\begin{equation}
    \Tilde{\Pspace}_\kord \plbr{\boundcurv} = \cubr{ \Tilde{\qpoly} = \qpoly \circ \gamma_e^{-1}: \qpoly \in \Pspace_\kord \plbr{I_e}},
    \label{eq:curvedpoly_space}
\end{equation}
so that the functions on the curved edges are polynomials through the parametrization $\gamma_e$.

The VEM space has dimension $\spacedim$, and $\unkud \in \Vspaced \plbr{\el}$ is uniquely identified by the following unisolvent set of degrees of freedom:
\begin{itemize}
\item $\textbf{DOF}_{1}$: the values at the $\nv$ vertices of $\el$,
\item $\textbf{DOF}_{2}$: for $\kord > 1$, the values at the $\kord-1$ internal Gauss-Lobatto quadrature points on each $\bound$ and the values of at the $\kord-1$ points on each $\boundcurv$ that are images through $\gamma_e$ of the $\kord-1$ internal Gauss-Lobatto quadrature points on $I_e$ \cite{daveiga2019virtual2},
\item $\textbf{DOF}_{3}$ : for $\kord > 1$, the internal moments up to order $\kord-2$:
\begin{equation}
    \frac{1}{\shbr{\el}} \int_\el \qpoly \plbr{x,y} \unkud \plbr{x,y} \de \el \quad \forall \qpoly \in \Qpolyspace_{\kord-2}\plbr{\el},
    \label{eq:VEMlocalmoments}
\end{equation}
where $\shbr{\el}$ is the area of the element $\el$ and $\Qpolyspace_{\kord}\plbr{\el}$ is a basis for the polynomial space $\Pspace_{\kord}\plbr{\el}$.
\end{itemize}

The global virtual element space is then obtained by gluing all the local spaces by continuity:
\begin{equation}
    \Vspaced = \cubr{ \unkud \in H_0^1\plbr{\surom} : \unkud{\big|_\el} \in \Vspaced \plbr{\el} \quad \forall \el \in \globdom }.
    \label{eq:VEMglobalspace}
\end{equation}

The generic function $v_h \in \Vspaced \plbr{\el}$ can be expressed via standard Lagrangian  basis functions $\trialfcn_i$ as:
\begin{equation}
    v_h = \sum_{i=1}^{\spacedim} \text{dof}_i \plbr{v_h} \trialfcn_i \quad \forall v_h \in \Vspaced \plbr{\el},
    \label{eq:lagrangeinterpolation}
\end{equation}
where $\text{dof}_j \plbr{\trialfcn_i} = \delta_{ij}$ for $i,j=1,\dots,\spacedim$.
It is clear that the discrete functions are not explicitly known in the interior of $\el$. Therefore, the polynomial projection $\Pinabla : \Vspaced \plbr{\el} \rightarrow \Pspace_\kord \plbr{\el}$ is required to construct the discrete bilinear form. In particular, $\Pinabla \trialfcn_i$ is defined as the solution of:
\begin{equation}
    \begin{cases}
    \begin{aligned}
        &\bila^\el \plbr{ \Pinabla \trialfcn_i, \qpoly } = \bila^\el \plbr{ \trialfcn_i, \qpoly } && \forall \qpoly \in \Pspace_\kord \plbr{\el},\\
        &P_0^\el \plbr{ \Pinabla \trialfcn_i, \qpoly } = P_0^\el \plbr{ \trialfcn_i, \qpoly } && \forall \qpoly \in \Pspace_0 \plbr{\el},
    \end{aligned}    
    \end{cases}
    \label{eq:VEMH1projection}
\end{equation}
for all $i=1,\dots,\spacedim$, where:
\begin{equation}
    P_0^\el \plbr{ \trialfcn_i, \qpoly } = \frac{1}{\nv} \sum_{j=1}^{\nv} \text{dof}_j \plbr{ \trialfcn_i } \text{dof}_j \plbr{ \qpoly },
    \label{eq:laplace_invcondition2}
\end{equation}
is required to prevent the system to be singular.

The left hand side of \eq{VEMH1projection} is the product of known polynomials, so it is readily available. 
On the contrary, integration by parts is applied to evaluate the right hand side of the equation:
\begin{equation}
\begin{aligned}
    \bila^\el \plbr{ \trialfcn_i, \qpoly } &= - \int_\el \trialfcn_i \Delta \qpoly \; \de \el +  \sum_{e=1}^{\nedges} \int_{\bound} \trialfcn_i \cdot \nabla \qpoly \, \normedge \; \de \bound\\
    &+\sum_{e=1}^{\nedgescurv} \int_{\boundcurv} \trialfcn_i \cdot \nabla \qpoly \,\normedgecurv \; \de \boundcurv,
\end{aligned}    
    \label{eq:B_matrix_integrationbyparts}
\end{equation}
The volume integral of the equation above is computed from the internal degrees of freedom, while the boundary integrals on the straight edges are evaluated by standard quadrature rules for polynomials.
On the other hand, the discrete functions $\trialfcn_i$ are not polynomials on curved edges, so that the parametrization $\gamma_e$ defined in \eq{mapping_gamma} is used to recover the polynomial structure on the reference interval $I_e$ \cite{daveiga2019virtual2}:
\begin{equation}
    \begin{aligned}
        \sum_{e=1}^{\nedgescurv} \int_{\boundcurv} \trialfcn_i \cdot \nabla \qpoly \,\normedgecurv \; \de \boundcurv  &= \sum_{e=1}^{\nedgescurv} \int_{\boundcurv} \sqbr{\Tilde{\trialfcn_i} \cdot \nabla \Tilde{\qpoly} \; \normedgecurvtilde} \circ \gamma_e^{-1}   \de \boundcurv \\
        &= \sum_{e=1}^{\nedgescurv} \int_{I_e} \sqbr{\Tilde{\trialfcn_i} \cdot \nabla \Tilde{\qpoly} \;\normedgecurvtilde}  \norm{\gamma_e^{'}}   \de I_e,\\
    \end{aligned}
    \label{eq:B_matrix_curve}
\end{equation}
where $\Tilde{\trialfcn_i}$ is the image of $\trialfcn_i$ through $\gamma_e$.
Once the projection $\Pinabla$ is computed, the following decomposition holds:
\begin{equation}
    \trialfcn_i = \Pinabla \trialfcn_i + \plbr{ I - \Pinabla } \trialfcn_i.
    \label{eq:trial_fcn_subdivision}
\end{equation}
By leveraging the definition of $\Pinabla$, the bilinear form reads:
\begin{equation}
    \begin{aligned}
        \bila^\el \plbr{\trialfcn_i, \trialfcn_j } = 
        \bila^\el \plbr{ \Pinabla \trialfcn_i, \Pinabla \trialfcn_j } + \bila^\el \plbr{ \plbr{ I - \Pinabla } \trialfcn_i, \plbr{ I - \Pinabla } \trialfcn_j },
    \end{aligned}
    \label{eq:el_stiffness}
\end{equation}
where the first term at the right-hand side is fully computable through the degrees of freedom, while the second one is usually replaced by a suitable stabilization term, based on a symmetric and semi-positive definite bilinear form $\KEscal$, satisfying:
\begin{equation}
    \alpha_*(p)\,\bila^\el\plbr{\testud,\testud} \leq \KEscal \plbr{\testud,\testud} \leq \alpha^*(p)\, \bila^\el\plbr{\testud,\testud},
    \label{eq:stability}
\end{equation}
for all $\testud \in \Vspaced \plbr{\el}$ with $\Pinabla \testud = 0$ and for positive constants $\alpha_*$ and $\alpha^*$, independent of $\diammax$.
The local virtual element bilinear form is then defined as:
\begin{equation}
    \begin{aligned}
        \bilad^\el \plbr{\trialfcn_i, \trialfcn_j } = 
        \bila^\el \plbr{ \Pinabla \trialfcn_i, \Pinabla \trialfcn_j } + \tau \KEscal \plbr{ \plbr{ I - \Pinabla } \trialfcn_i, \plbr{ I - \Pinabla } \trialfcn_j },
    \end{aligned}
    \label{eq:el_stiffness_discrete}
\end{equation}
where $\tau$ is a user-defined parameter.

Due to the arbitrary nature of the stabilization term, several definitions have been introduced in the literature, also depending on the problem under consideration. The most popular ones are now recalled in view of the numerical tests presented later.

The first three stabilization techniques are the most trivial to be implemented and their expression reads:
\begin{equation}
\begin{aligned}
     &\KEscal_1 \plbr{\trialfcn_i,\trialfcn_j} = \sum_{r=1}^{\spacedim} \text{dof}_r \plbr{\plbr{I-\Pinabla}\trialfcn_i} \text{dof}_r \plbr{\plbr{I-\Pinabla}\trialfcn_j}, \\
     &\KEscal_2 \plbr{\trialfcn_i,\trialfcn_j} = \sum_{r=1}^{\kord \nv} \text{dof}_r \plbr{\plbr{I-\Pinabla}\trialfcn_i} \text{dof}_r \plbr{\plbr{I-\Pinabla}\trialfcn_j}, \\
     &\KEscal_3 \plbr{\trialfcn_i,\trialfcn_j} =  \sum_{r=1}^{\spacedim}  \omega_r \, \text{dof}_r \plbr{\plbr{I-\Pinabla}\trialfcn_i} \text{dof}_r \plbr{\plbr{I-\Pinabla}\trialfcn_j},
\end{aligned}     
    \label{eq:stabilization_S123}
\end{equation}
where $\omega_r=\max \cubr{1,\bila^\el\plbr{\Pinabla\trialfcn_r,\Pinabla\trialfcn_r}}$. Observe that $\KEscal_2$ only accounts for the boundary degrees of freedom, whereas $\KEscal_3$ scales as the consistency part, making it more robust than $\KEscal_1$ and $\KEscal_2$.

Another scheme specifically designed for $\kord$-VEM is:
\begin{equation}
\begin{aligned}
  \KEscal_4 \plbr{\trialfcn_i,\trialfcn_j} &=  \frac{\kord}{\diam} \sum_{j=1}^{\nedges} \int_{\bound} \trialfcn_i \trialfcn_j \; \de \bound + \frac{\kord}{\diam} \sum_{j=1}^{\nedgescurv} \int_{\boundcurv} \trialfcn_i \trialfcn_j \; \de \boundcurv\\ &+  \frac{\kord^2}{\diam^2} \int_{\el} \Pioktwo \trialfcn_i \Pioktwo \trialfcn_j \; \de \el,
    \label{eq:stabilization_S4}
\end{aligned}    
\end{equation}
where $\Piok : \Vspaced \plbr{\el} \rightarrow \Pspace_\kord \plbr{\el}$ is the $L^2$  projection onto polynomials.

The last stabilization considered here refers to \cite{artioli2017arbitrary,daltri2022first} and is particularly effective in linear elasticity problems. It is usually expressed in matrix form as:
\begin{equation}
    \KEs_5 = \bm{I} - \Dm \plbr{\Dm^T\Dm}^{-1}\Dm^T,
    \label{eq:stabilization_S5}
\end{equation}
where $\Dm_{i,\alpha} = \text{dof}_i \plbr{ \qpoly_{\alpha} }$ for $i=1,\dots,\spacedim$ and for all $\qpoly_\alpha \in \Qpolyspace \plbr{\el}$.

To conclude the construction of the discrete problem \eq{discrete_problem}, the right hand side can be computed as:
\begin{equation}
    \plbr{\testud, \forced}_\surom = \sum_{\el \in \globdom} \plbr{\testud, \Piok \force}_\el.
    \label{eq:bodyforces_vector}
\end{equation}

\subsection{Self-stabilized VEM}\label{subsec:Laplace_selfstabilization}

As detailed in the previous section, the arbitrary nature of stabilization terms motivates the seek for more general formulations, not relying upon heuristic terms. Thus, self-stabilized formulations have been proposed as a viable strategy to prevent the use of the stabilization term. These approaches are particularly useful when dealing with high-order approximations and more complex problems.

The construction of a self-stabilized bilinear form requires in general an enlargement of the underlying virtual element space, so that a higher-order polynomial projection $\Pi_{\kord^\star}$ can be computed, with $\kord^\star > \kord$. The discrete bilinear form is then defined either as:
\begin{equation}
    \begin{aligned}
        \bilad^\el \plbr{\trialfcn_i, \trialfcn_j } = 
        \plbr{ \Pi_{\kord^\star} \nabla \trialfcn_i, \Pi_{\kord^\star} \nabla \trialfcn_j }_\el,
    \end{aligned}
    \label{eq:el_stiffness_discrete_ss}
\end{equation}
or:
\begin{equation}
    \begin{aligned}
        \bilad^\el \plbr{\trialfcn_i, \trialfcn_j } = 
        \bila^\el \plbr{ \Pi_{\kord^\star} \trialfcn_i, \Pi_{\kord^\star} \trialfcn_j },
    \end{aligned}
    \label{eq:el_stiffness_discrete_ss_L2}
\end{equation}
depending on the definition of the projection. Conversely, the right hand side preserved the expression provided in \eq{bodyforces_vector}.

The six different self-stabilization strategies described hereafter are summarized in \tab{selfstabVEM}. The first three and the last one rely on the discrete bilinear form defined in \eq{el_stiffness_discrete_ss_L2}, while the fourth and the fifth one refer to \eq{el_stiffness_discrete_ss}.

\renewcommand{\arraystretch}{1.25}
\begin{table}[]
    \centering
    \begin{tabular}{c l l l}
         \textbf{Version} &  \textbf{Reference} & \textbf{VEM space type} & \textbf{Notes}\\
         \hline
         V1 & Berrone et al.~\cite{berrone2025lowest} & Enlarged enhanced&\\
         V2 & Lamperti et al.~\cite{lamperti2023hu} & Augmented  &Additional dofs\\
         V3 & Berrone et al.~\cite{berrone2025stabilization} & Standard& Div-free polynomials\\
         V4 & \textit{inspired by}~\cite{berrone2025lowest} & Enlarged enhanced&\\
         V5 & \textit{inspired by}~\cite{lamperti2023hu} & Augmented &Additional dofs \\
         V6 & \textit{inspired by}~\cite{berrone2025stabilization} & Standard&Div-free polynomials\\
         \hline
         &&&\\
         & \textbf{Projection type} & \multicolumn{2}{c}{\textbf{Bilinear form type}}\\
         \hline
         V1--V3,V6&$L^2$ projection&\multicolumn{2}{c}{$\bilad^\el \plbr{\trialfcn_i, \trialfcn_j } = 
        \plbr{ \Pi_{\kord^\star} \nabla \trialfcn_i, \Pi_{\kord^\star} \nabla \trialfcn_j }_\el$}\\
         V4--V5&Elliptic projection&\multicolumn{2}{c}{$\bilad^\el \plbr{\trialfcn_i, \trialfcn_j } = 
        \bila^\el \plbr{ \Pi_{\kord^\star} \trialfcn_i, \Pi_{\kord^\star} \trialfcn_j }$}\\
        \hline
    \end{tabular}
    \caption{Summary of the considered self-stabilized VEM formulations.}
    \label{tab:selfstabVEM}
\end{table}
\renewcommand{\arraystretch}{1}

\subsubsection*{V1: Enlarged enhanced space with $L^2$ projection}

The first self-stabilized technique was introduced in \cite{berrone2025lowest} for the lowest-order VEM and here extended to the generic order $\kord$.

The original virtual element space in \eq{VEMlocalspace_enhanced} is enlarged as:
\begin{equation}
    \begin{aligned}
    \VspacedenhV{1} \plbr{\el} = \left\{ \right. &\unkud \in H^1 \plbr{\el} \cap C^0 \plbr{\el}  :\Delta \unkud \big|_\el \in \Pspace_{\kord+\laug} \plbr{\el},\\
    & \unkud \big|_{\bound} \in \Pspace_\kord \plbr{\bound} \forall e=1,\dots,\nedges, \\
    &\unkud \big|_{\boundcurv} \in \Tilde{\Pspace}_\kord \plbr{\boundcurv} \forall e=1,\dots,\nedgescurv\\
    &\plbr{\unkud,\qpoly}_\el = \plbr{\Pinabla \unkud,\qpoly}_\el \;\forall \qpoly \in \Pspace_{\kord+\laug} \plbr{\el} \setminus \Pspace_{\kord-2} \plbr{\el}\left. \right\},
    \end{aligned}
    \label{eq:VEMlocalspace_selfstab_enhanced_V1}
\end{equation} 
where $\laug \geq 1$ is the additional polynomial order. 
In this case, rather than computing the polynomial projection of the trial functions $\trialfcn_i$, the projection of $\nabla \trialfcn_i$ is considered instead.
The projection, $\PinablaStwo{1}$, is defined as follows:
\begin{equation}
\begin{aligned}      
    \plbr{\PinablaStwo{1}\nabla \trialfcn_i,\qpolyvecunc}_\el = \plbr{\nabla \trialfcn_i,\qpolyvecunc}_\el \quad \forall \qpolyvecunc \in \sqbr{ \Pspaceunc_{\kordaug-1} \plbr{\el}}^2,
\end{aligned}
\label{eq:SS4_projections}
\end{equation}
for all $i=1,\dots,\spacedim$. Note, the space $\VspacedenhV{1} \plbr{\el}$ has the same dimension of $\Vspaced \plbr{\el}$.

\subsubsection*{V2: Augmented space with $L^2$ projection}

This formulation refers to the one presented in \cite{lamperti2023hu} for the lowest-order VEM for linear elasticity problems and here extended to generic order.

If the additional polynomial order is $\laug > 1$, the augmented space is defined as:
\begin{equation}
    \begin{aligned}
    \VspacedenhV{2} \plbr{\el} = \left\{ \right. &\unkud \in H^1 \plbr{\el} \cap C^0 \plbr{\el}  :\Delta \unkud \big|_\el \in \Pspace_{\kord+\laug-2} \plbr{\el},\\
    & \unkud \big|_{\bound} \in \Pspace_\kord \plbr{\bound} \forall e=1,\dots,\nedges, \\
    &\unkud \big|_{\boundcurv} \in \Tilde{\Pspace}_\kord \plbr{\boundcurv} \forall e=1,\dots,\nedgescurv\left.\right\},
    \end{aligned}
    \label{eq:VEMlocalspace_selfstab_enhanced_V2}
\end{equation} 
so that its dimension is:
\begin{equation}
    \spacedimaug = \text{dim} \plbr{\VspacedenhV{2} \plbr{\el}} = \kord \nv + \frac{\plbr{\kordaug-1}\plbr{\kordaug} }{2},
    \label{eq:VEMlocalspacedim_V2}
\end{equation}
and the following additional degrees of freedom are considered:
\begin{itemize}
\item $\textbf{DOF}_{4}$ :
\begin{equation}
    \frac{1}{\shbr{\el}} \int_\el \qpoly \plbr{x,y} \unkud \plbr{x,y} \de \el \quad \forall \qpoly \in \Qpolyspace_{\kordaug-2}\plbr{\el} \setminus \Qpolyspace_{\kord-2}\plbr{\el}  .
    \label{eq:VEMlocalmoments_aug}
\end{equation}
\end{itemize}

Note, if $\laug=1$, from \eq{VEMlocalspace_selfstab_enhanced_V2} one obtains that $\Delta \unkud \big|_\el \in \Pspace_{\kord-1} \plbr{\el}$, which results in a weaker condition than that of the standard enhanced space in~\eq{VEMlocalspace_enhanced}. Hence, for the particular case of $\laug=1$, the following variant of $\VspacedenhV{2} \plbr{\el}$ is considered to guarantee the computability of the $L^2$ product with polynomials of degree $\kord$:
\begin{equation}
    \begin{aligned}
    \VspacedenhV{2} \plbr{\el} = \left\{ \right. &\unkud \in H^1 \plbr{\el} \cap C^0 \plbr{\el}  :\Delta \unkud \big|_\el \in \Pspace_{\kord} \plbr{\el},\\  & \unkud \big|_{\bound} \in \Pspace_\kord \plbr{\bound} \forall  e=1,\dots,\nedges, \\
    &\left. \unkud \big|_{\boundcurv} \in \Tilde{\Pspace}_\kord \plbr{\boundcurv} \forall e=1,\dots,\nedgescurv, \right.\\
    &\plbr{\unkud,\qpoly}_\el = \plbr{\Pinabla \unkud,\qpoly}_\el \; \forall \qpoly \in \Pspace_\kord \plbr{\el} \setminus \Pspace_{\kord-1} \plbr{\el}\left. \right\}.
    \end{aligned}
    \label{eq:VEMlocalspacedim_V2_l1}
\end{equation} 
As a consequence, internal moments of order $\kord-1$ are added in combination with the standard enhanced VEM condition for computing the $L^2$ product with polynomials of degree $\kord$.

The same projection of $\mathrm{V1}$ is now computed by exploiting the additional dofs defined in \eq{VEMlocalspacedim_V2}:
\begin{equation}
\begin{aligned}     
        \plbr{\PinablaStwo{2}\nabla \trialfcn_i,\qpolyvecunc}_\el = \plbr{\nabla \trialfcn_i,\qpolyvecunc}_\el \quad \forall \qpolyvecunc \in \sqbr{ \Pspaceunc_{\kordaug-1} \plbr{\el}}^2,
\end{aligned}
\label{eq:SS5_projections}
\end{equation}
for all $i=1,\dots,\spacedimaug$.

\subsubsection*{V3: Divergence-free space with $L^2$ projection}

The third self-stabilized formulation is taken from \cite{berrone2025stabilization} and refers to the same VEM space of the classical stabilized VEM.
The polynomial projection is then computed onto the space: $$\sqbr{\Pspace_{\kord-1} \plbr{\el}}^2 \oplus \curl \plbr{\Pspace_{\kordaug} \plbr{\el} \setminus \Pspace_{\kord} \plbr{\el}},$$ which has a divergence-free component.
The projection reads:
\begin{equation}
\begin{aligned}      
    &\begin{cases}
        \begin{aligned}
            &\plbr{\PinablaStwo{3}\nabla \trialfcn_i,\qpolyvecunc}_\el = \plbr{\nabla \trialfcn_i,\qpolyvecunc}_\el  && \forall \qpolyvecunc \in \sqbr{ \Pspaceunc_{\kord-1} \plbr{\el}}^2,\\
            &\plbr{\PinablaStwo{3}\nabla \trialfcn_i,\curl \plbr{\qpoly}}_\el = \plbr{\nabla \trialfcn_i,\curl \plbr{\qpoly}}_\el && \forall \qpoly \in \Pspace_{\kordaug} \plbr{\el} \setminus \Pspace_{\kord} \plbr{\el},
        \end{aligned}
    \end{cases}  
\end{aligned}    
\label{eq:SS6_projections}
\end{equation}
for all $i=1,\dots,\spacedim$. The divergence-free property is convenient since the volume integral vanishes after integrating by parts  $\plbr{\nabla \trialfcn_i,\curl \plbr{ \qpoly }}_\el$.

\subsubsection*{V4: Enlarged enhanced space with elliptic projection}

In the fourth formulation, the enhanced space defined for $\mathrm{V1}$ is used, i.e. $\VspacedenhV{4} \plbr{\el} = \VspacedenhV{1} \plbr{\el}$.
Therefore, the following high-order projection is computed:
\begin{equation}
\begin{aligned}
    \begin{cases}
        \begin{aligned}
            &\bila^\el \plbr{\PinablaS{4} \trialfcn_i,\qpoly} = \bila^\el \plbr{\trialfcn_i,\qpoly} && \forall \qpoly \in \Pspace_{\kordaug} \plbr{\el},\\
            &P_0^\el \plbr{ \PinablaS{4} \trialfcn_i, \qpoly } = P_0^\el \plbr{ \trialfcn_i, \qpoly } && \forall \qpoly \in \Pspace_0 \plbr{\el},
        \end{aligned}
\end{cases}
\end{aligned}
\label{eq:SS1_projections}
\end{equation}
for $i=1,\dots,\spacedim$, where $P_0^\el$ is defined in \eq{laplace_invcondition2}. 

\subsubsection*{V5: Augmented space with elliptic projection}

This version uses the same space as $\mathrm{V2}$, i.e. $\VspacedenhV{5} \plbr{\el} = \VspacedenhV{2} \plbr{\el}$.
Therefore, it has an increased number of internal degrees of freedom.

The following projection is then computed by using the available information:
\begin{equation}
\begin{aligned}  
    \begin{cases}
        \begin{aligned}
            &\bila^\el \plbr{\PinablaS{5} \trialfcn_i,\qpoly} = \bila^\el \plbr{\trialfcn_i,\qpoly} && \forall \qpoly \in \Pspace_{\kordaug} \plbr{\el},\\
            &P_0^\el \plbr{\PinablaS{5} \trialfcn_i, \qpoly } = P_0^\el \plbr{ \trialfcn_i, \qpoly } && \forall \qpoly \in \Pspace_0 \plbr{\el}.
        \end{aligned}
    \end{cases}
\end{aligned}
\label{eq:SS2_projections}
\end{equation}
for $i=1,\dots,\spacedimaug$.

\subsubsection*{V6: Divergence-free space with $L^2$ projection}

In this case, the VEM space is the standard one, defined in \eq{VEMlocalspace_enhanced}, while the polynomial projection is computed in the space: $$\nabla \plbr{\Pspace_{\kord} \plbr{\el}} \oplus \curl \plbr{\Pspace_{\kordaug} \plbr{\el} \setminus \Pspace_{\kord} \plbr{\el}},$$ as:
\begin{equation}
\begin{aligned}      
    \begin{cases}
        \begin{aligned}
            &\plbr{\PinablaStwo{6} \nabla \trialfcn_i,\nabla\qpoly}_\el = \plbr{\nabla \trialfcn_i,\nabla\qpoly}_\el  && \forall \qpoly \in \Pspace_{\kord} \plbr{\el},\\
            &\plbr{\PinablaStwo{6}\nabla \trialfcn_i,\curl \plbr{ \qpoly }}_\el = \plbr{\nabla \trialfcn_i,\curl \plbr{ \qpoly }}_\el && \forall \qpoly \in \Pspace_{\kordaug} \plbr{\el} \setminus \Pspace_{\kord} \plbr{\el},\\
            &P_0^\el \plbr{ \PinablaStwo{6} \trialfcn_i, \qpoly } = P_0^\el \plbr{ \trialfcn_i, \qpoly } && \forall \qpoly \in \Pspace_0 \plbr{\el},
        \end{aligned}
    \end{cases}
\end{aligned}
\label{eq:SS3_projections}
\end{equation}
for $i=1,\dots,\spacedim$. This formulation refers to an $L^2$ projection. However, the operator $P_0^\el$ is needed since the space $\nabla \plbr{\Pspace_{\kord} \plbr{\el}}$ vanishes for constant polynomials.

\subsubsection{Algorithm for computing $\laug$ on a given element $\el$}

The choice of $\laug$ is crucial for the definition of the proposed schemes. However, a general theory for the efficient choice of $\laug$ is still lacking; thus, an iterative procedure is carried out to determine its value~\cite{berrone2025stabilization}. The following procedure applies to all the aforementioned self-stabilized formulations. The following necessary condition must be satisfied~\cite{berrone2025stabilization}:
\begin{equation}
    \dim \plbr{\Pspace_{\kordaug} \plbr{\el}} \geq \spacedim - 1.
    \label{eq:necessary_condition_selfstab}
\end{equation}

The inequality of \eq{necessary_condition_selfstab} is necessary but not sufficient. The order $\kordaug$ will be chosen so that the local stiffness matrix $\KE$ is full rank. 
The order $\laug$ can be determined by applying the procedure outlined in \alg{polyorder_selfstab}, which differs from the one introduced in~\cite{berrone2025lowest}, in which a QR decomposition is employed. Notice that \alg{polyorder_selfstab} is preferred solely for simplicity of implementation.

\begin{algorithm}
\begin{algorithmic}
    \State{Let $\laug$ be the smallest number satisfying \eq{necessary_condition_selfstab}}
    \State{compute local stiffness matrix $\KE$ }
    \State{evaluate $r_K = \mathtt{rank} \plbr{\KE,\mathtt{tol}}$}
    \While{$r_K < \spacedim - 1$}
    \State{$\laug = \laug + 1$}
    \State{update local stiffness matrix $\KE$ }
    \State{evaluate $r_K = \mathtt{rank} \plbr{\KE,\mathtt{tol}}$}
    \EndWhile
\end{algorithmic}
\caption{Algorithm for computing $\laug$ on a given element $\el$}
\label{alg:polyorder_selfstab}
\end{algorithm}
The tolerance $\mathtt{tol}$ is a user-defined parameter, whose effect on the solution is investigated in the numerical section.

\section{Linear Elasticity problem} \label{sec:LinearElasticity_SelfStab}

The stabilized and self-stabilized formulations presented for the Laplace problem can be extended to linear elasticity. Without loss of generality, the 2D problem is considered.
The two planar displacement components are defined as:
\begin{equation}
    \disp \plbr{x,y} = 
    \begin{bmatrix}
        u\plbr{x,y}\\
        v\plbr{x,y}
    \end{bmatrix},
    \label{eq:kinematics}
\end{equation}
while the strain-displacement relation, under the assumptions of infinitesimal displacements, reads:
\begin{equation}
\strain \plbr{\disp} =
\begin{bmatrix}
    \varepsilon_{xx}\\
    \varepsilon_{yy}\\
    \gamma_{xy}
\end{bmatrix}
=
\begin{bmatrix}
    u_{,x}\\
    v_{,y}\\
    u_{,y} + v_{,x}
\end{bmatrix},
\label{eq:strains}
\end{equation}
where $(\cdot)_{,x}$ and $(\cdot)_{,y}$ denote derivatives with respect to the coordinates $x$ and $y$.
The plane stress constitutive law is formulated as:
\begin{equation}
\begin{aligned}
\bm{\sigma} = 
\begin{bmatrix}
    \sigma_{xx}\\
    \sigma_{yy}\\
    \sigma_{xy}
\end{bmatrix}
&=
\begin{bmatrix}
    \frac{E_1 \plbr{x,y}} {1 - \nu_{12} \plbr{x,y} \nu_{21} \plbr{x,y}} & \frac{\nu_{12} \plbr{x,y} E_2\plbr{x,y}}{1 - \nu_{12} \plbr{x,y} \nu_{21} \plbr{x,y}} & 0\\
    \frac{\nu_{12} \plbr{x,y} E_2\plbr{x,y}}{1 - \nu_{12} \plbr{x,y} \nu_{21} \plbr{x,y}} & \frac{E_2\plbr{x,y} \plbr{x,y}}{1 - \nu_{12} \plbr{x,y} \nu_{21} \plbr{x,y}} & 0\\
    0 & 0 & G_{12}\plbr{x,y}
\end{bmatrix}
\begin{bmatrix}
    \varepsilon_{xx}\\
    \varepsilon_{yy}\\
    \gamma_{xy}   
\end{bmatrix}\\
&= \cost \plbr{x,y} \, \strain \plbr{\disp},
    \label{eq:const_law}
\end{aligned}    
\end{equation}
where $E_1 \plbr{x,y}$ and $E_2 \plbr{x,y}$ are the elastic moduli, $\nu_{12} \plbr{x,y}$ and $\nu_{21} \plbr{x,y}$ are the Poisson's coefficients, and $G_{12} \plbr{x,y}$ is the in-plane tangential modulus. The properties are a function of the position in view of the future extension to the linear elasticity problem with variable coefficients. 

The weak form of the linear elasticity problem reads: 
\begin{equation}
    \begin{cases}
        \text{find } \disp \in \sqbr{H_0^1 \plbr{\surom}}^2 \text{ such that:}\\
        \bila \plbr{\testdisp,\disp} = \plbr{\testdisp,\forcevec}_\surom \quad \forall \testdisp \sqbr{H_0^1 \plbr{\surom}}^2.\\
    \end{cases}
\label{eq:linearelasticity_governingeqs}
\end{equation}
In the above equation, $\bila \plbr{\testdisp,\disp}$ is the internal virtual work, and $\plbr{\testdisp,\forcevec}$ is the external virtual work due to the body force $\forcevec$, which are defined as:
\begin{equation}
    \begin{aligned}
        &\bila \plbr{\testdisp,\disp} = \int_\surom \strain \plbr{\testdisp}^T \, \cost \plbr{x,y} \,\strain \plbr{\disp} \, \de \surom, \quad
        &\plbr{\testdisp,\forcevec}_\surom = \int_\surom \testdisp \cdot \forcevec \, \de \surom.
    \end{aligned}
\label{eq:linearelasticity_bilabilf}
\end{equation}
To introduce the standard VEM formulation, $\cost$ is assumed to be constant within the elements in this initial section. The extension to $\cost = \cost\plbr{x,y}$ is discussed later, in Section \ref{sec:VC}.

The local virtual element space, with order of accuracy $\kord \geq 1$ is defined for $p \leq m$, as \cite{mengolini2019engineering}:
\begin{equation}
        \begin{aligned}\Vspacedvec \plbr{\el} = \left\{ \right.&\dispd \in \sqbr{H^1 \plbr{\el} \cap C^0 \plbr{\el}}^2  :  \bm{L} \sqbr{ \bm{\sigma} \plbr{\dispd} } \big|_\el \in \sqbr{ \Pspace_{\kord-2}\plbr{\el} }^2 ,\\ 
         & \dispd \big|_{\bound} \in \sqbr{ \Pspace_{\kord}\plbr{\bound} }^2 \forall e=1,\dots,\nedges, \\
         &\dispd \big|_{\boundcurv} \in \sqbr{ \Tilde{\Pspace}_{\kord}\plbr{\boundcurv} }^2 \forall e=1,\dots,\nedgescurv,\\
         &\plbr{\dispd,\qpolyvec}_\el = \plbr{\PiK \dispd,\qpolyvec}_\el \; \forall \qpolyvec  \in \sqbr{ \Pspace_{\kord}\plbr{\el} }^2 \setminus\sqbr{ \Pspace_{\kord-2}\plbr{\el} }^2 \left. \right\},\end{aligned}
    \label{eq:VEMspaceLinearElasticity_local}
\end{equation}
where the operator $\bm{L} \sqbr{\cdot}$ is the differential matrix defined as:
\begin{equation}
    \bm{L} \sqbr{\cdot}= 
    \begin{bmatrix}
    ,x &  0 & ,y \\ 
    0  & ,y & ,x 
    \end{bmatrix}.
    \label{eq:Ldiff_operator}
\end{equation}
The space has dimension $\spacedim$ and the degrees of freedom for $\dispd \in \Vspacedvec \plbr{\el}$ are defined as:
\begin{itemize}
    \item $\textbf{DOF}_{1}$: the values at the $\nv$ vertices of $\el$,
    \item $\textbf{DOF}_{2}$: for $\kord > 1$, the values of at the $\kord-1$ internal Gauss-Lobatto quadrature points on each $\bound$ and the values at the $\kord-1$ points on each $\boundcurv$ that are images through $\gamma_e$ of the $\kord-1$ internal Gauss-Lobatto quadrature points on $I_e$,
    \item $\textbf{DOF}_{3}$: the internal moments:
    \begin{equation}
        \frac{1}{\shbr{\el}} \int_\el \qpolyvec \plbr{x,y} \cdot \dispd\plbr{x,y} \de \el \quad \forall \qpolyvec \in \sqbr{ \Qpolyspace_{\kord-2}\plbr{\el} }^2.\\
        \label{eq:localmoments}
    \end{equation}
\end{itemize}

Following the previous definitions, it is now possible to define the global virtual element space as:
\begin{equation}
    \Vspacedvec = \cubr{ \dispd \in \sqbr{H_0^1\plbr{\Omega}}^2 : \dispd{\big|_\el} \in \Vspacedvec \plbr{\el} \quad \forall \el \in \globdom }.
    \label{eq:VEMspaceLinearElasticity_global}
\end{equation}

The generic function $\bm{v}_h \in \Vspacedvec \plbr{\el}$ is expressed as a linear combination of the Lagrangian basis functions $\trialfcnvec_i$ as in \eq{lagrangeinterpolation}. The discrete counterpart of the problem defined by \eq{linearelasticity_governingeqs} is presented, in the stabilized VEM framework, similarly to the Laplace problem. Thus, the following $\bila^\el$-orthogonal polynomial projection is considered:
\begin{equation}
    \begin{cases}
    \begin{aligned}
        &{\bila}^\el \plbr{\PiK \trialfcnvec_i,\qpolyvec} = {\bila}^\el \plbr{\trialfcnvec_i,\qpolyvec} && \forall \qpolyvec \in \sqbr{ \Pspace_{\kord}\plbr{\el} }^2,\\
        &P_0^\el \plbr{\PiK \trialfcnvec_i, \qpolyvec } = P_0^\el \plbr{\trialfcnvec_i, \qpolyvec } && \forall \qpolyvec \in RM\plbr{\el},
    \end{aligned}    
    \end{cases}
    \label{eq:energyproj_linearelasticity}
\end{equation}
for all $i=1,\dots,\spacedim$, with:
\begin{equation}
    \begin{aligned}
        P_0^\el \plbr{\trialfcnvec_i, \qpolyvec} = \frac{1}{\nv} \sum_{j=1}^{2\nv} \text{dof}_j \plbr{\trialfcnvec_i} \text{dof}_j \plbr{\qpolyvec},
    \end{aligned}
    \label{eq:constenergyproj_linearelasticity}
\end{equation}
and $RM\plbr{\el}$ being the set of rigid body motions defined as
\begin{equation}
    RM\plbr{\el} = \cubr{\mathbf{v}\in\sqbr{H^1\plbr{\el}}^2:\,\strain \plbr{\mathbf{v}}=0} = \cubr{ \sqbr{c_1,c_2}^T + c_3\sqbr{-y,x}^T }.
\end{equation}

where $c_i$ are constants representing the planar rigid body motions. The stabilization forms listed in Section \ref{subsec:Laplace_stabilization} are still valid options for the problem under consideration. In this case, $\KEs_1$, $\KEs_2$, $\KEs_4$ and $\KEs_5$ should be multiplied by the trace of the consistency matrix \cite{mengolini2019engineering}. 

The procedure to derive the self-stabilized formulations mirrors the one used for the Laplace problem, after suitable modification of the discrete spaces. In particular, the projection operators for $\mathrm{V1}$ and $\mathrm{V2}$ are obtained by replacing the gradient with the strain operator $\strain$ in \eqs{SS4_projections}{SS5_projections}, respectively. The formulations $\mathrm{V4}$ and $\mathrm{V5}$ require an adapted definition of the bilinear form.
The formulation $\mathrm{V3}$ is defined in the same space of $\mathrm{V1}$ and the associated polynomial projection is now computed in $$\sqbr{\Pspace_{\kord-1} \plbr{\el}}^3 \oplus \strain^\curlcross \plbr{\sqbr{\Pspace_{\kordaug} \plbr{\el} \setminus \Pspace_{\kord} \plbr{\el}}^2},$$ where the $\curl$ operator was replaced by:
\begin{equation}
    \strain^\curlcross \plbr{\cdot}=
    \begin{bmatrix}
        ,y  & 0   \\
        0   & -,x \\
        -,x & ,y  
    \end{bmatrix}.
    \label{eq:espilon_operator_divfree}
\end{equation}
Lastly, the formulation $\mathrm{V6}$ uses the same space of $\mathrm{V4}$ with polynomial projection onto the space: $$\strain \plbr{\sqbr{\Pspace_{\kord} \plbr{\el}}^2} \oplus \strain^\curlcross \plbr{\sqbr{\Pspace_{\kordaug} \plbr{\el} \setminus \Pspace_{\kord} \plbr{\el}}^2}.$$

\section{Stokes problem} \label{sec:Stokes_SelfStab}
The Stokes problem is a mixed problem with velocity $\vel$ and pressure $\pres$ as unknowns, whose weak form reads: 
\begin{equation}
    \begin{cases}
    \begin{aligned}
        &\text{find } \plbr{\vel,\pres} \in \sqbr{H_0^1 \plbr{\surom}}^2\times L_0^2 \plbr{\surom} \text{ such that:}\\
        &\bila \plbr{\testvel,\vel} + \bilb \plbr{\testvel,\pres} = \plbr{\testvel,\forcevec}_\surom && \forall \testvel \in \sqbr{H_0^1 \plbr{\surom}}^2,\\
        &\bilb \plbr{\testpres,\vel}= 0 && \forall \testpres \in L_0^2 \plbr{\surom},
    \end{aligned}    
    \end{cases}
\label{eq:governingeqs_Stokes}
\end{equation}
with:
\begin{equation}
    \begin{aligned}
        & \bila \plbr{\testvel,\vel} = \int_\surom \nu \bm{\nabla} \vel  \bm{:} \bm{\nabla} \testvel \; \de \surom, \quad
        & \bilb \plbr{\testvel,\pres} = \int_\surom \diverg \plbr{\testvel} \pres \; \de \surom,\\
        & \plbr{\testvel,\forcevec}_\surom = \int_\surom \testvel \cdot \forcevec \; \de \surom,
    \end{aligned}
    \label{eq:bilinearforms_Stokes}
\end{equation}
where $\nu$ is the viscosity and $L_0^2 \plbr{\surom}$ is the subspace of $L^2 \plbr{\surom}$ of null-mean valued functions.
In this work, the divergence-free approximation detailed in \cite{daveiga2017divergence,dassi2020bricks} is considered. Therefore, for $\kord \geq 2$ and $\kord \leq m$, the local velocity space is defined as:
\begin{equation}
    \begin{aligned}
        \Vspacedvec \plbr{\el} = \left\{ \right. \veld \in &\sqbr{H^1 \plbr{\el} \cap C^0 \plbr{\el}}^2 : \diverg \plbr{\veld} \in \Pspace_{\kord-1} \plbr{\el},\\
        &- \nu \bm{\Delta} \veld - \nabla s \in \sqbr{ \Gspace_{\kord-2}^{\perp} \plbr{\el} }^2  \text{for some } s \in L^2 \plbr{\el},\\
        & \veld \big|_{\bound} \in \sqbr{\Pspace_k \plbr{\bound}}^2 \forall e=1,\dots,\nedges,\\
        & \veld \big|_{\boundcurv} \in \sqbr{\Tilde{\Pspace}_\kord \plbr{\boundcurv}}^2 \forall e=1,\dots,\nedgescurv,\\
        &\left.\plbr{\veld,\qpolyvec}_\el=\plbr{\Pinablavec\veld,\qpolyvec}_\el \forall \qpolyvec \in \sqbr{ \Pspace_{\kord} \plbr{\el} }^2 \setminus \sqbr{ \Pspace_{\kord-2} \plbr{\el} }^2\right\},
    \end{aligned}
    \label{eq:vspace_stokes}
\end{equation}
while the local pressure space is:
\begin{equation}
    \Qspaced \plbr{\el} = \Pspace_{\kord-1} \plbr{\el}.
    \label{eq:pspace_stokes}
\end{equation}
More precisely $\sqbr{ \Gspace_{\kord} \plbr{\el} }^2 = \nabla \plbr{\Pspace_\kord\plbr{\el}}$ and $\sqbr{ \Gspace_{\kord}^{\perp} \plbr{\el} }^2$ is its $L^2$-orthogonal complement.

The degrees of freedom for the velocity $\veld \in \Vspacedvec$ are defined as:
    \begin{itemize}
        \item $\textbf{DOF}_{1}$: the values at the $\nv$ vertices of $\el$,
        \item $\textbf{DOF}_{2}$: the values at the $\kord-1$ internal Gauss-Lobatto quadrature points on each $\bound$ and the values at the $\kord-1$ points on each $\boundcurv$ that are images through $\gamma_e$ of the $\kord-1$ internal Gauss-Lobatto quadrature points on $I_e$,
        \item $\textbf{DOF}_{3}$: the internal moments:
        \begin{equation}
            \frac{1}{\shbr{\el}} \int_\el \veld \plbr{x,y} \cdot \gpolyperp \plbr{x,y} \de \el \quad \forall \gpolyperp \in \sqbr{ \Gspace_{\kord-2}^{\perp} \plbr{\el} }^2,
            \label{eq:vinternalmomentsperp_stokes}
        \end{equation}
        \item $\textbf{DOF}_{4}$: the divergence moments up to order $\kord-1$:
        \begin{equation}
            \frac{1}{\shbr{\el}} \int_\el \diverg \plbr{\veld \plbr{x,y}} \qpoly \plbr{x,y} \de \el \quad \forall \qpoly \in \Qpolyspace_{\kord-1} \plbr{\el} \setminus \Qpolyspace_0 \plbr{\el}.
            \label{eq:vinternalmomentsscalar_stokes}
        \end{equation}
    \end{itemize}
For the pressure $\presd \in \Qspaced \plbr{\el}$ the degrees of freedom are defined as:   
    \begin{itemize}
    \item $\textbf{DOF}_{5}$: the internal moments up to order $\kord-1$:
        \begin{equation}
            \frac{1}{\shbr{\el}} \int_\el \presd \plbr{x,y} \qpoly \plbr{x,y} \de \el \quad \forall \qpoly \in \Qpolyspace_{\kord-1} \plbr{\el}.
            \label{eq:pinternalmoments_stokes}
        \end{equation}
    \end{itemize}

And the global virtual element spaces are defined as:
\begin{equation}
    \begin{aligned}
        &\Vspacedvec = \cubr{\veld \in \sqbr{H_0^1\plbr{\Omega}}^2 : \veld{\big|_\el} \in \Vspacedvec \plbr{\el} \quad \forall \el \in \globdom },\\
        &\Qspaced = \cubr{ \presd \in L_0^2\plbr{\Omega} : \presd{\big|_\el} \in \Qspaced \plbr{\el} \quad \forall \el \in \globdom }.\\
    \end{aligned}
    \label{eq:VEMspaceStokes_global}
\end{equation}

The bilinear form $\bilb \plbr{\cdot,\cdot}$ is computable from the degrees of freedom. On the contrary, a projection is needed to obtain a computable bilinear form $\bilad \plbr{\cdot,\cdot}$. 
In case of a stabilized formulation, $\bilad \plbr{\cdot,\cdot}$ is constructed by using the vector counterpart of the operator $\Pinabla$, defined in \eq{VEMH1projection} for the Laplace problem, together with one of the stabilization forms listed in Section \ref{subsec:Laplace_stabilization}.
Analogously, the self-stabilized formulations are obtained by considering the vector counterpart of those for the Laplacian, see Section \ref{subsec:Laplace_selfstabilization}, after suitable modifications of the local velocity space.

\section{Numerical Results: Benchmarking stabilized and self-stabilized VEM} \label{sec:NumericalResults_SelfStab}

In this section,  the performance of different stabilization techniques and self-stabilized formulations is compared in terms of accuracy and conditioning. 

It is well known that the stabilization term may alter the results for high-order approximations, giving sub-optimal convergence~\cite{daveiga2016basic}. As the $\kord$-VEM is the focus of the present investigation, this section aims at clarifying the influence of the stabilization parameter and stabilization techniques on the accuracy and stability of the method. 

\begin{figure}[htbp]
    \centering
    \subfigure[Quadrilateral mesh.\label{fig:quad_mesh}]{
        \includegraphics[width=0.22\textwidth]{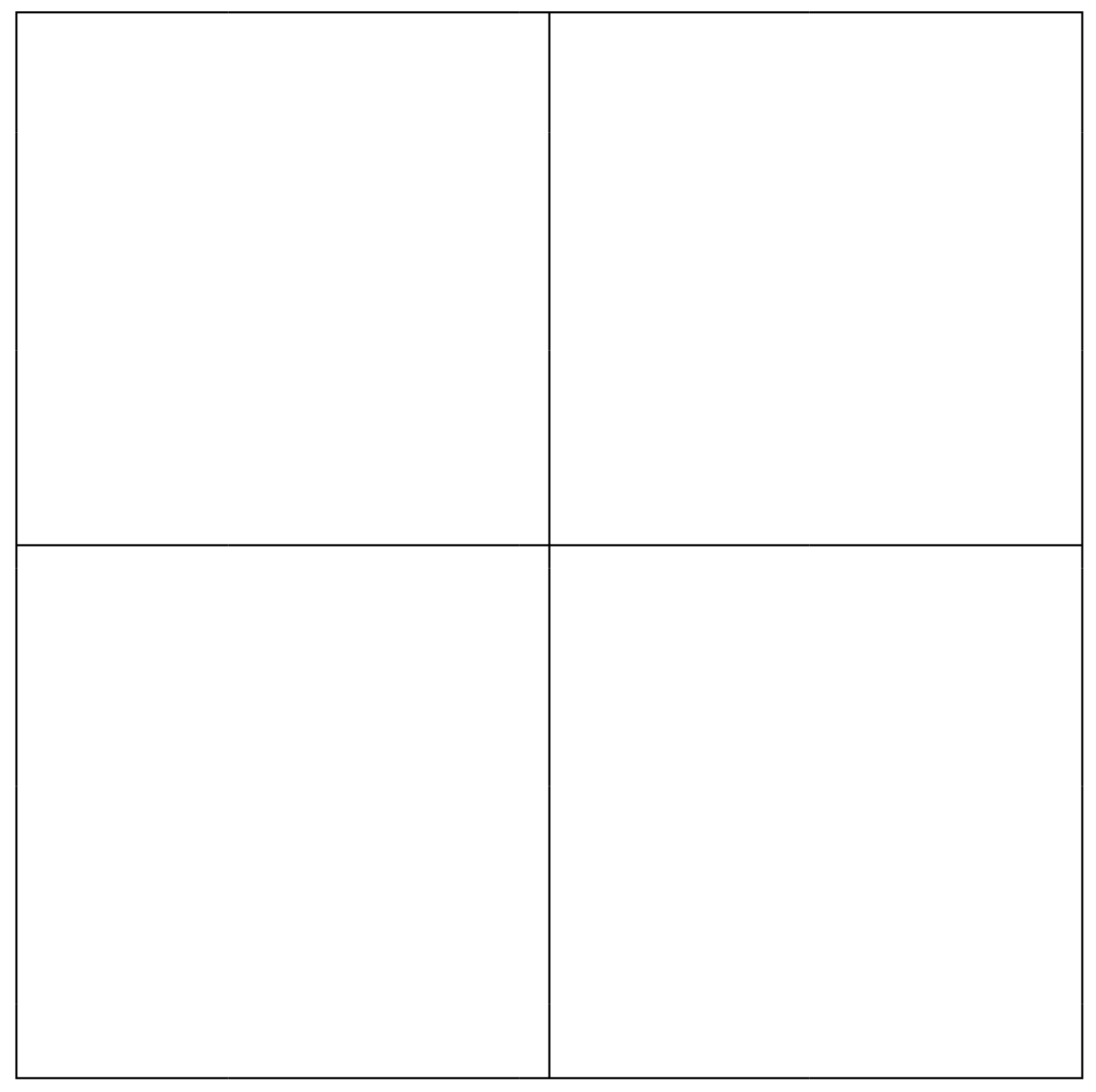}
    }
    \hfill
    \subfigure[Voronoi mesh.\label{fig:voronoi_mesh}]{
        \includegraphics[width=0.22\textwidth]{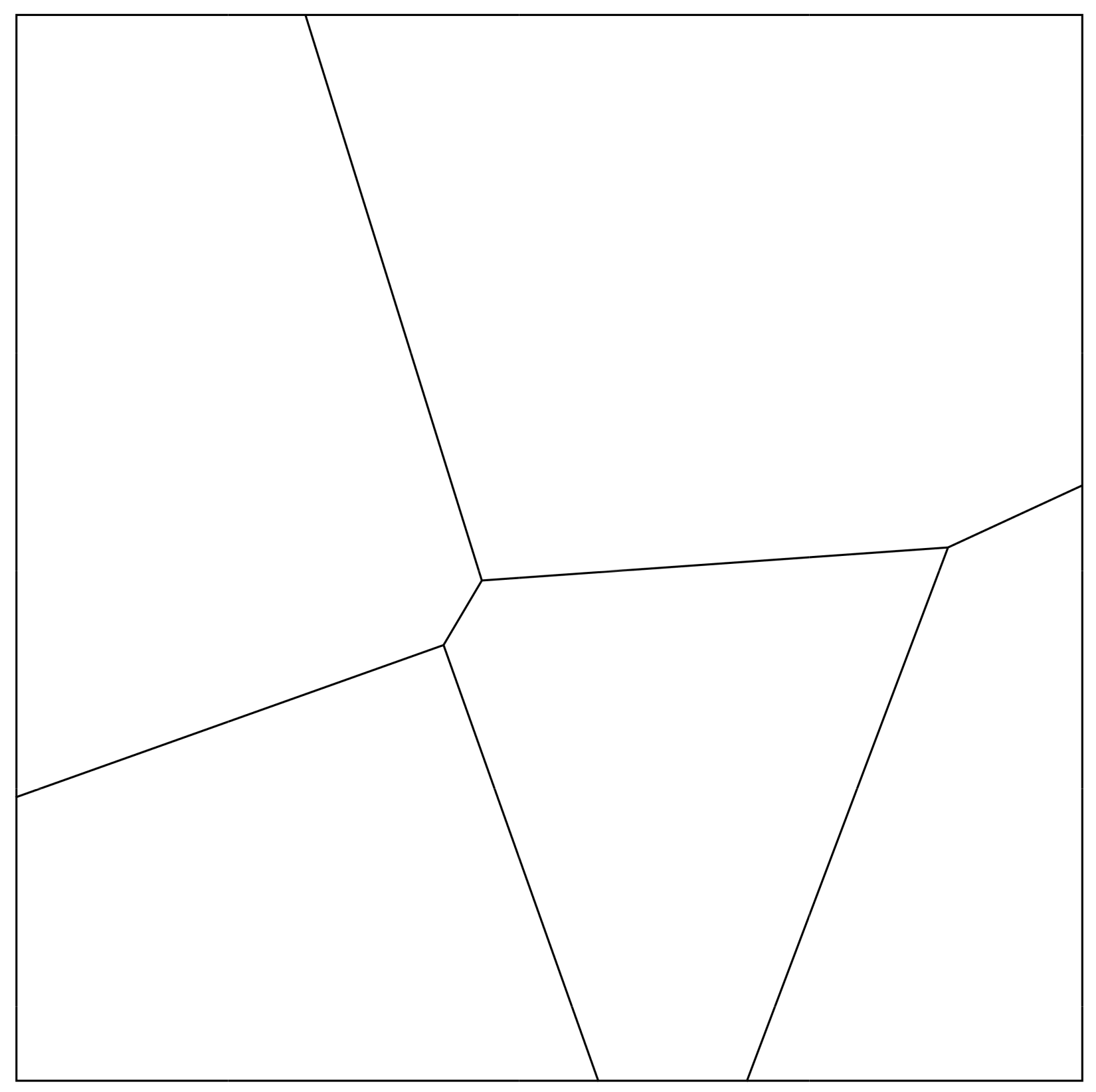}
    }
    \hfill
    \subfigure[Octagonal mesh.\label{fig:ottagoni_mesh}]{
        \includegraphics[width=0.22\textwidth]{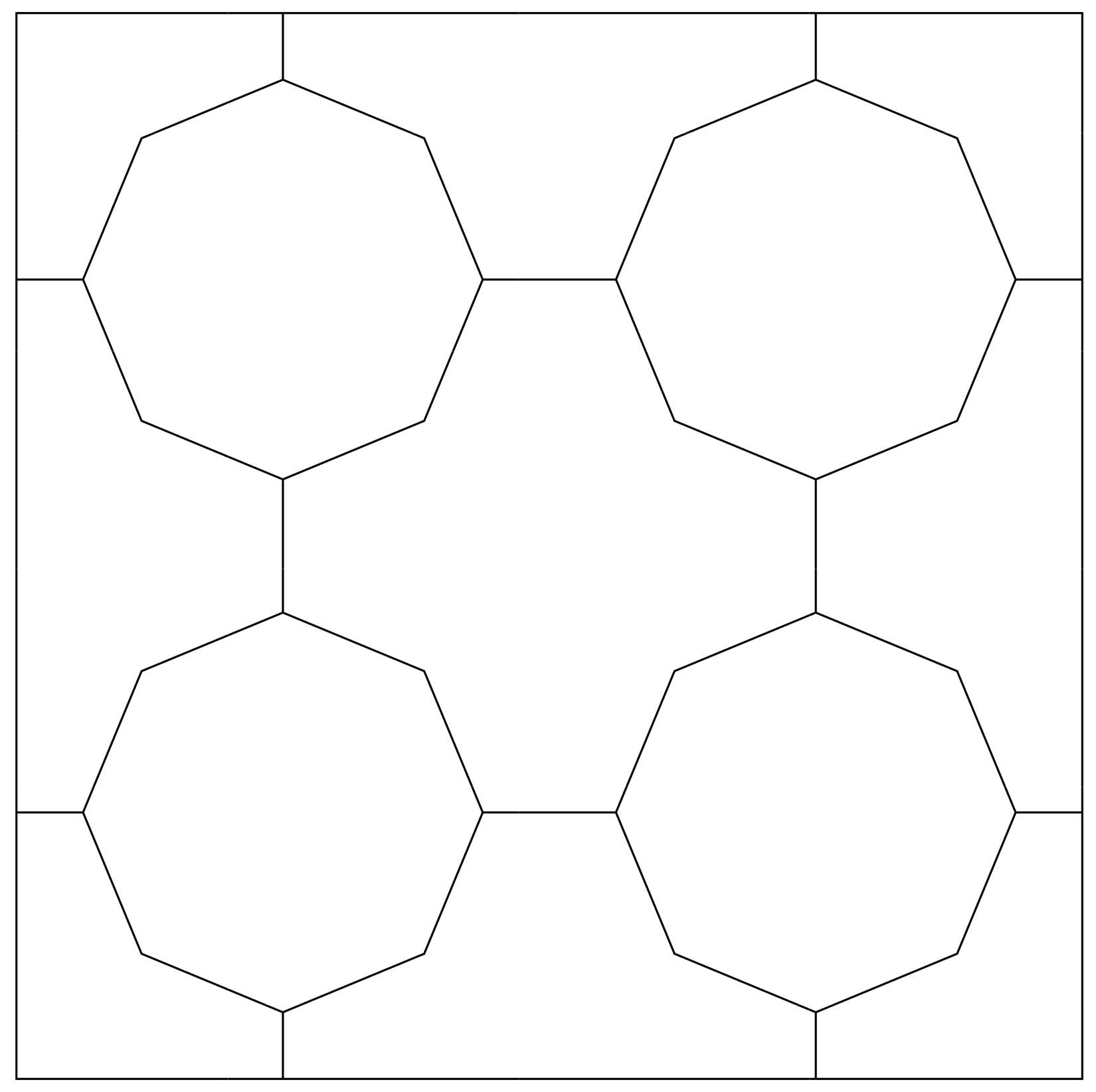}
    }
    \hfill
    \subfigure[B\'ezier-edge mesh.\label{fig:curv_mesh}]{
        \includegraphics[width=0.22\textwidth]{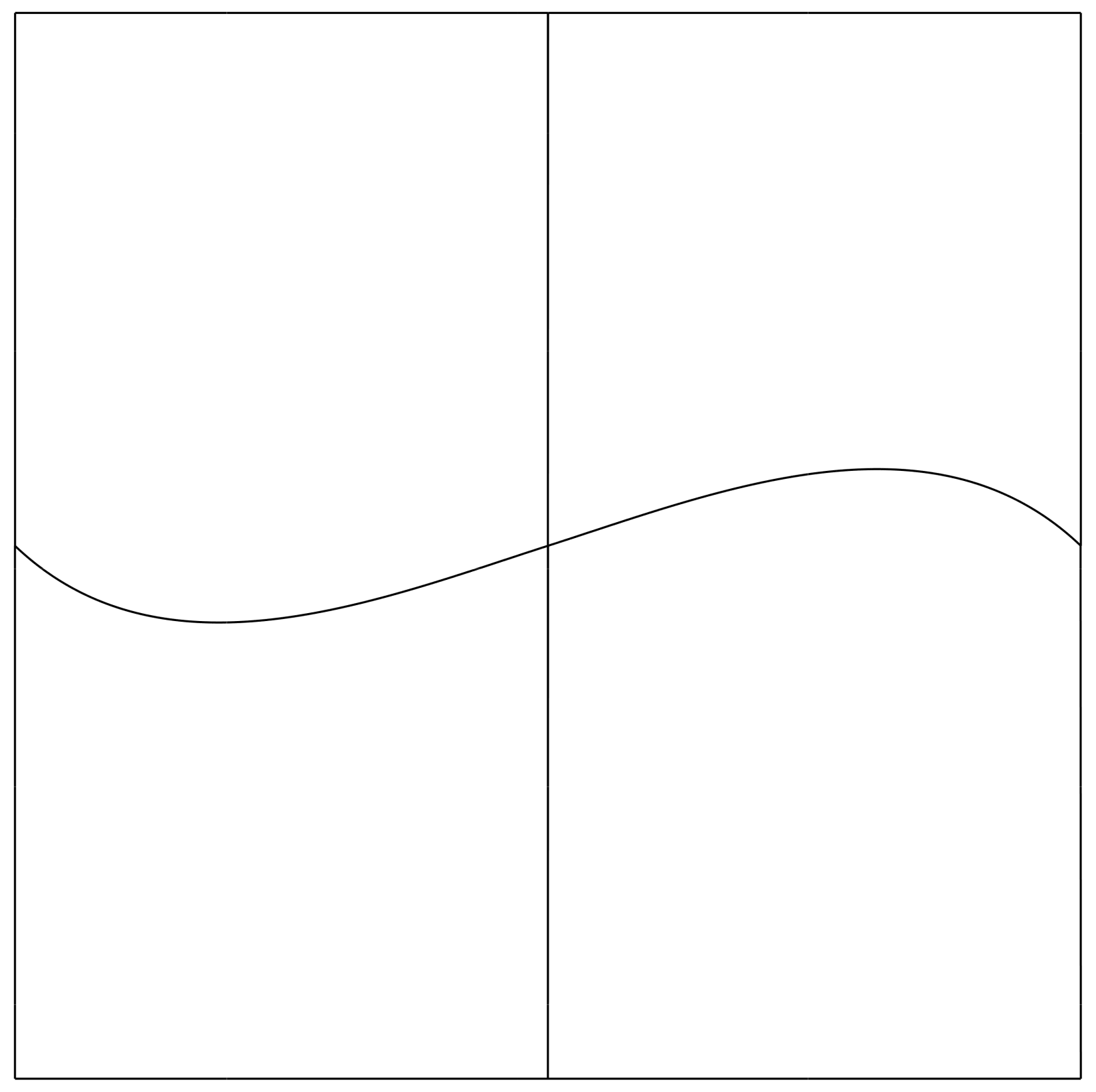}
    }
    \caption{Mesh types (Control points for B\'ezier-edge mesh in $\Omega = \plbr{0,1}^2$: $\bm{P}_0 = (0,0.5),\bm{P}_1 = (0.25,0.25),\bm{P}_2 = (0.75,0.75),\bm{P}_3 = (1,0.5)$).\label{fig:meshes}}
\end{figure}

Four different meshes are considered in this study. They are shown in \fig{meshes}. The first one is composed by four squares, the second by five Voronoi elements, the third by octagons and concave polygons, and the last one by four elements, each characterized by a third-order B\'ezier-type edge. 

A comprehensive study is conducted for the Laplace problem on the first three meshes, while for Stokes only the Voronoi one is considered. Regarding the elasticity problem, both the quadrilateral and the B\'ezier-edge mesh are employed, as the latter is more relevant in problems involving curvilinear reinforcements on plates.

For the polynomial basis $\Qpolyspace_\kord \plbr{\el}$,
the Modified Gram-Schmidt (MGS) orthonormalization \cite{mascotto2018ill} is used to obtain $L^2\plbr{\el}$ orthonormal polynomials. Indeed, they are more suitable in a $\kord$-VEM context than standard monomials, which suffer from numerical instabilities.

\subsection{Laplace problem} \label{subsec:Laplace_results_selfstab}

The load term and the boundary conditions are determined through an inverse approach, starting from the exact solution:
\begin{equation}
    \unku \plbr{x,y} = \sin\plbr{\pi x} \sin\plbr{\pi y},
    \label{eq:laplace_exactsolution}
\end{equation}
in $\surom = \plbr{0,1}^2$.
The error in the energy norm is computed as:
\begin{equation}
    e_{\nabla} =  \frac{ \plbr{\sum_{\el \in \globdom}  \norm{ \nabla \plbr{ \unku - \Piok \unkud } }_{0,\el}^2  }^{\frac{1}{2}} }{ \norm{\nabla \unku}_{0,\el} }.
    \label{eq:laplace_energyerror}
\end{equation}

The first part of this study focuses on stabilized VEM and aims at investigating the effect of $\tau$ on the accuracy and conditioning of the discrete problem. Then, for self-stabilized formulations, the effect of the rank estimation to tolerance variations in \alg{polyorder_selfstab} is analyzed. Finally, stabilized and self-stabilized formulations are compared in terms of accuracy and conditioning, by varying the order $\kord$ from 1 to 10.

\subsubsection{Parametric study on the stabilization parameter}

The first part of this study deals with the analysis of the stabilization parameter $\tau$ within the stabilized VEM. The most common stabilization $\KEs_1$ is used and the parameter $\tau$ is varied from $10^{-10}$ to $10^{10}$ with steps of $10^{2}$. Moreover, $\tau_{\mathrm{mean}} = \text{mean} \plbr{\text{eig} \plbr{\KEc}}$ is considered as well, where $\KEc$ is the consistency part of the stiffness matrix.

\begin{figure}[htbp]
        \centering
        \includegraphics[width=0.5\textwidth]{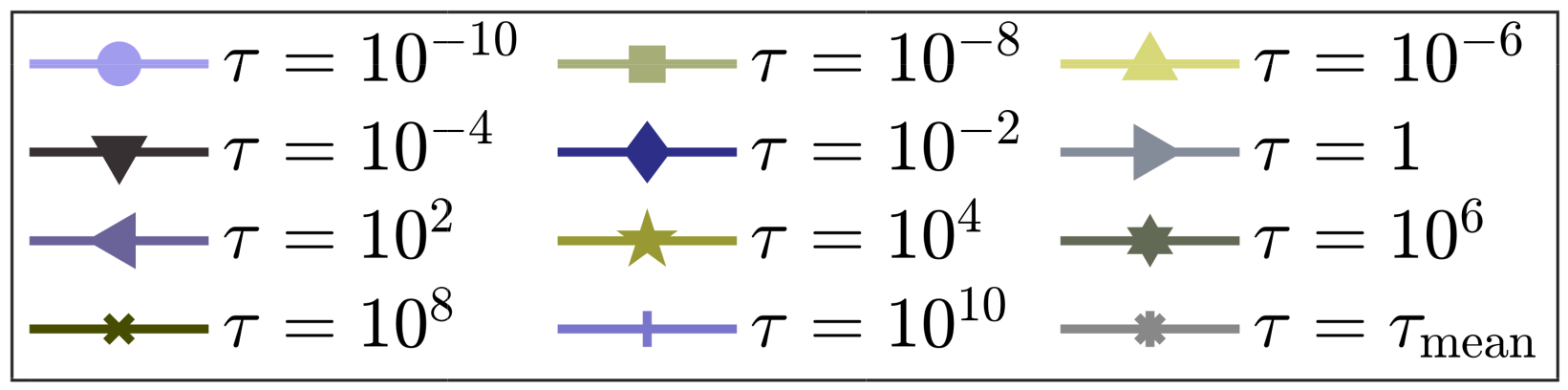}\\
	\subfigure[Energy norm error. \label{fig:tau_err_quad}]{
		\includegraphics[width=0.47\textwidth]{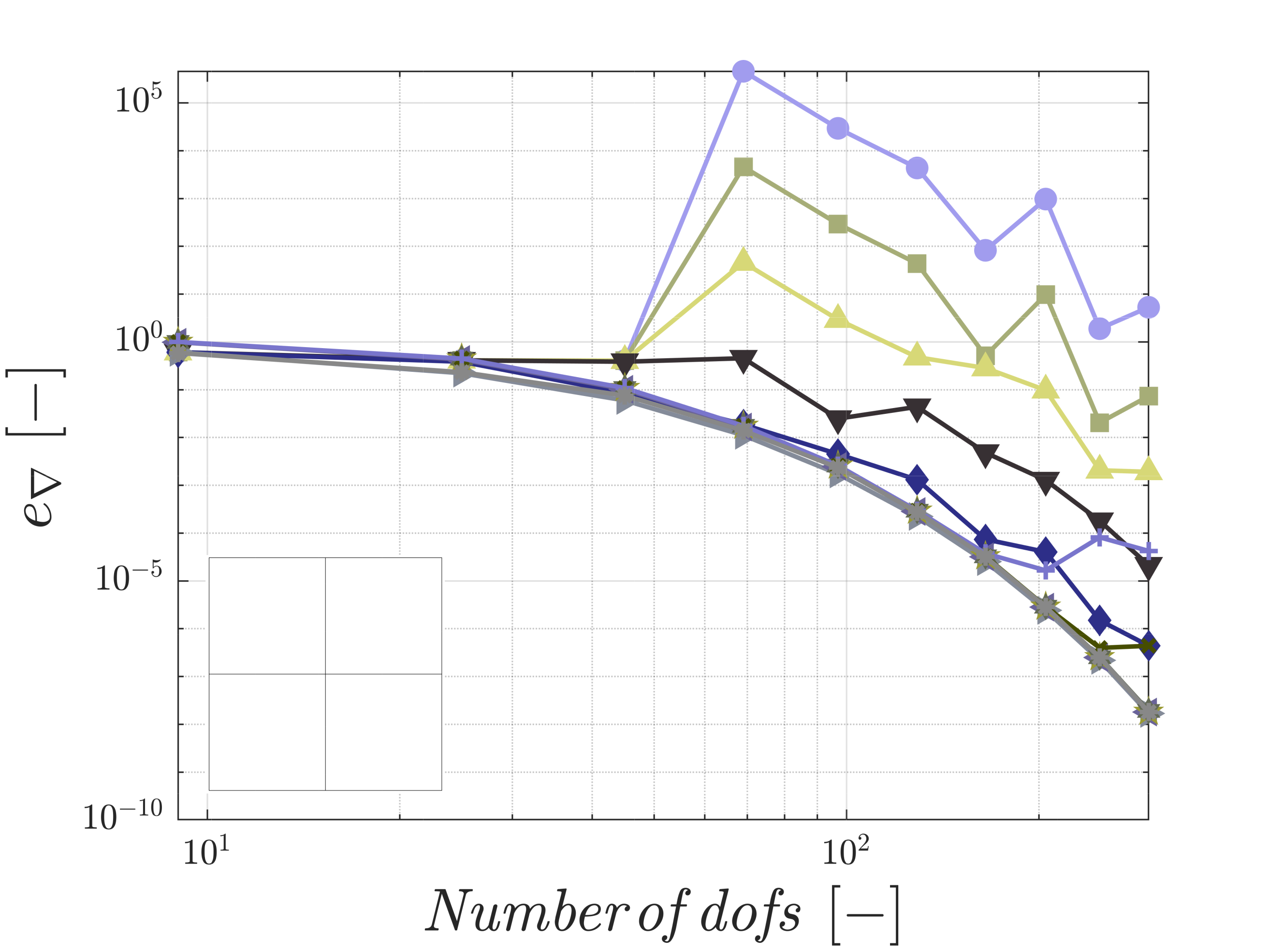}
	}
        \centering
        \hfill     
	\subfigure[Condition number. \label{fig:tau_cond_quad}]{
		\includegraphics[width=0.47\textwidth]{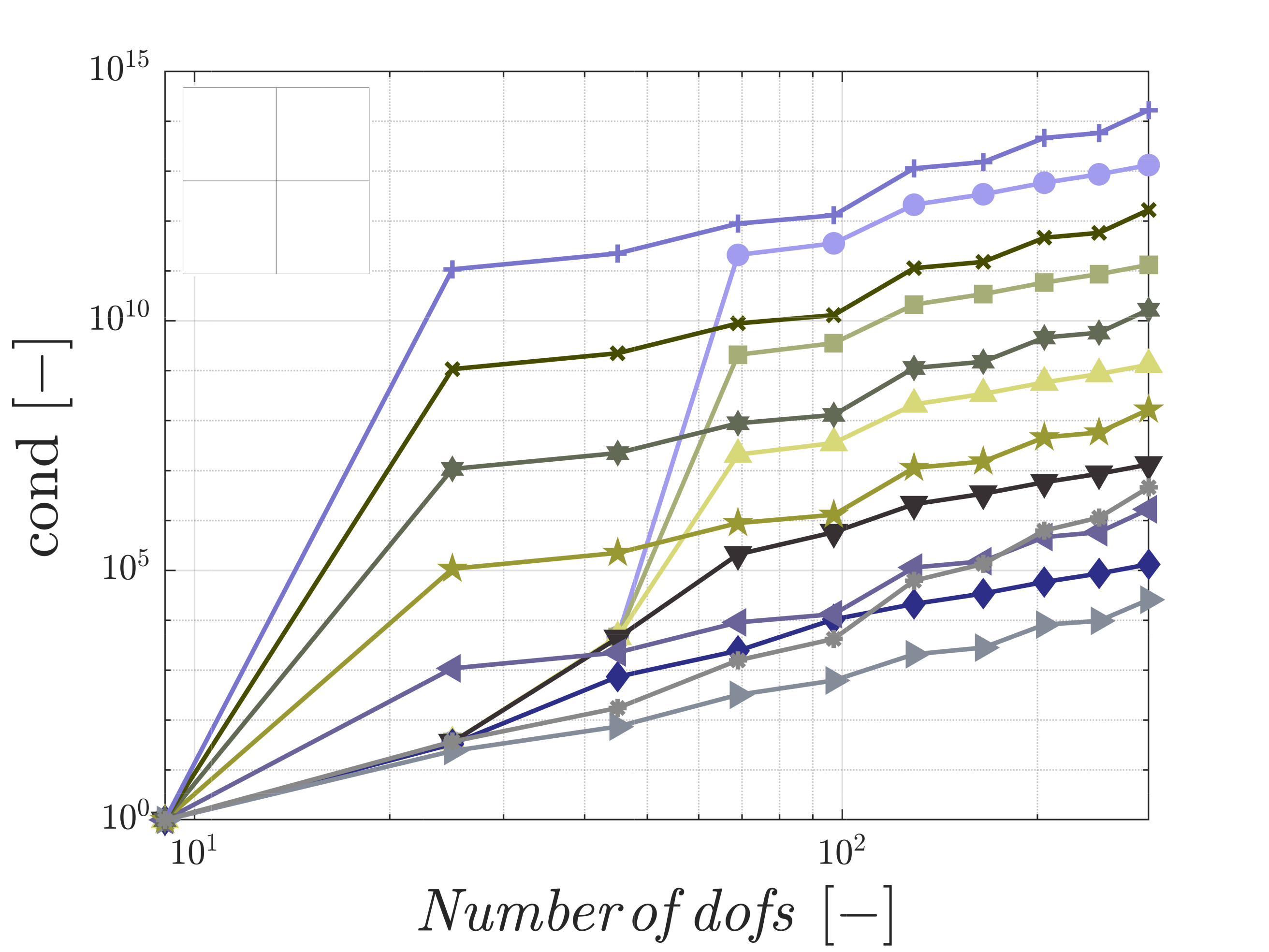}
	}     
    \caption{Laplace problem: parametric study on $\tau$ for quadrilateral mesh.\label{fig:tau_quad}}
 \end{figure}

The energy norm and the conditioning of the stiffness matrix are reported in \fig{tau_quad} for the quadrilateral mesh. 

The analysis of the error is conducted by inspection of  \fig{tau_err_quad}. No significant differences are noted up to order $\kord = 3$. So, the choice of the stabilization parameter $\tau$ has a limited impact for low-order approximations.
For order $\kord > 3$, the undesired behavior of the curves associated with $\tau \leq 10^{-2}$ is evident: the magnitude of the error makes the method completely unreliable. In contrast, a reduction of the error is achieved for $10^{-2}<\tau <10^{10}$. In this range of values the error is essentially the same irrespective of $\tau$. The errors observed for the largest values of $\tau$ are ascribed to worse conditioning of the matrix. This behavior is clearly illustrated in \fig{tau_cond_quad}. In particular, the condition number is shown to be high for very small and very large values. Optimal conditioning is observed for $\tau=1$, while for $\tau_{\mathrm{mean}}$ one can observe a condition number close to that of $\tau = 10^{2}$.

\begin{figure}[htbp]
        \centering
        \includegraphics[width=0.5\textwidth]{Images/patch_test_laplaciano_tau_legend.png}\\
	\subfigure[Energy norm error. \label{fig:tau_err_voronoi}]{
		\includegraphics[width=0.47\textwidth]{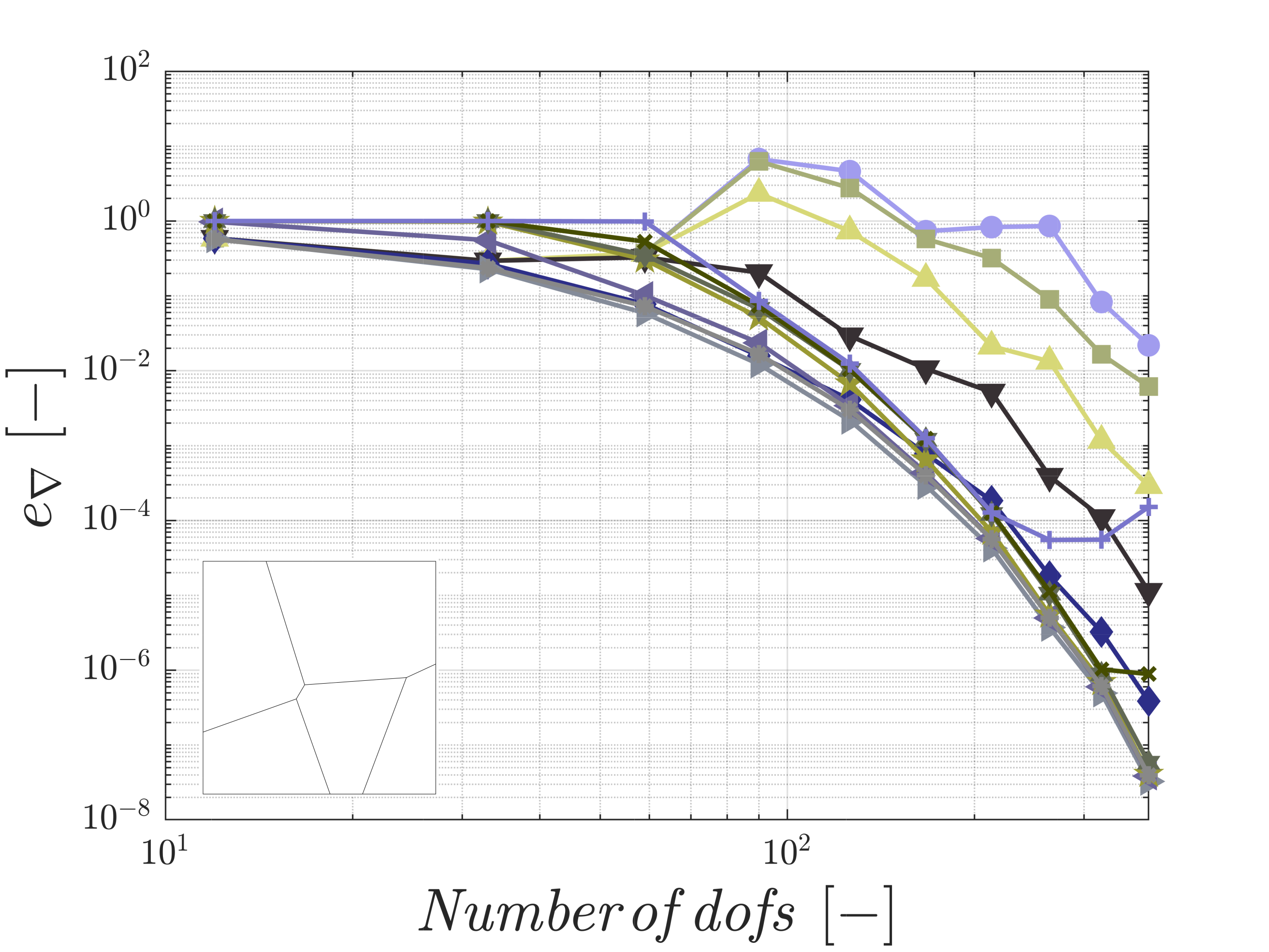}
	}
        \centering
        \hfill      
	\subfigure[Condition number. \label{fig:tau_cond_voronoi}]{
		\includegraphics[width=0.47\textwidth]{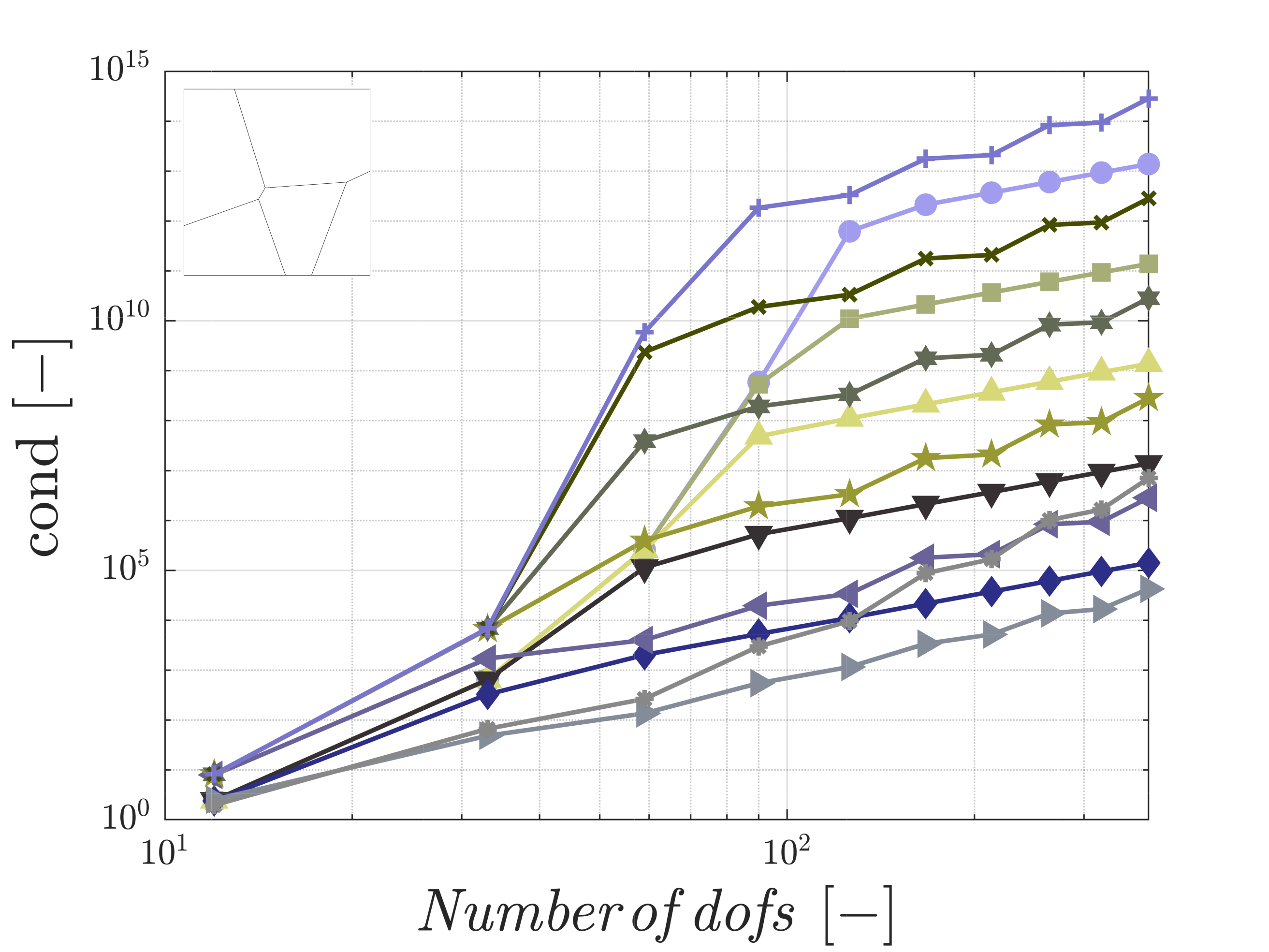}
	}     
    \caption{Laplace problem: parametric study on $\tau$ for Voronoi mesh.\label{fig:tau_voronoi}}
 \end{figure}

 The analysis is repeated in \fig{tau_voronoi}, now considering the Voronoi mesh. For the error, similar conclusions are drawn as in the case of quadrilateral mesh, i.e. lower $\tau$ values are associated to worse accuracy. One meaningful difference regards the trend of the corresponding curves, which are slightly smoother with respect to the quadrilateral case. A similar behavior is shown by the condition number reported in \fig{tau_cond_voronoi}.
 
\begin{figure}[htbp]
        \centering
        \includegraphics[width=0.5\textwidth]{Images/patch_test_laplaciano_tau_legend.png}\\
	\subfigure[Energy norm error. \label{fig:tau_err_ottagoni}]{
		\includegraphics[width=0.47\textwidth]{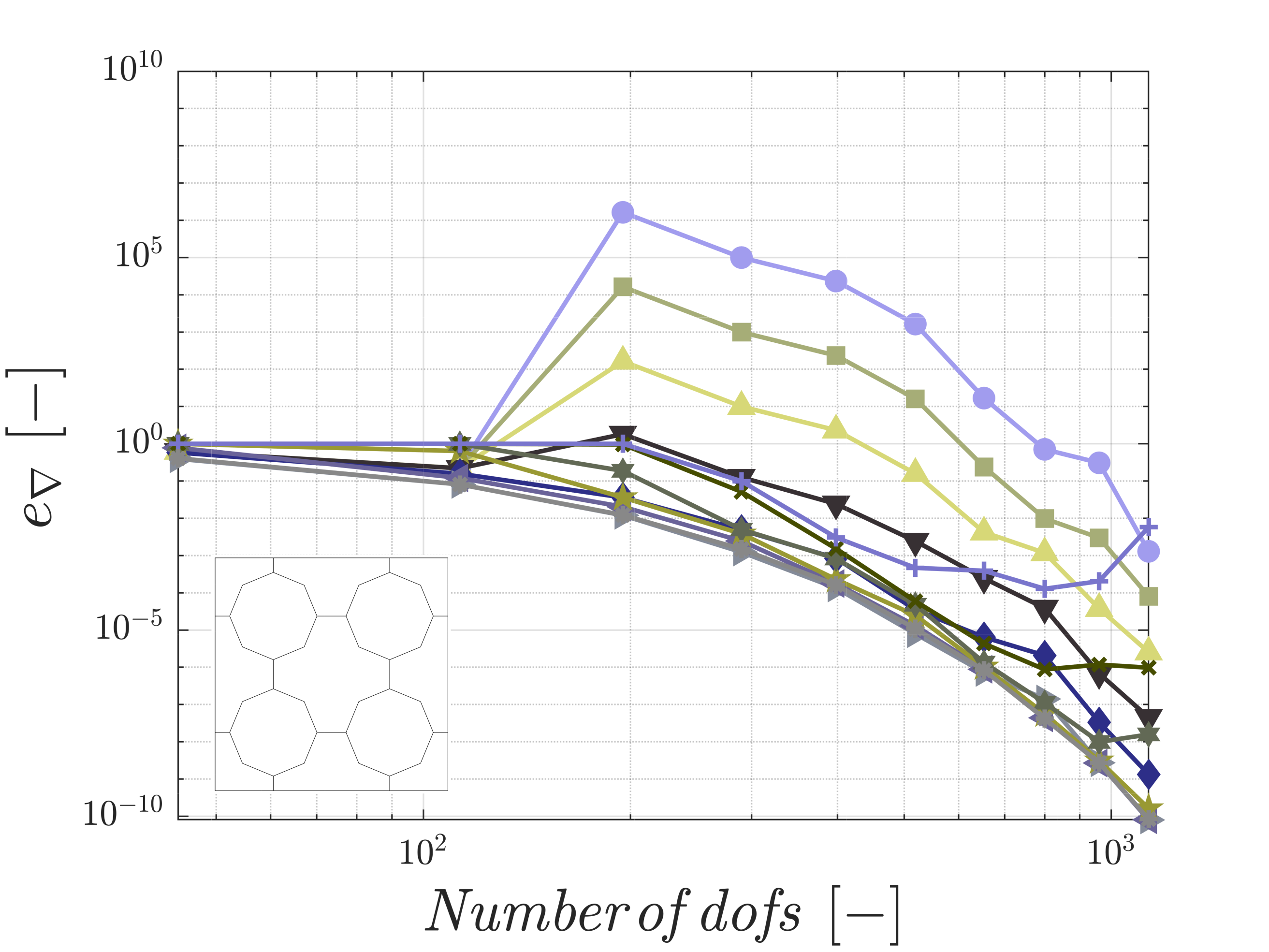}
	}
        \centering
        \hfill       
	\subfigure[Condition number. \label{fig:tau_cond_ottagoni}]{
		\includegraphics[width=0.47\textwidth]{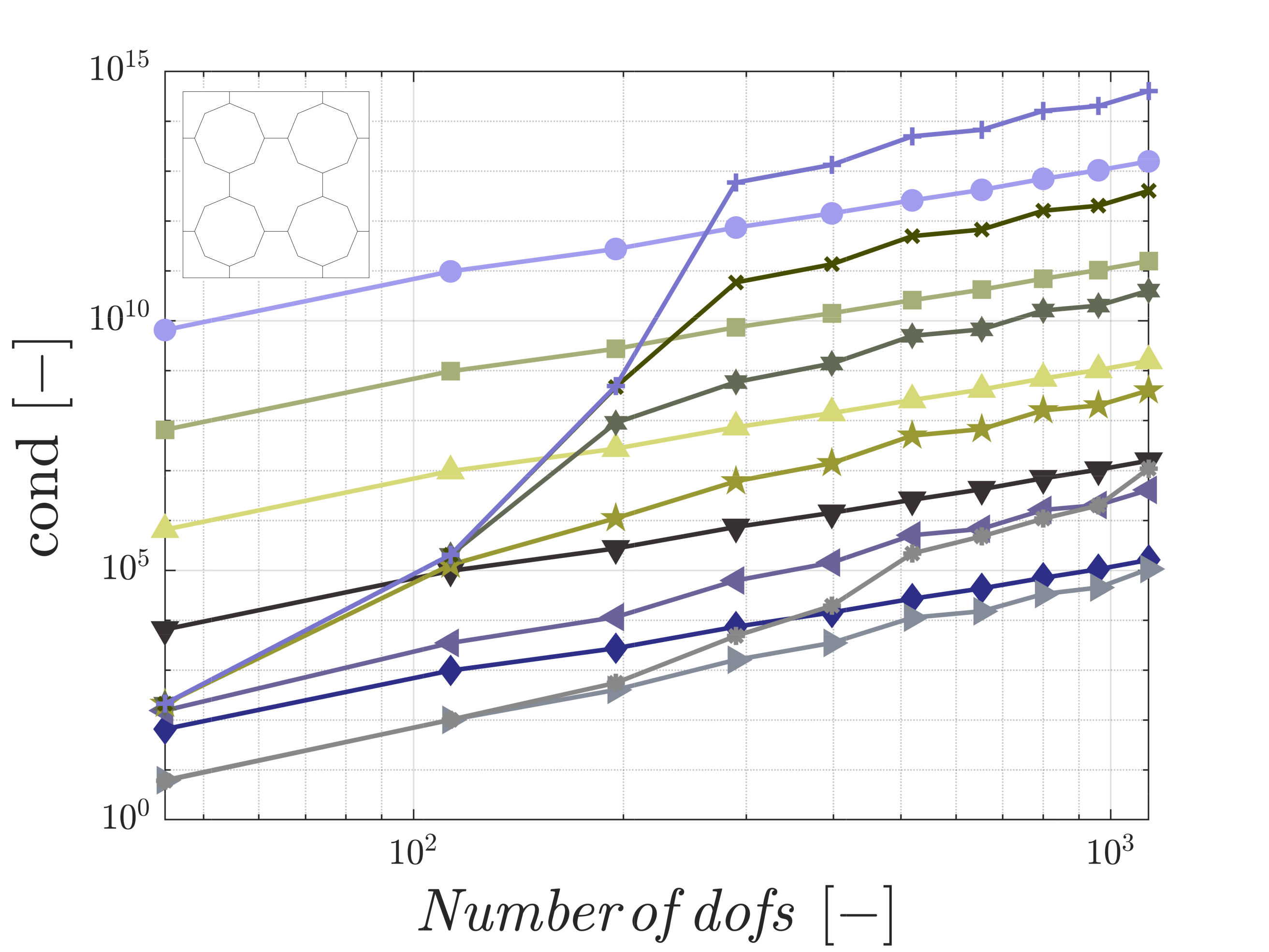}
	}     
    \caption{Laplace problem: parametric study on $\tau$ for octagonal mesh.\label{fig:tau_ottagoni}}
 \end{figure}

The third assessment regards the  octagonal mesh. The results are reported in \fig{tau_ottagoni}. Even in this case, the conclusions are similar to those reported for the previous cases. 

\subsubsection{Parametric study on the tolerance of the rank estimation}

A second assessment regards the effect of the tolerance value used for rank estimation in \alg{polyorder_selfstab}. The study is conducted by investigating the condition number and polynomial orders for the different mesh types. 

The tolerance is varied starting from the Matlab default value: $${\mathrm{tol}_0 = \text{max} \plbr{\text{size} \plbr{\KE}} \times \epsilon \plbr{\norm{\KE}}},$$ and progressively increased by powers of ten. As a reference, the self-stabilized formulation $\mathrm{V3}$ \cite{berrone2025stabilization} is considered.

\begin{figure}[htbp]
        \centering
	\subfigure[Quadrilateral mesh. \label{fig:kp_quad}]{
		\includegraphics[width=0.3\textwidth]{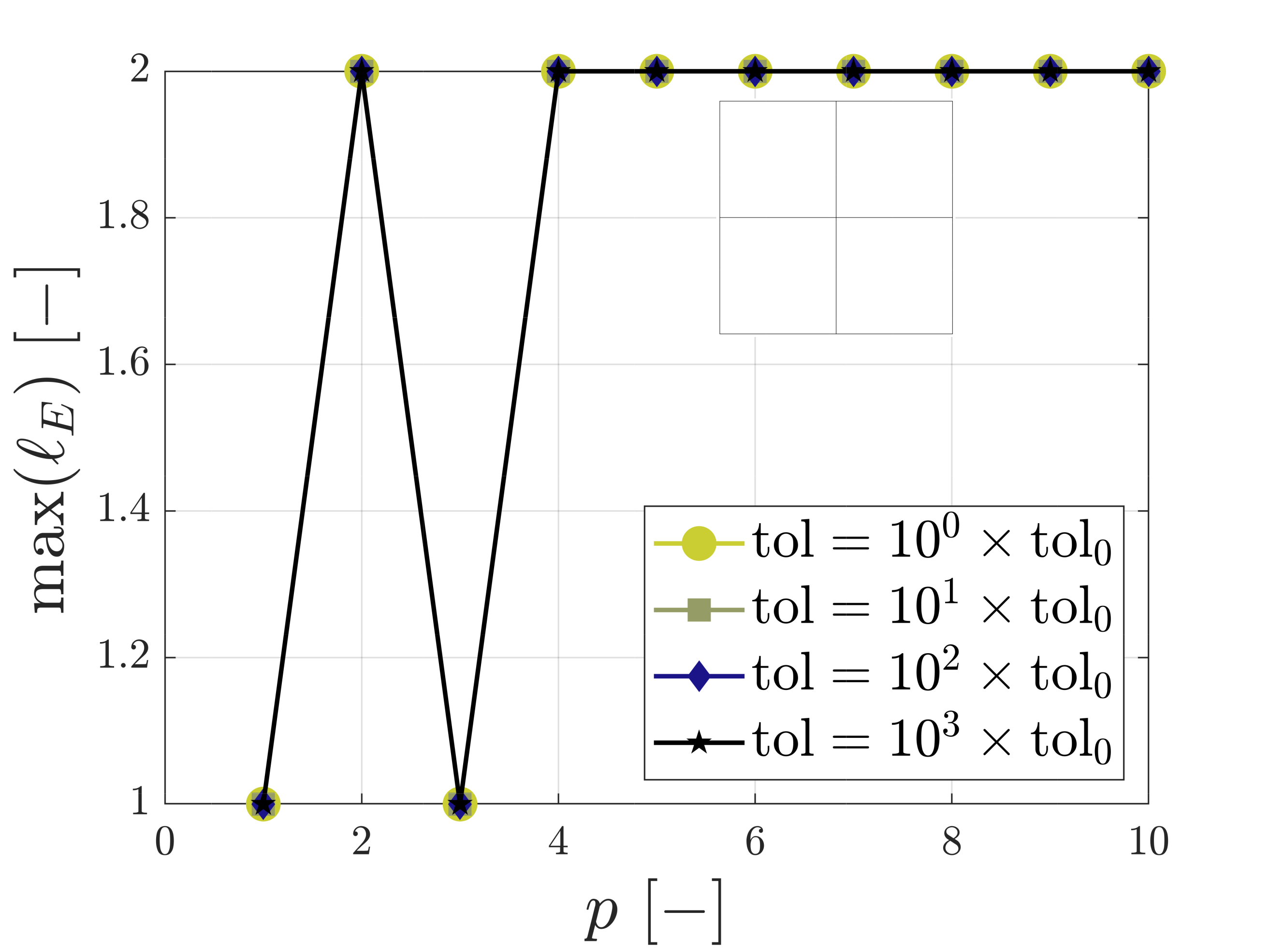}
	}
        \centering
	\subfigure[Voronoi mesh. \label{fig:kp_voronoi}]{
		\includegraphics[width=0.3\textwidth]{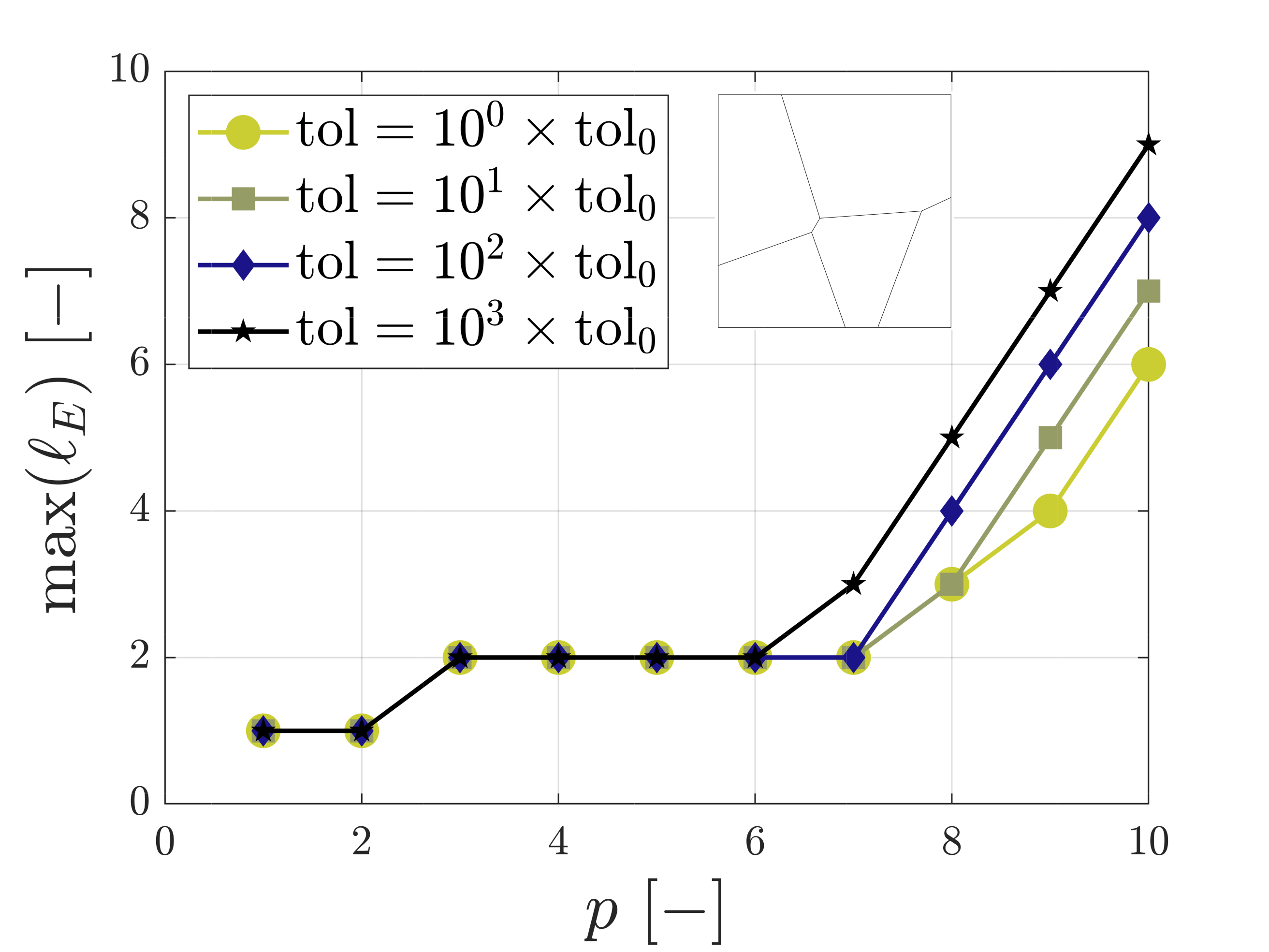}
	}     
    \centering
	\subfigure[Octagonal mesh. \label{fig:kp_ottagoni}]{
		\includegraphics[width=0.3\textwidth]{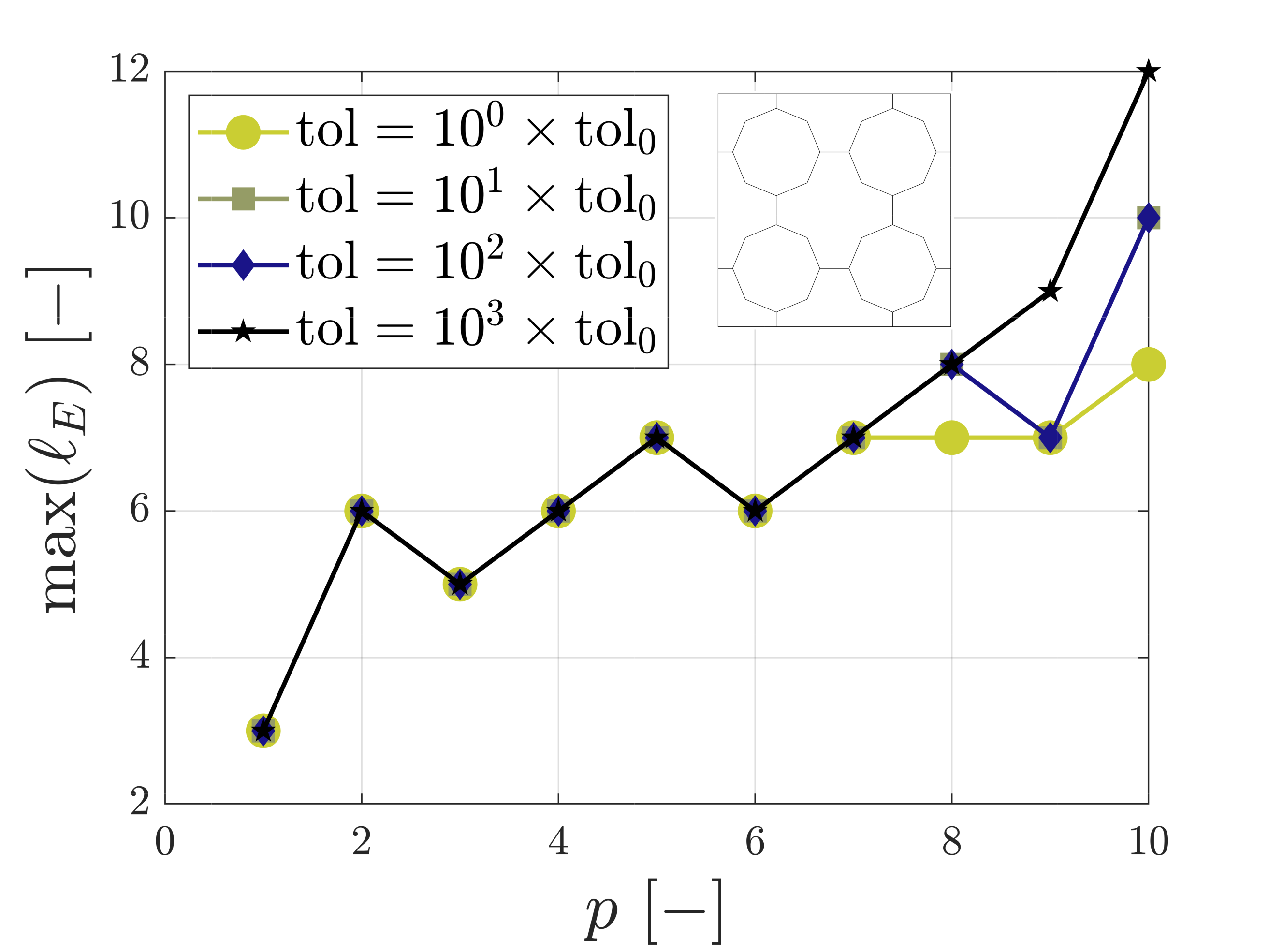}
	}    
    \caption{Laplace problem: maximum polynomial order for different tolerance values.\label{fig:svd_poly}}
 \end{figure}

For all mesh types, the detected polynomial order is reported in \fig{svd_poly}. The plots illustrate the maximum value for each mesh as, in general, each element is characterized by different $\laug$.

\begin{figure}[htbp]
        \centering
	\subfigure[Energy norm error. \label{fig:svd_err_quad}]{
		\includegraphics[width=0.47\textwidth]{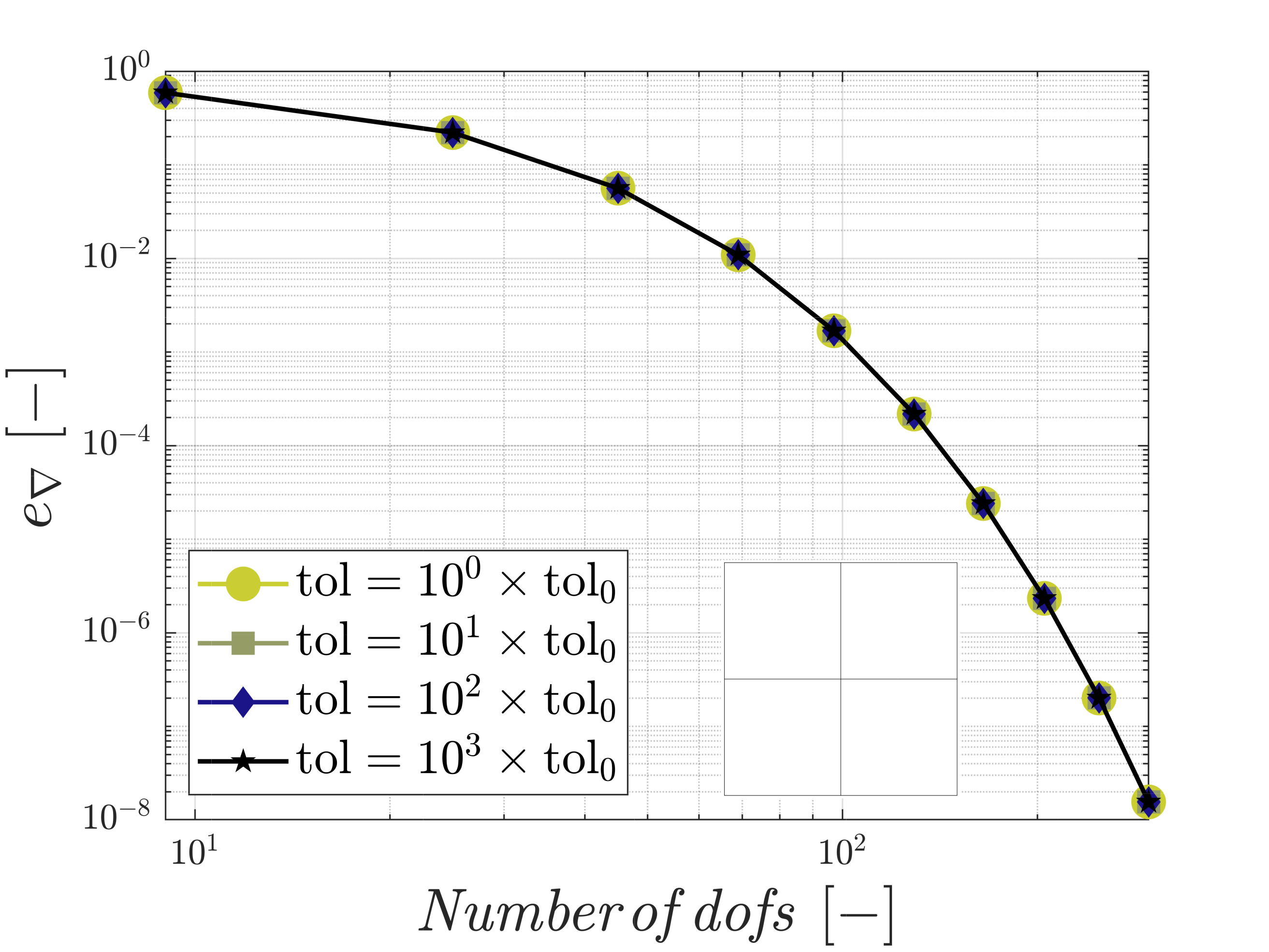}
	}
        \centering
        \hfill      
	\subfigure[Condition number. \label{fig:svd_cond_quad}]{
		\includegraphics[width=0.47\textwidth]{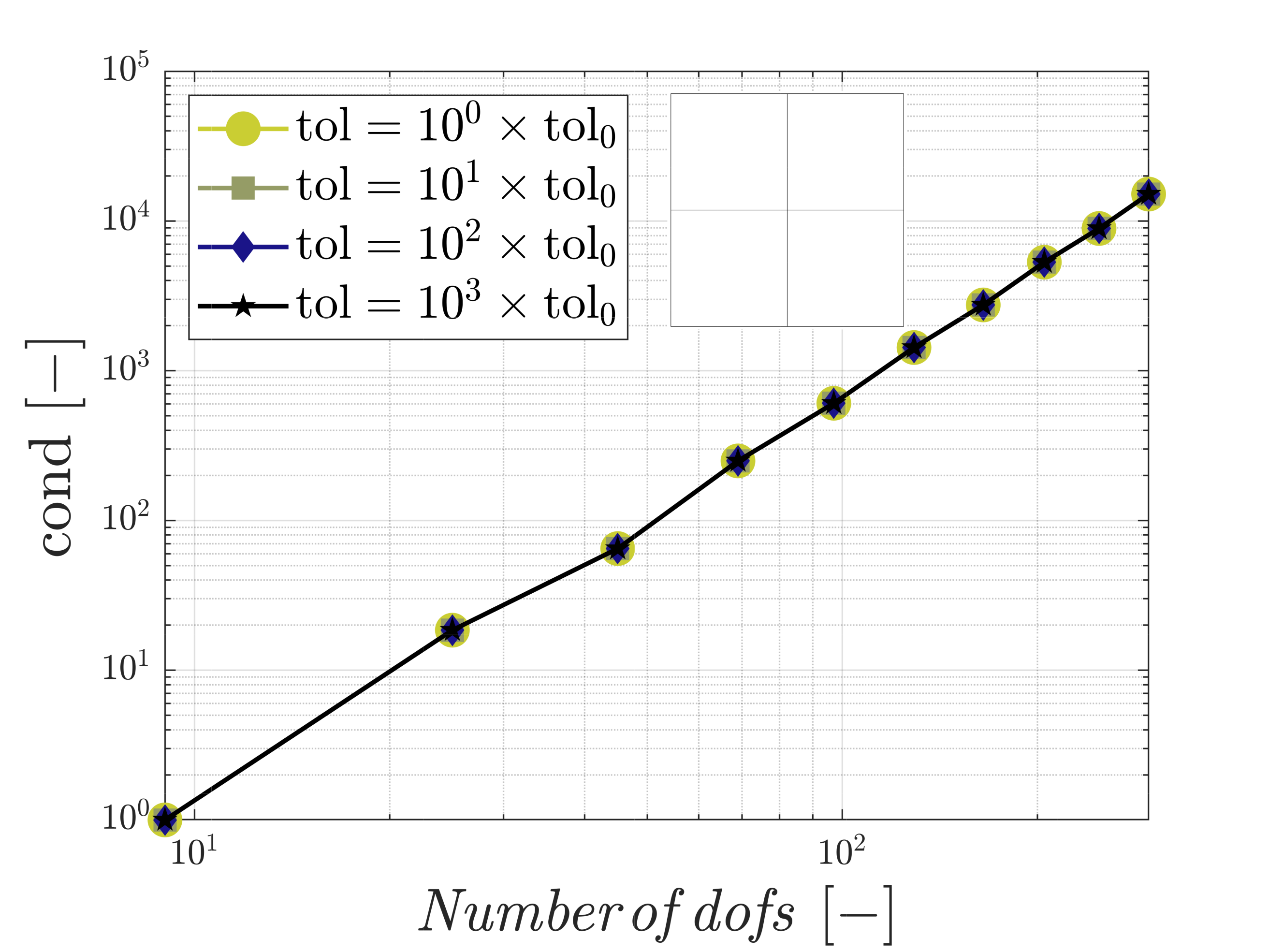}
	}     
    \caption{Laplace problem: parametric study on the tolerance for quadrilateral mesh.\label{fig:svd_quad}}
 \end{figure}

For the quadrilateral mesh shown in \fig{kp_quad}, the rank estimation tolerance has no effect on the detected polynomial order. Even if increased, the detected order remains $\laug \leq 2$.
As a consequence, both the energy norm error and the condition number are unaffected by this choice, as shown in \fig{svd_quad}. 

\begin{figure}[htbp]
        \centering
	\subfigure[Energy norm error. \label{fig:svd_err_voronoi}]{
		\includegraphics[width=0.47\textwidth]{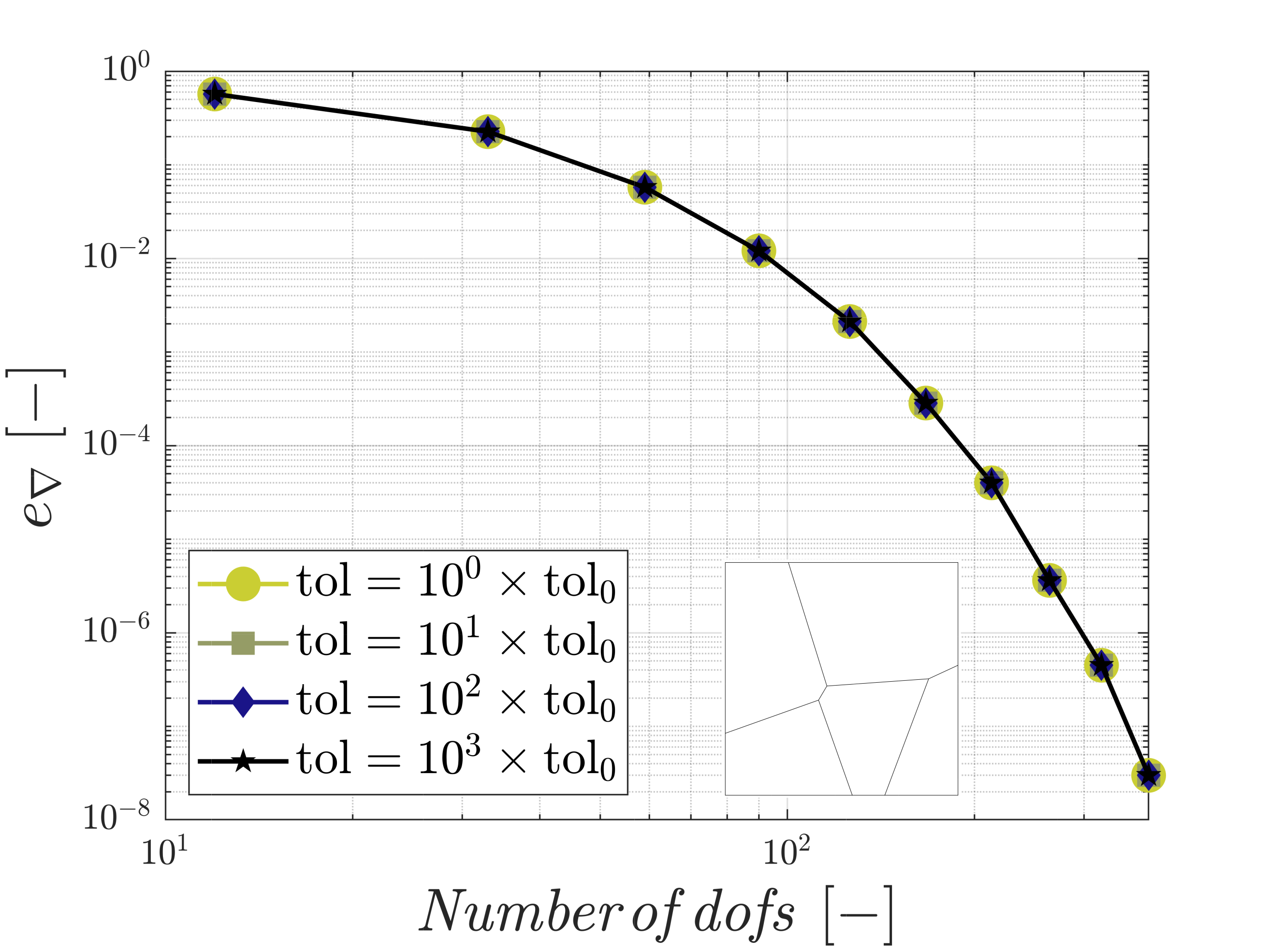}
	}
        \centering
        \hfill      
	\subfigure[Condition number. \label{fig:svd_cond_voronoi}]{
		\includegraphics[width=0.47\textwidth]{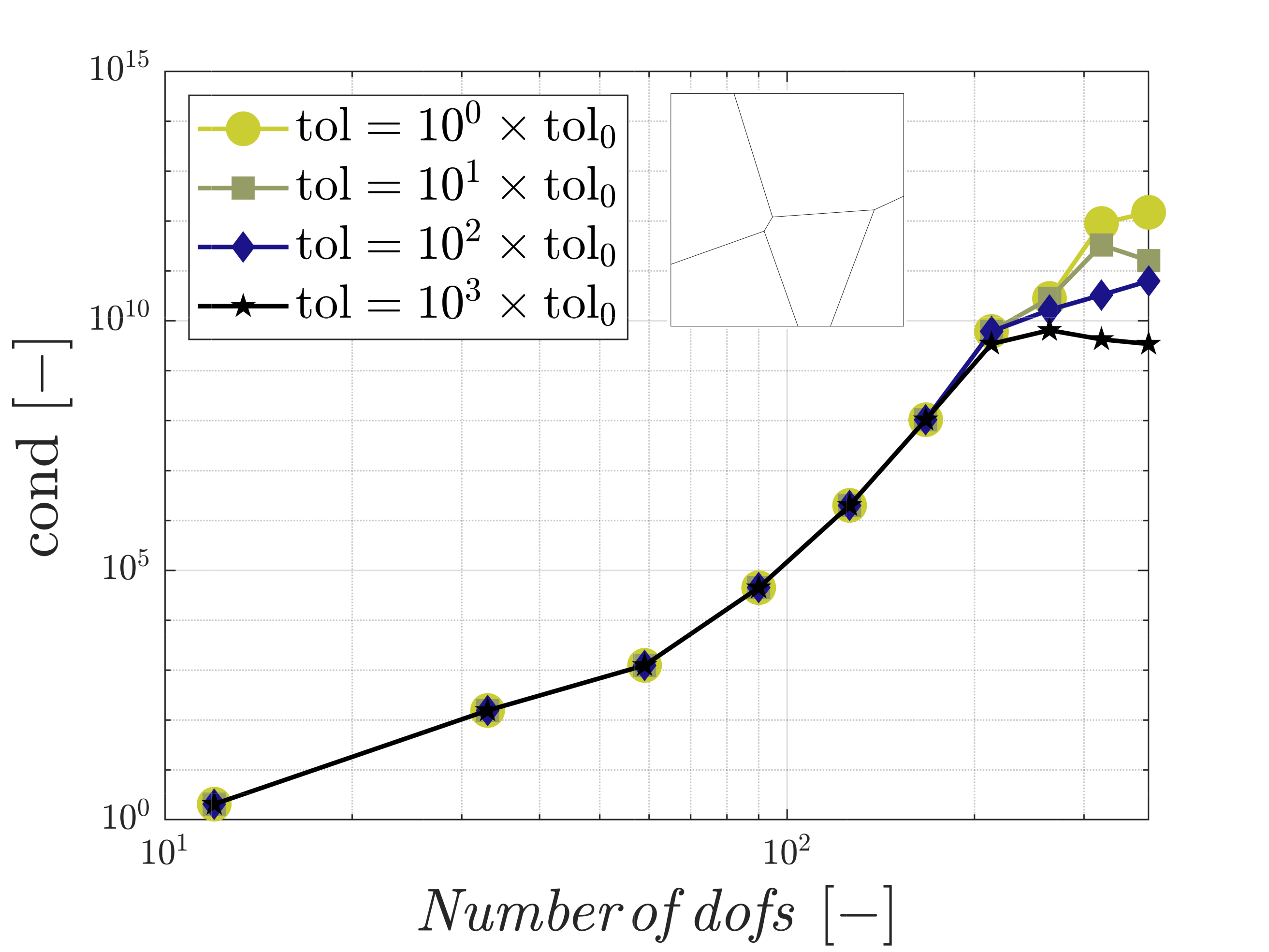}
	}     
    \caption{Laplace problem: parametric study on the tolerance for Voronoi mesh.\label{fig:svd_voronoi}}
 \end{figure}

The Voronoi mesh is investigated in \fig{svd_voronoi}. As seen, all the tests yield comparable accuracy and condition number up to the order $\kord=6$. 
On the contrary, for $\kord > 6$, improved conditioning can be observed for increasing tolerance values. This behavior is explained by referring back to \fig{kp_voronoi}. The increase in the order $\laug$ has the effect of shifting the system away from the nearly singular regime, resulting in improved conditioning. 

\begin{figure}[htbp]
        \centering
	\subfigure[Energy norm error. \label{fig:svd_err_ottagoni}]{
		\includegraphics[width=0.47\textwidth]{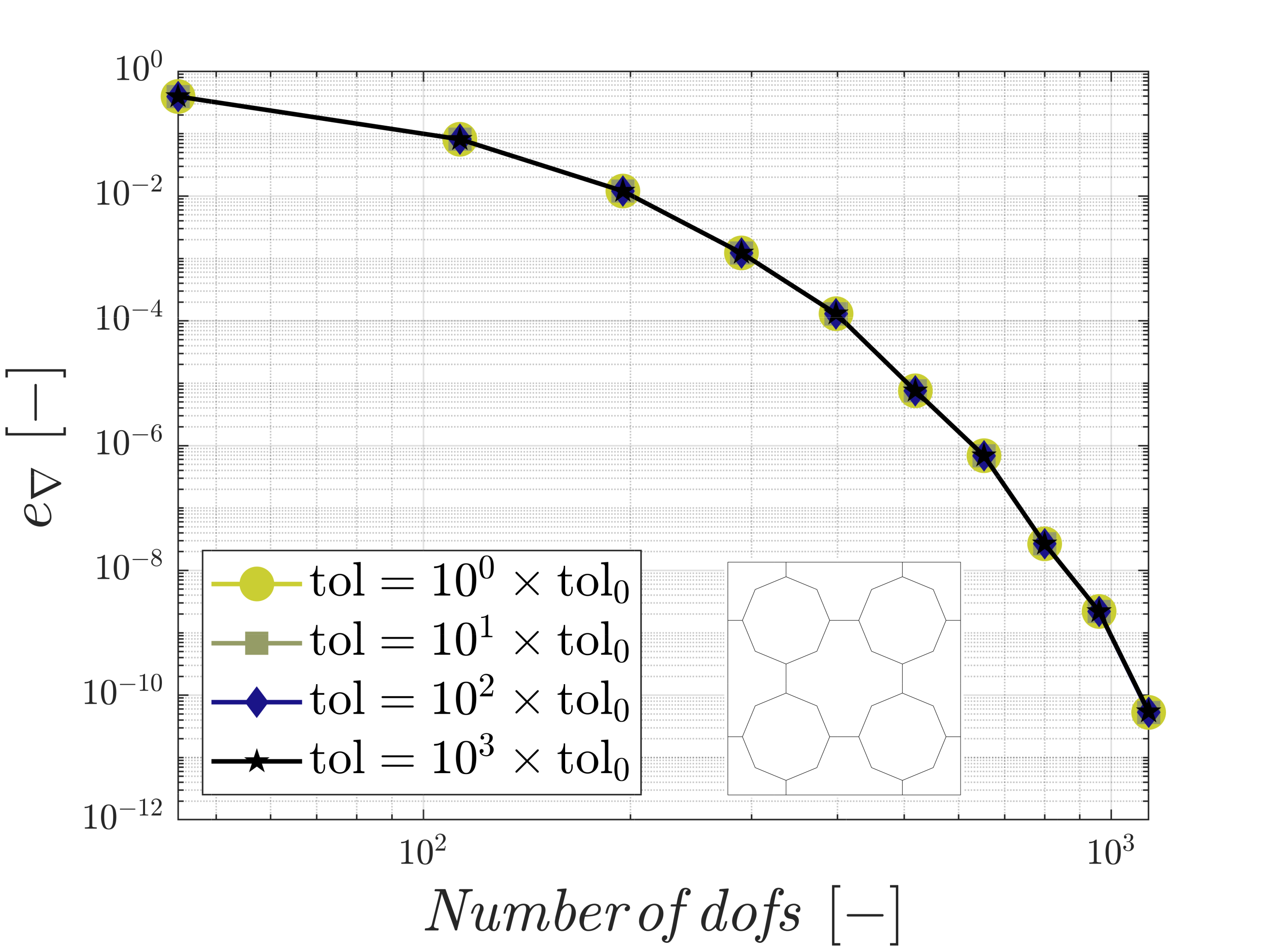}
	}
        \centering
        \hfill      
	\subfigure[Condition number. \label{fig:svd_cond_ottagoni}]{
		\includegraphics[width=0.47\textwidth]{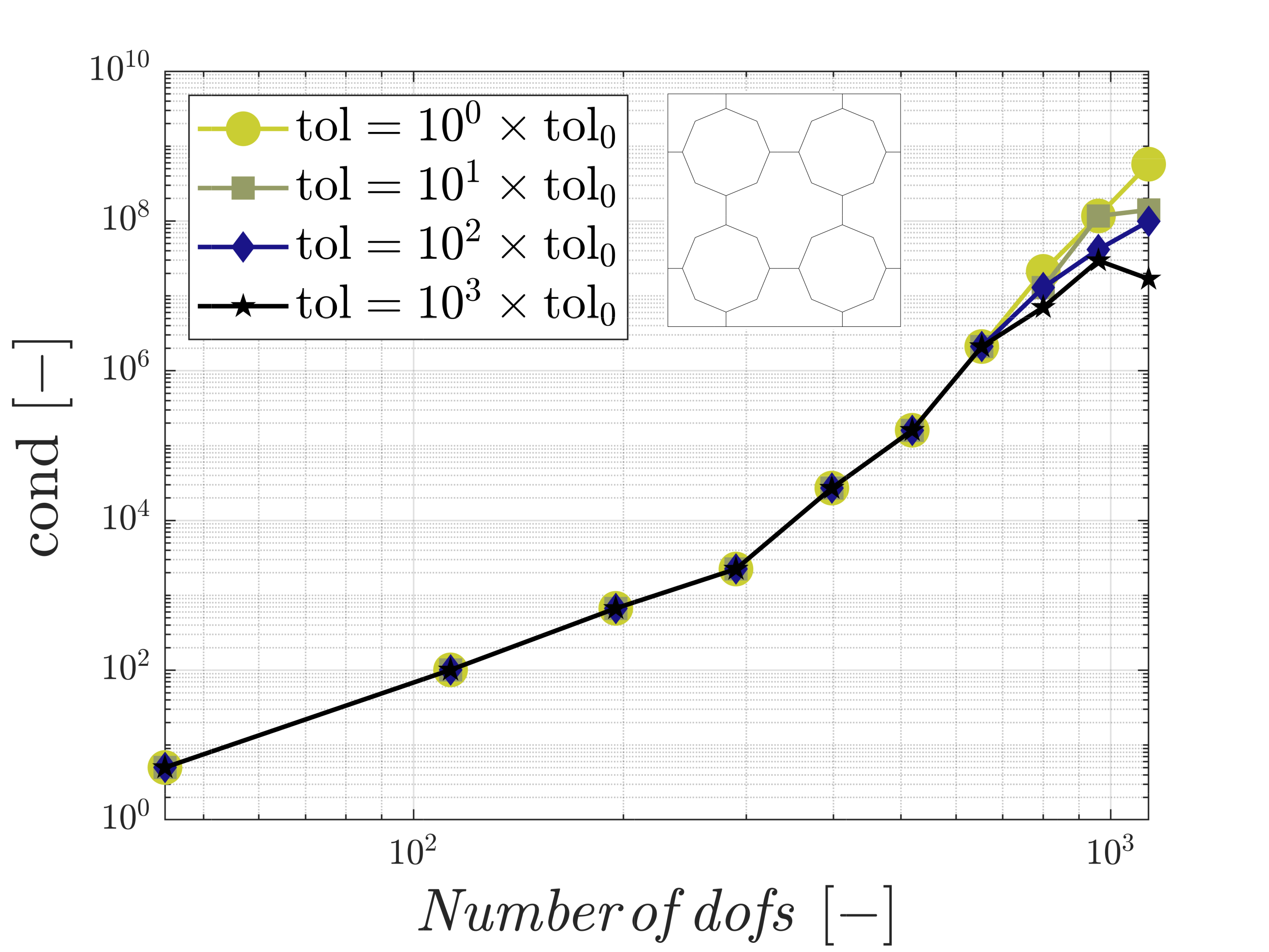}
	}     
    \caption{Laplace problem: parametric study on the tolerance for octagonal mesh.\label{fig:svd_ottagoni}}
 \end{figure}

The octagonal mesh is the subject of the last investigation. The results are presented in \fig{svd_ottagoni}. As seen, all the simulations share similar accuracy. On the contrary, some deviations can be noted for $\kord>7$ in terms of conditioning. As for Voronoi meshes, a higher tolerance corresponds to a higher polynomial order and so to improved conditioning.

The maximum $\laug$ in \fig{kp_ottagoni} corresponds to the central element. Indeed, in general, the higher the number of edges, the higher the polynomial order to guarantee full rank \cite{berrone2025stabilization}.

The choice of a coarse tolerance can mitigate the increase of the condition number just observed. Indeed, the obtained value of $\laug$ is such that the resulting stiffness matrix is far from being singular. Worth of mention is, in any case, the higher computational cost.

\subsubsection{Performance of stabilized and self-stabilized formulations}

An assessment is now presented to compare the stabilized and self-stabilized formulations. For the former, the parameter $\tau=1$ is used. For the latter, the tolerance in the rank estimation is set to $\mathrm{tol}=\mathrm{tol}_0$. This value is found as a good tradeoff between computational cost and accuracy.

\begin{figure}[htbp]
        \centering
	\subfigure[Energy norm error. \label{fig:stability_err_quad}]{
		\includegraphics[width=0.47\textwidth]{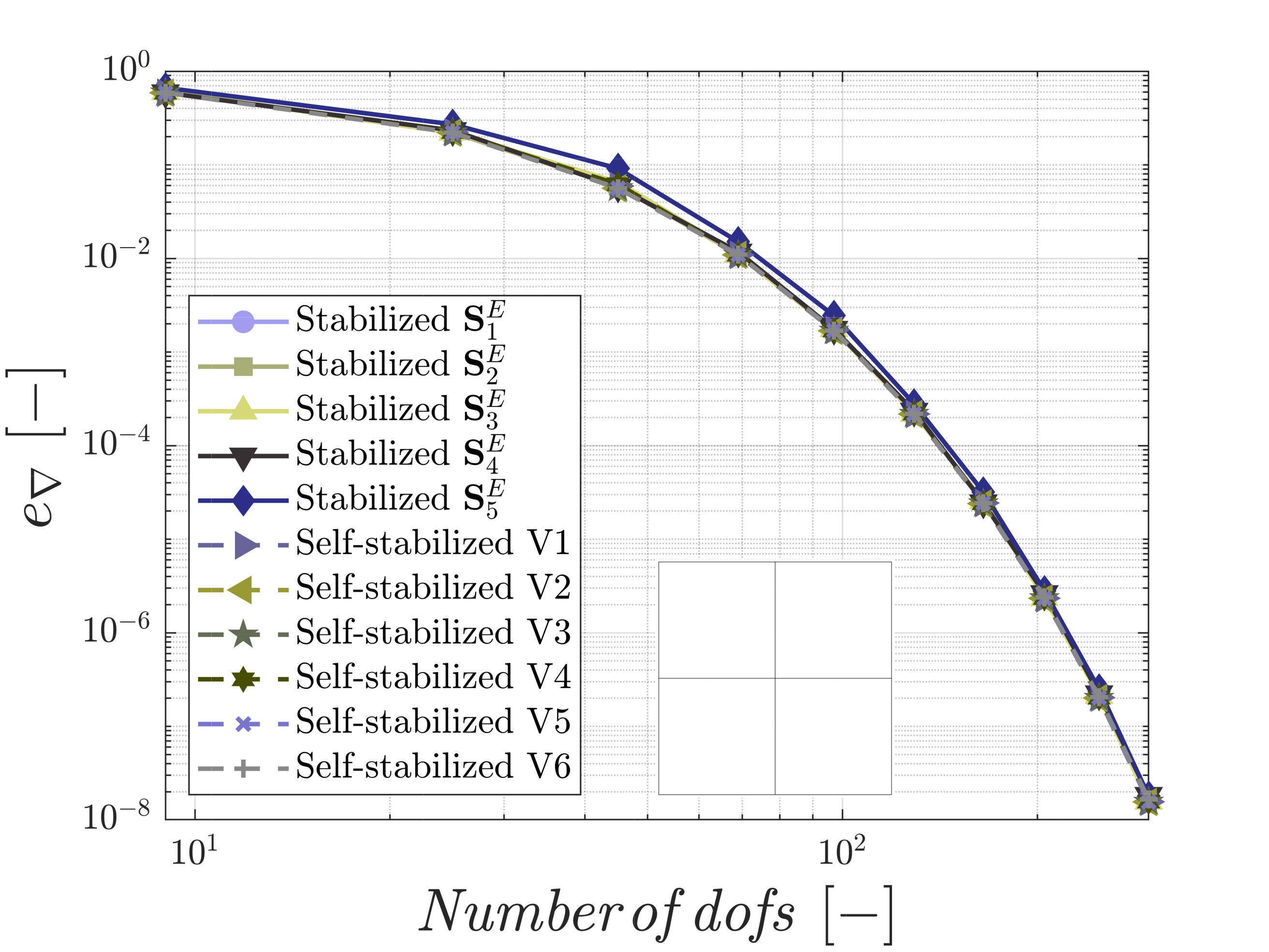}
	}
        \centering
        \hfill   
	\subfigure[Condition number. \label{fig:stability_cond_quad}]{
		\includegraphics[width=0.47\textwidth]{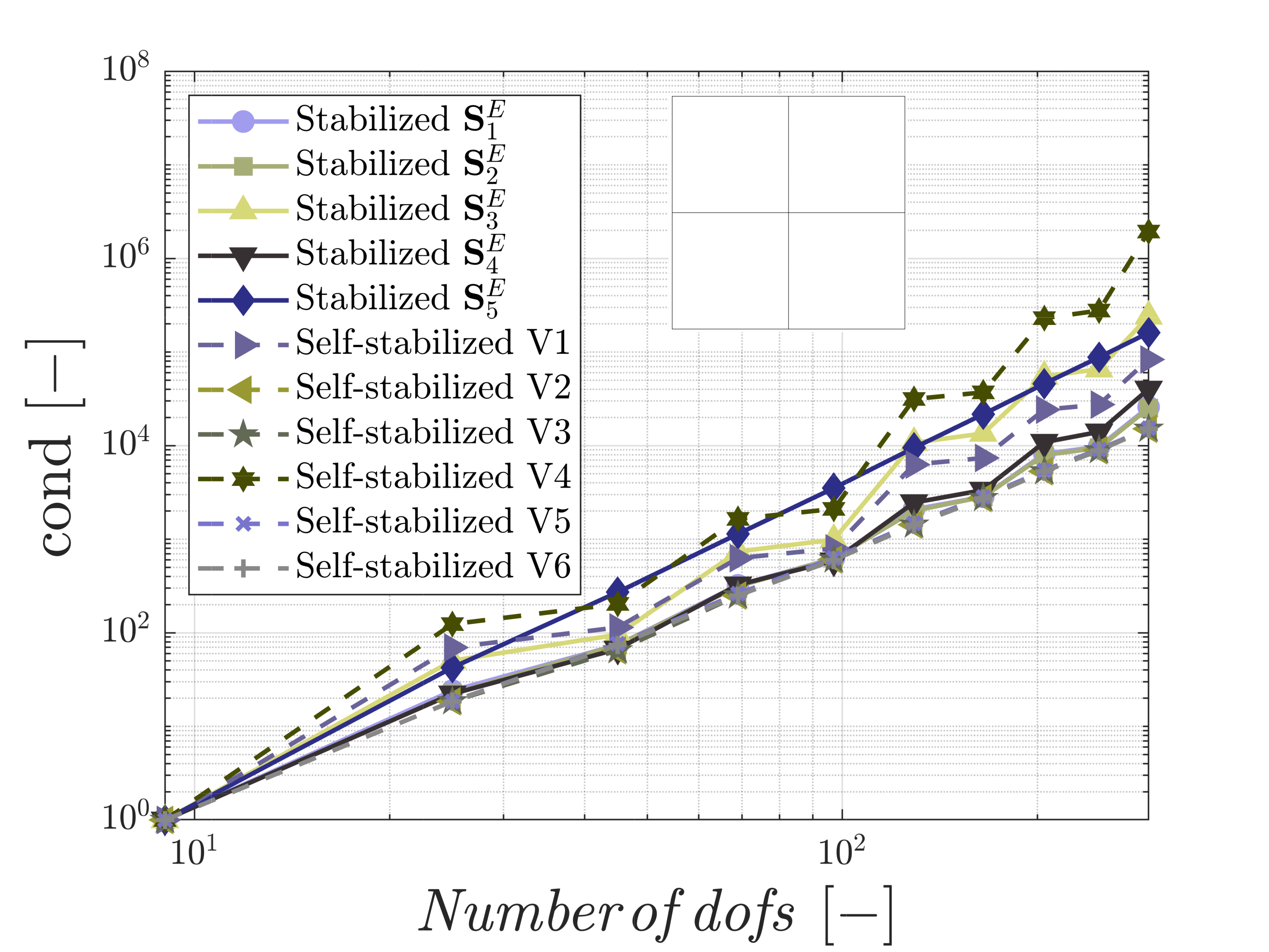}
	}     
    \caption{Laplace problem: performance of stabilized and self-stabilized formulations for quadrilateral mesh.\label{fig:stability_quad}}
 \end{figure}

The results obtained using the quadrilateral mesh are reported in \fig{stability_quad}. In this case, all the VEM variants lead to similar accuracy and condition number.

\begin{figure}[htbp]
        \centering
	\subfigure[Energy norm error. \label{fig:stability_err_voronoi}]{
		\includegraphics[width=0.47\textwidth]{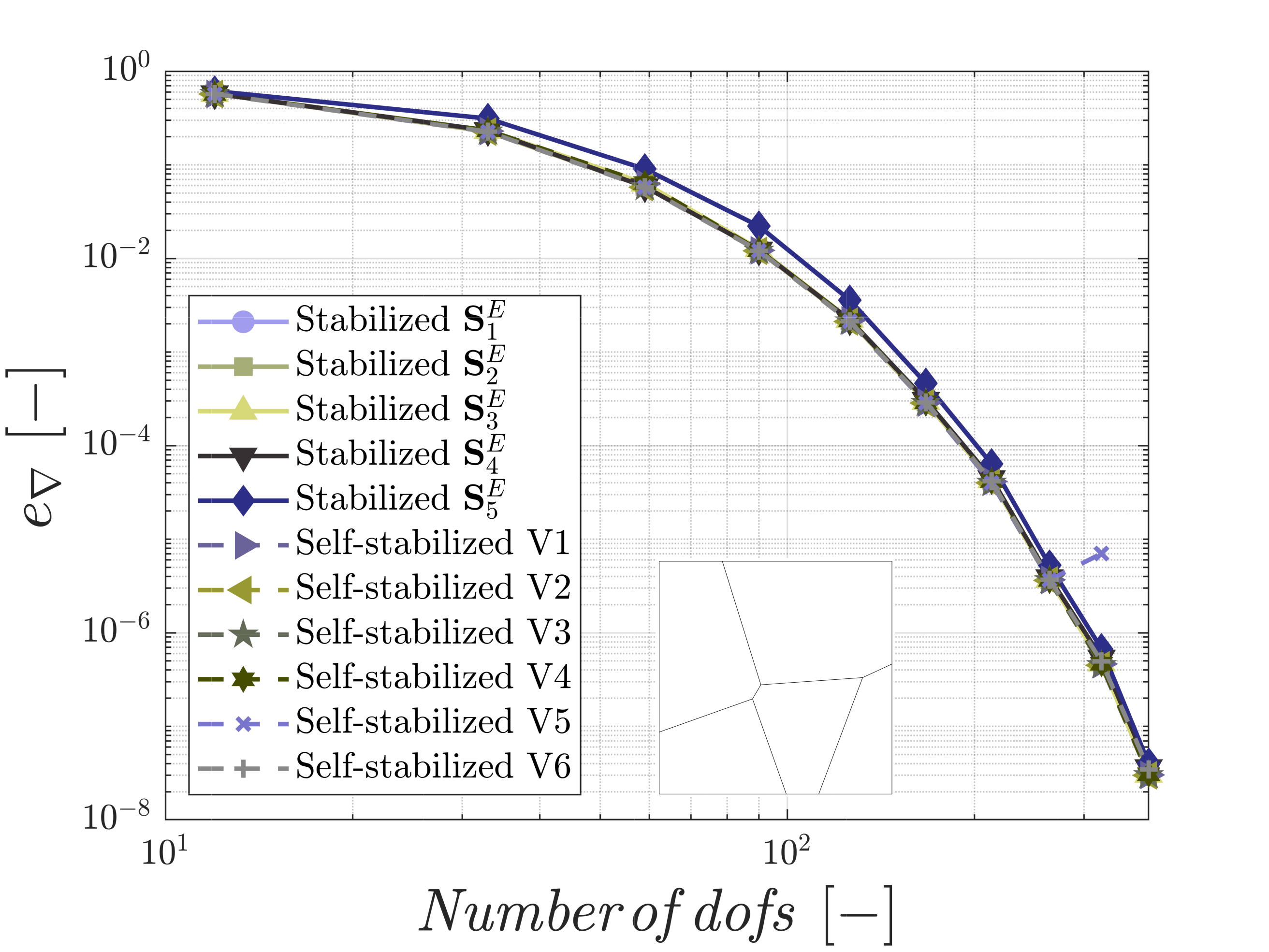}
	}
        \centering
        \hfill      
	\subfigure[Condition number. \label{fig:stability_cond_voronoi}]{
		\includegraphics[width=0.47\textwidth]{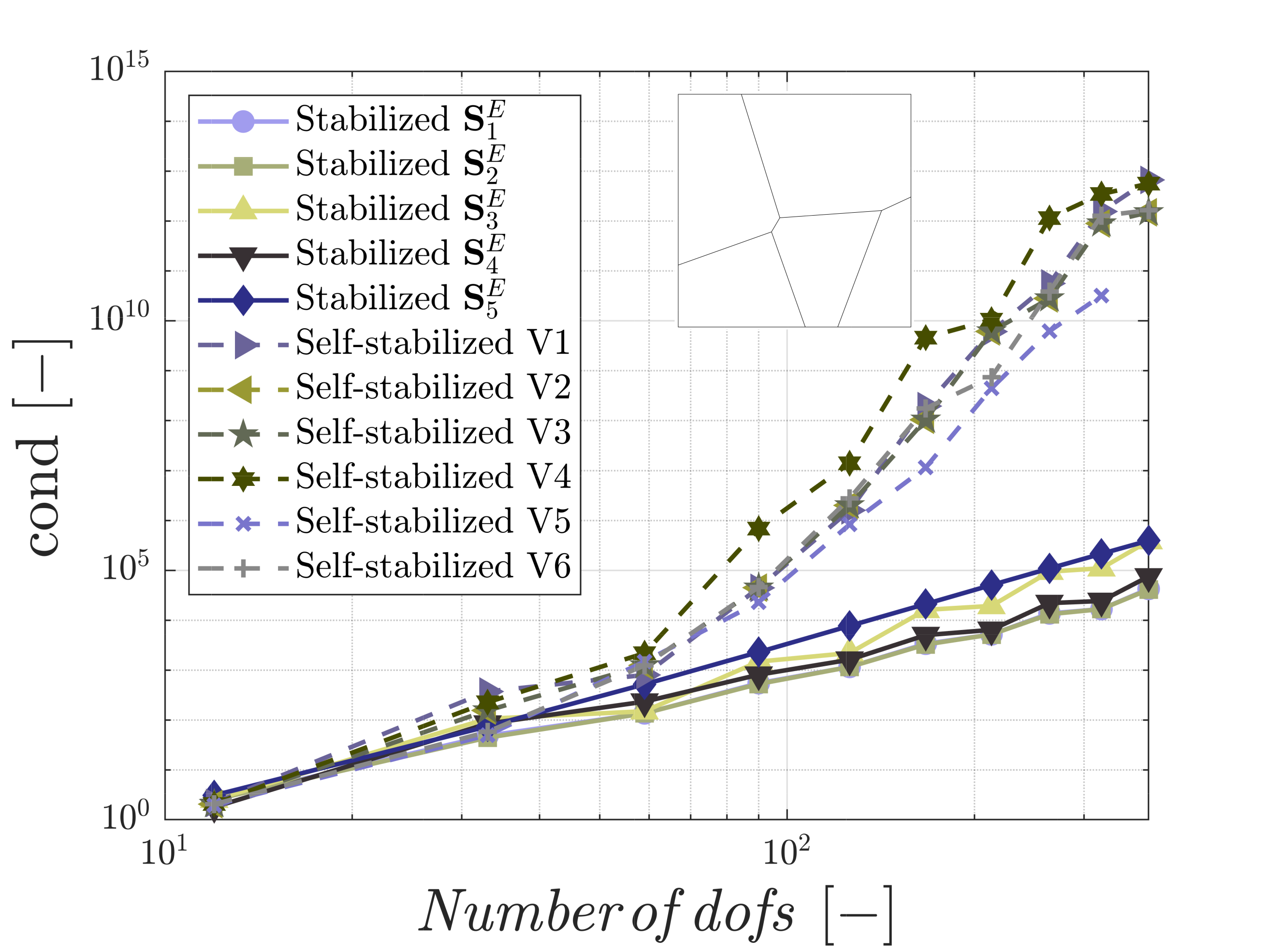}
	}     
    \caption{Laplace problem: performance of stabilized and self-stabilized formulations for Voronoi mesh.\label{fig:stability_voronoi}}
 \end{figure}

The results for the Voronoi mesh are presented in \fig{stability_voronoi}. All the test cases display similar accuracy. The only exception is the self-stabilized $\mathrm{V5}$ that diverges and is therefore stopped at $\kord=9$. For the curves of \fig{stability_cond_voronoi}, the self-stabilized formulations generally exhibit worse conditioning compared to the classical VEM.

\begin{figure}[htbp]
        \centering
	\subfigure[Energy norm error. \label{fig:stability_err_ottagoni}]{
		\includegraphics[width=0.47\textwidth]{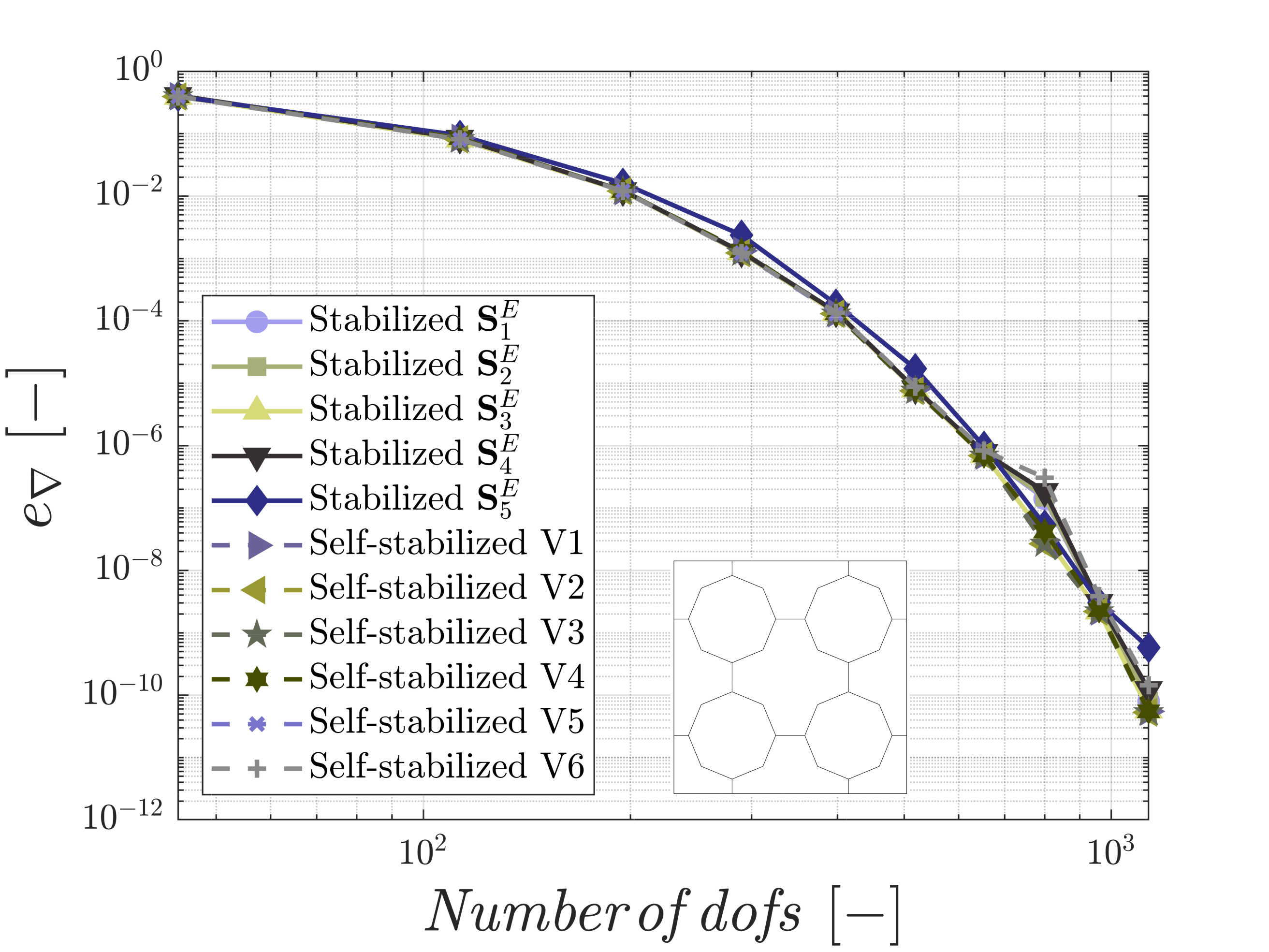}
	}
        \centering
        \hfill      
	\subfigure[Condition number. \label{fig:stability_cond_ottagoni}]{
		\includegraphics[width=0.47\textwidth]{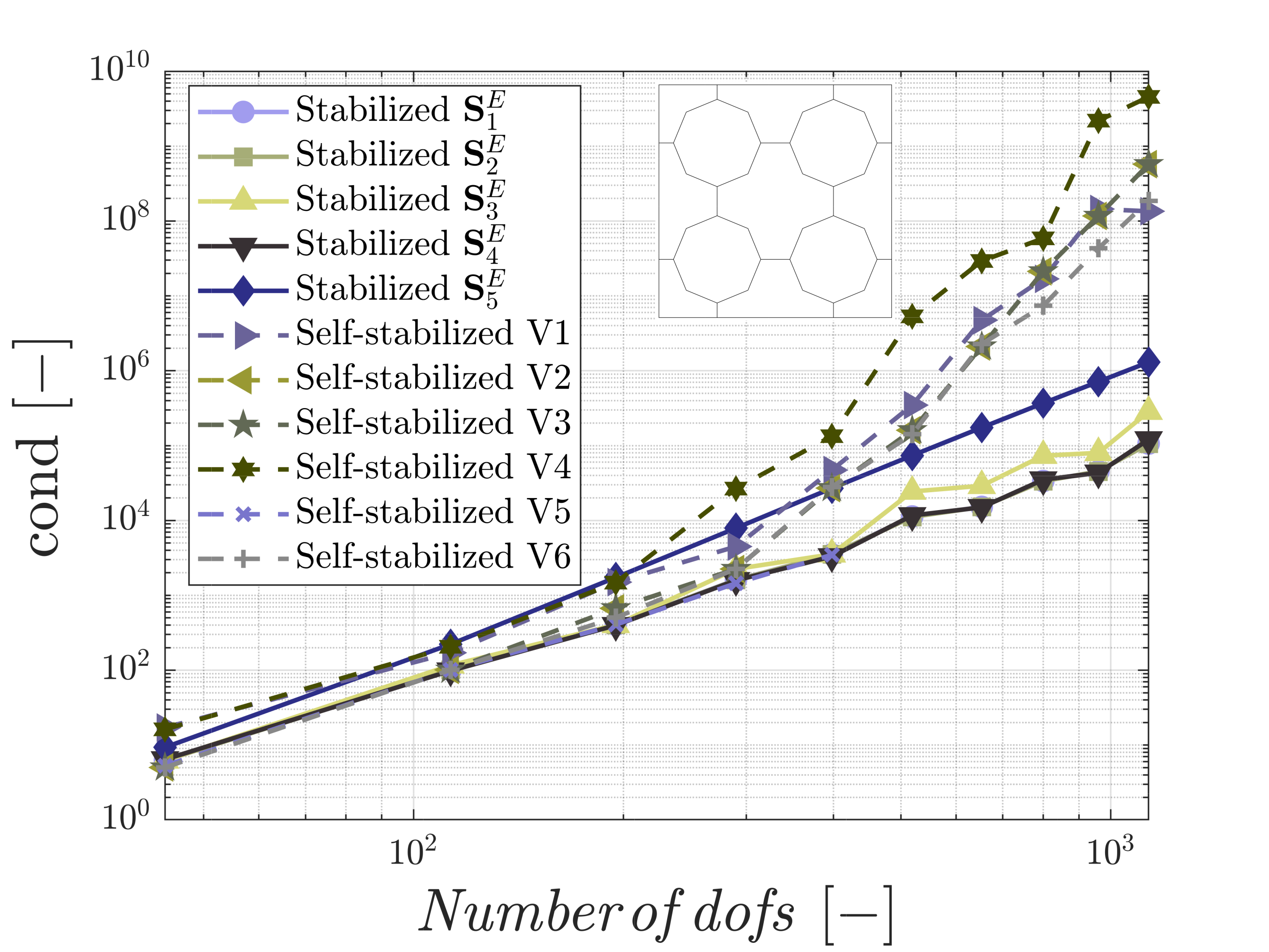}
	}     
    \caption{Laplace problem: performance of stabilized and self-stabilized formulations for octagonal mesh.\label{fig:stability_ottagoni}}
 \end{figure}

The results for the octagonal mesh are shown in \fig{stability_ottagoni}. Apart from similar accuracy exhibited by all cases, it is interesting to observe that the curve for the stabilization $\KEs_5$ features, for $\kord=10$, a slightly larger error. 
Also, the self-stabilized $\mathrm{V5}$ has been stopped at $\kord=5$ due to an excessive increase in the polynomial order, having already reached $\text{max}\plbr{\laug}=25$. 
The behavior in the condition number is the same observed for the Voronoi mesh, i.e. the self-stabilized versions feature worse conditioning than the stabilized ones, although the difference is less significant in this case.

\subsection{Linear elasticity problem} \label{subsec:Linerelasticity_results_selfstab}

The linear elasticity problem is now analyzed. The same inverse approach employed for the Laplace problem is considered. In particular, the load term and boundary conditions are selected such that the exact solution is:
\begin{equation}
    \disp \plbr{x,y} = 
    \begin{bmatrix}
    \sin\plbr{\pi x} \sin\plbr{\pi y}\\
    \sin\plbr{\pi x} \sin\plbr{\pi y}
    \end{bmatrix},
    \label{eq:linearelasticity_exactsolution}
\end{equation}
in $\surom = \plbr{0,1}^2$. The material is isotropic with $E=72000$ and $\nu=0.3$, thus $\cost$ is constant.

The convergence is now investigated both the $L^2$ and the energy norm, whose definitions read:
\begin{equation}
\begin{aligned}
    &e_{L^2}^{\disp} =  \frac{ \plbr{\sum_{\el \in \globdom}  \norm{ \disp - \PiokK \dispd }_{0,\el}^2  }^{\frac{1}{2}} }{ \norm{\disp}_{0,\el} },\\
    &e_{\bm{\varepsilon}}^{\disp} =  \frac{ \plbr{\sum_{\el \in \globdom}  \norm{ \sqrt{\cost} \strain \plbr{ \disp - \PiokK \dispd } }_{0,\el}^2  }^{\frac{1}{2}} }{ \norm{ \sqrt{\cost} \strain \plbr{\disp}}_{0,\el} }.
\end{aligned}    
    \label{eq:linearelasticity_l2error_dispenergy}
\end{equation}

The analysis is conducted on the quadrilateral and on the B\'ezier-edge meshes, with order $\kord=1,\dots,10$. 
The comparison between the stabilized and self-stabilized formulations is presented in terms of accuracy by setting $\tau=0.5$ -- following the guidelines of  Ref. \cite{artioli2017arbitrary}) -- and $\mathrm{tol}=\mathrm{tol}_0$. 

\begin{figure}[htbp]
        \centering
	\subfigure[$L^2$ error. \label{fig:stability_err_disp_quad}]{
		\includegraphics[width=0.47\textwidth]{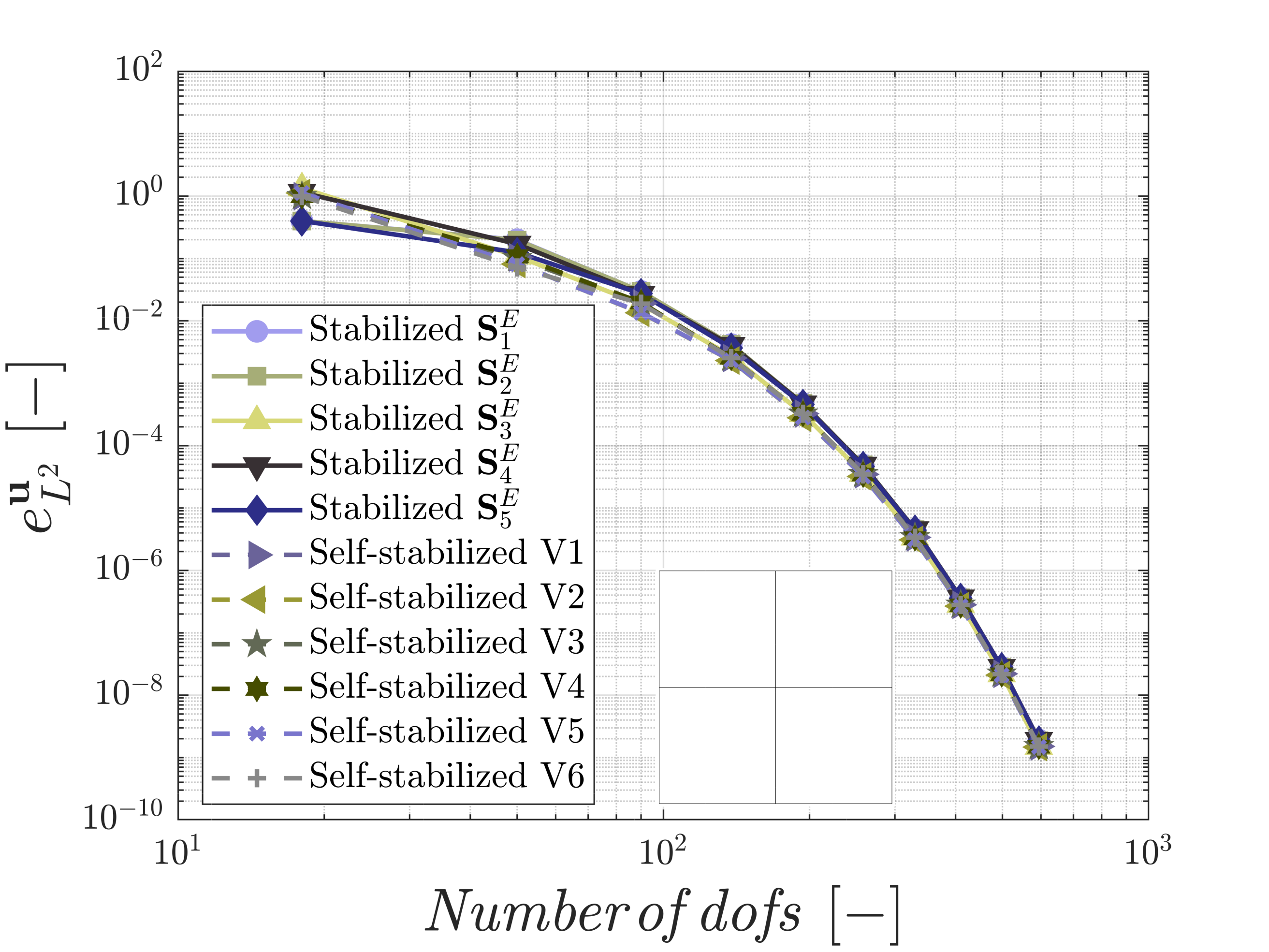}
	}
        \centering
        \hfill   
	\subfigure[Energy norm error. \label{fig:stability_err_stress_quad}]{
		\includegraphics[width=0.47\textwidth]{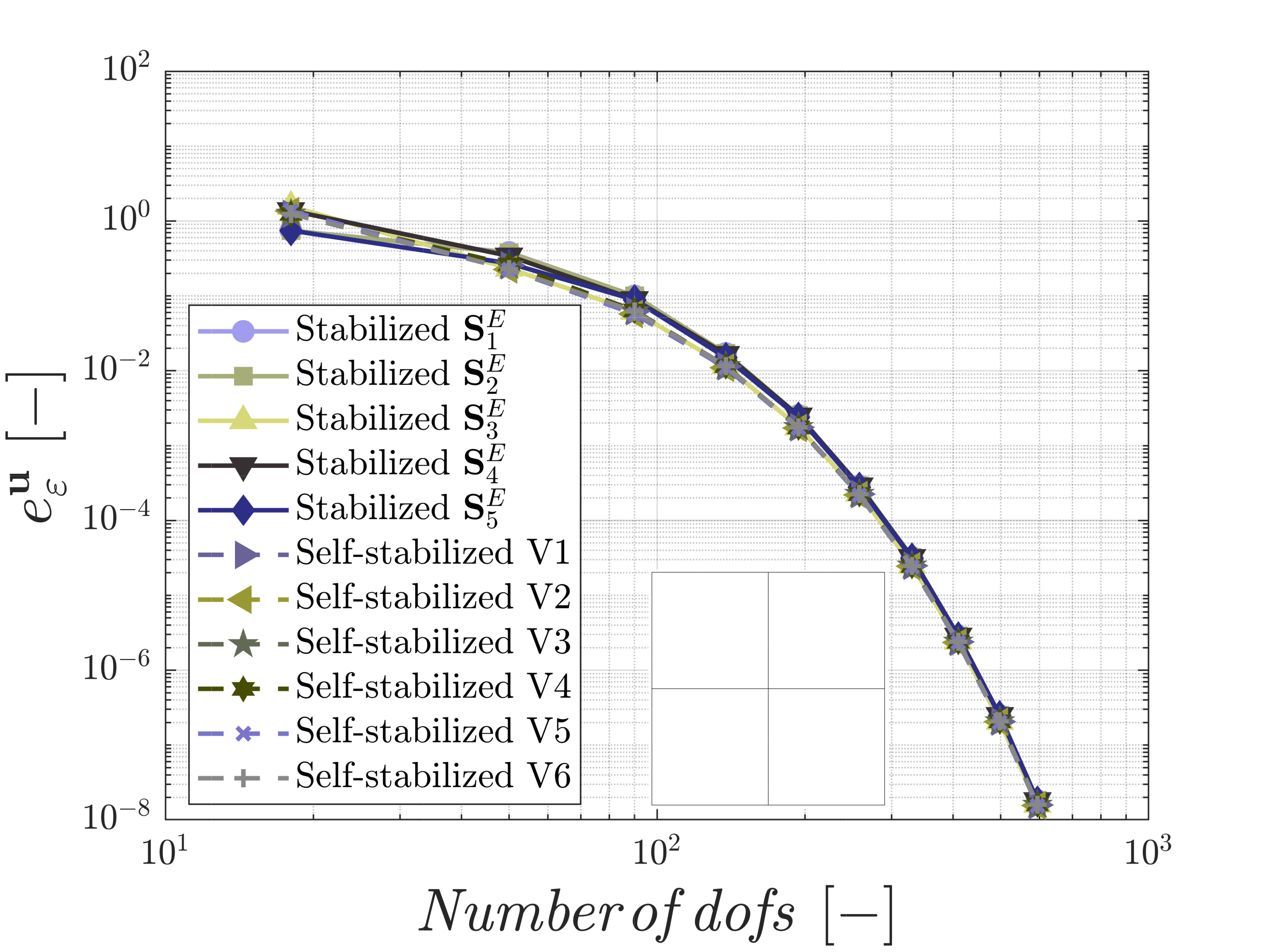}
	}     
    \caption{Linear elasticity problem: performance of stabilized and self-stabilized formulations for quadrilateral mesh.\label{fig:stability_linearelasticity_quad}}
\end{figure}

The results for the quadrilateral mesh are shown in \fig{stability_linearelasticity_quad}. All simulations have same accuracy both in the $L^2$ and in the energy norm.

\begin{figure}[htbp]
    \centering
    \subfigure[$L^2$ displacement error. \label{fig:stability_err_disp_curv}]{
        \includegraphics[width=0.47\textwidth]{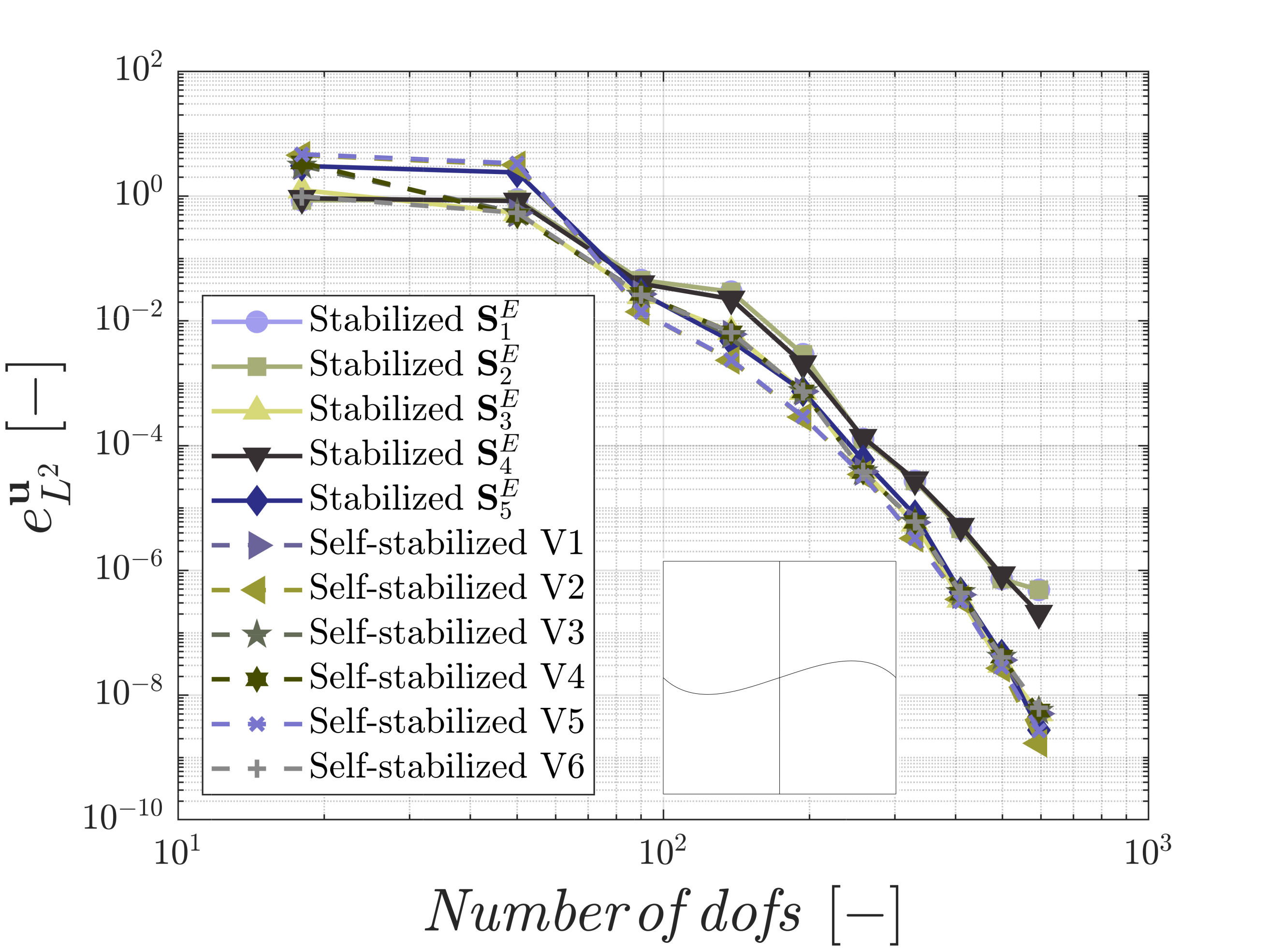}
    }
        \centering
        \hfill   
    \subfigure[Energy norm error. \label{fig:stability_err_stress_curv}]{
        \includegraphics[width=0.47\textwidth]{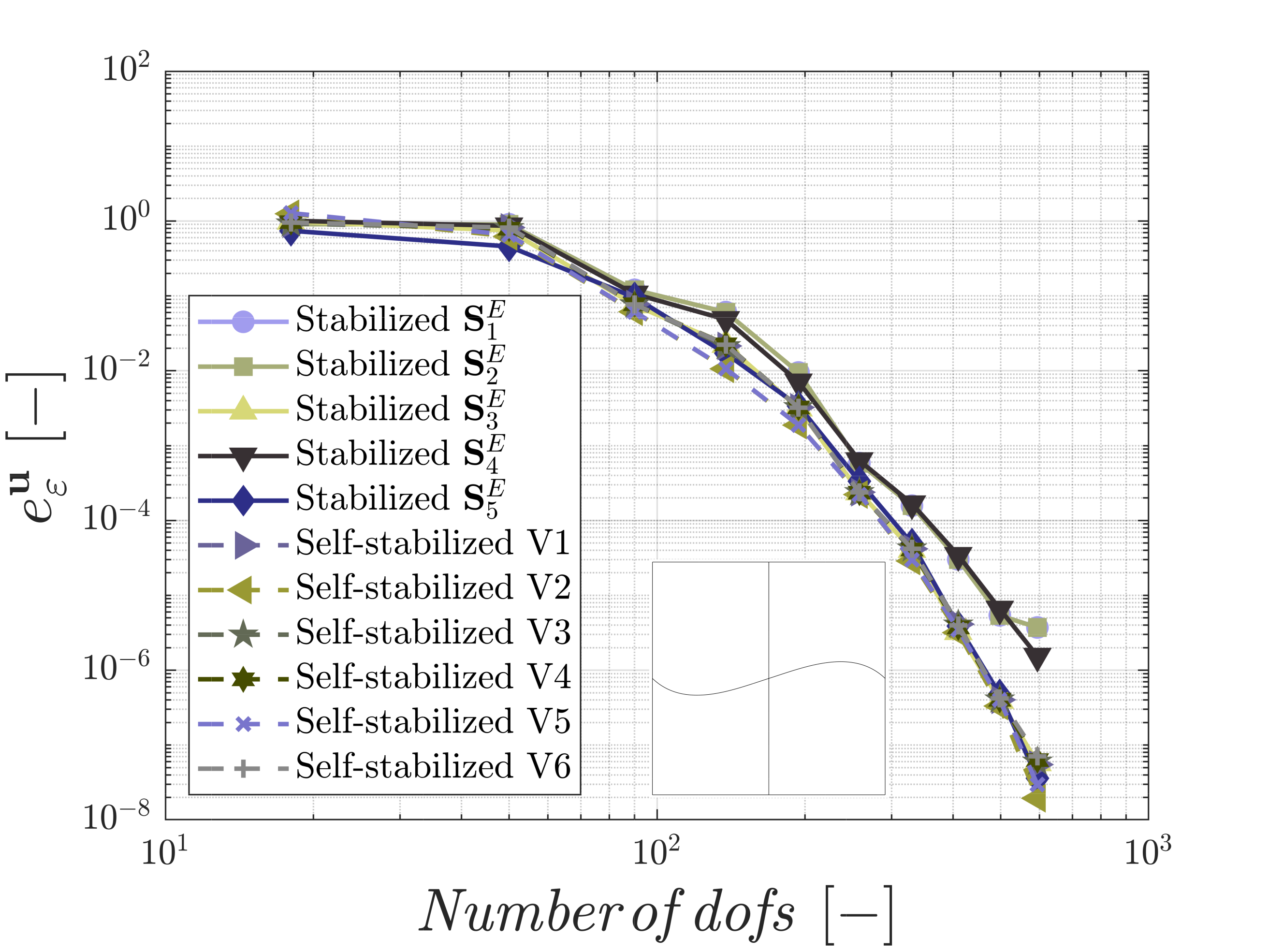}
    }     
    \caption{Linear elasticity problem: performance of stabilized and self-stabilized formulations for B\'ezier-edge mesh.\label{fig:stability_linearelasticity_curv}}
\end{figure}

The results for the B\'ezier-edge mesh are shown in \fig{stability_linearelasticity_curv}. The trends for both $L^2$ and energy errors are similar. The convergence rate is constant at lower-orders $\kord$, but improves as the approximation order increases. This may be ascribed to the presence of curved edges, which are not well captured by low-order approximations. It is clear that all the self-stabilized versions and the stabilized formulations $\KEs_3$ and $\KEs_5$ feature an optimal accuracy. Conversely, the stabilized formulations $\KEs_1$, $\KEs_2$ and $\KEs_4$ exhibit larger errors, indicating that the choice of stabilization has a more relevant impact on the results in the presence of curved edges.

\subsection{Stokes problem}

\begin{figure}[htbp]
        \centering
	\subfigure[Energy norm velocity error. \label{fig:stability_err_vel_voronoi}]{
		\includegraphics[width=0.47\textwidth]{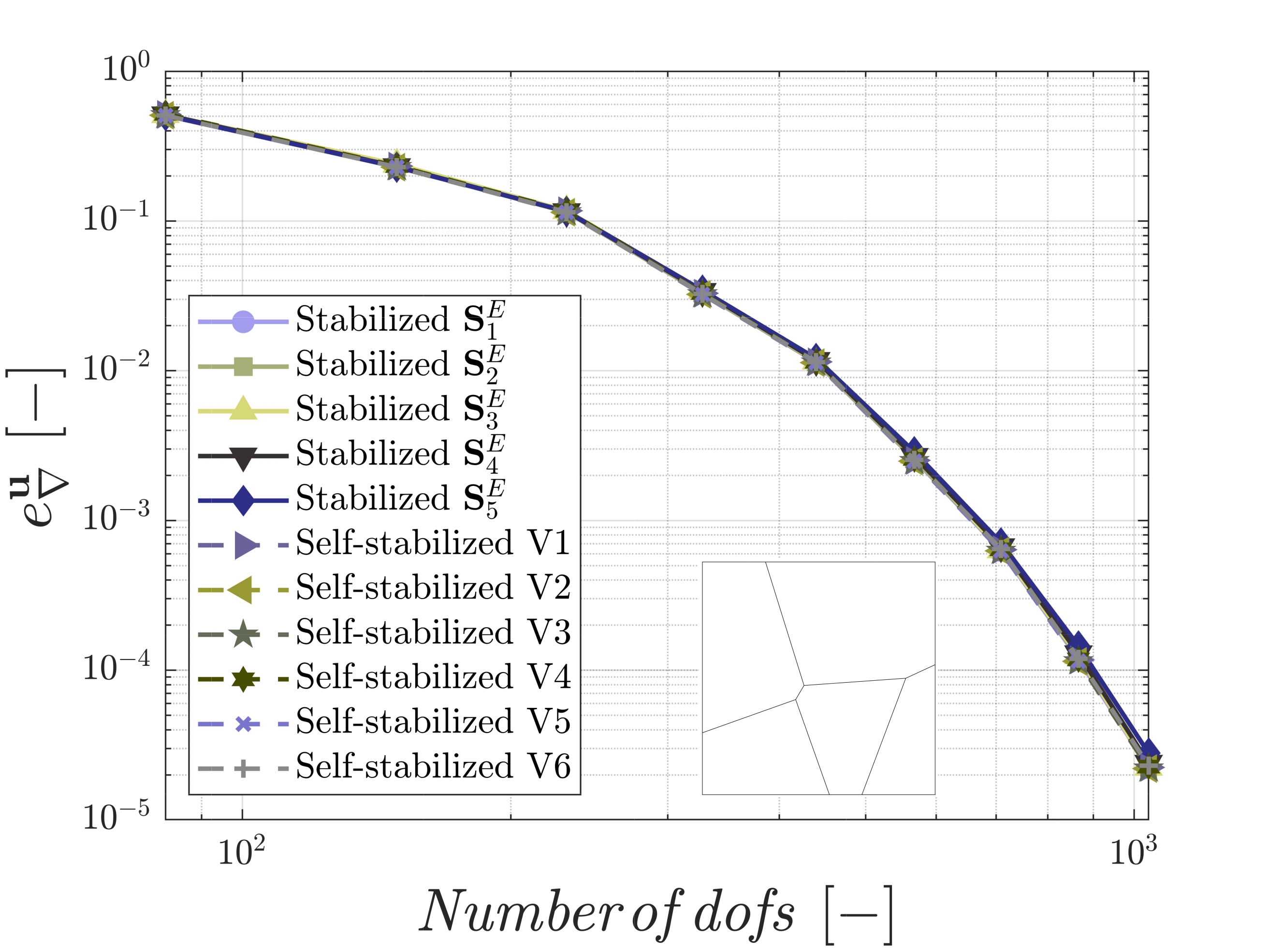}
	}
        \centering
        \hfill   
	\subfigure[$L^2$ pressure error. \label{fig:stability_err_pres_voronoi}]{
		\includegraphics[width=0.47\textwidth]{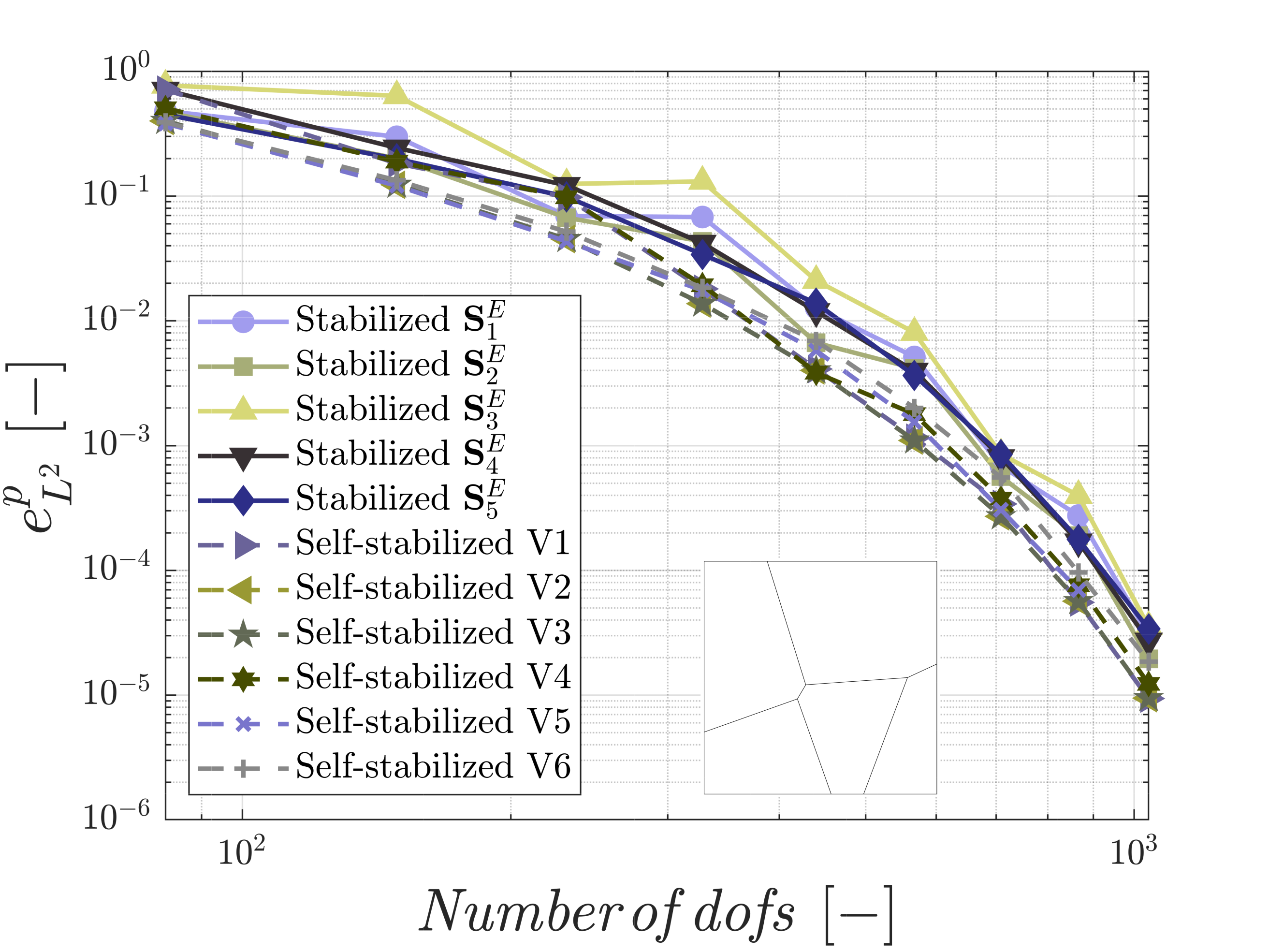}
	}     
    \caption{Stokes problem: performance of stabilized and self-stabilized formulations for Voronoi mesh.\label{fig:stability_stokes_voronoi}}
 \end{figure}

To complete the set of problems at hand, the investigation is now extended to the Stokes problem. The load term and the boundary conditions are chosen so that the exact solution is:
\begin{equation}
\begin{aligned}
    &\vel \plbr{x,y} = 
    \begin{bmatrix}
    +\sin\plbr{\pi x}^2 \sin\plbr{\pi y} \cos\plbr{\pi y}\\
    -\sin\plbr{\pi y}^2 \sin\plbr{\pi x} \cos\plbr{\pi x}
    \end{bmatrix},\\
    &\pres \plbr{x,y} = 
     \sin\plbr{\pi x}\cos\plbr{\pi y},
\end{aligned}    
    \label{eq:stokes_exactsolution}
\end{equation}
in the unit square $\surom = \plbr{0,1}^2$. The viscosity is set to $\nu=1$. Both the velocity and the pressure convergence are investigated:
\begin{equation*}
\begin{aligned}
    e_{\bm{\nabla}}^{\vel} =  \frac{ \plbr{\sum_{\el \in \globdom}  \norm{ \sqrt{\nu} \bm{\nabla} \plbr{ \vel - \PiokK \veld } }_{0,\el}^2  }^{\frac{1}{2}} }{ \norm{ \sqrt{\nu} \bm{\nabla} \vel }_{0,\el} },&&
    e_{L^2}^{\pres} =  \frac{ \plbr{\sum_{\el \in \globdom}  \norm{ \pres - \presd }_{0,\el}^2  }^{\frac{1}{2}} }{ \norm{\pres}_{0,\el} }.
\end{aligned}    
    \label{eq:stokes_error_velpres}
\end{equation*}
The analysis is carried out on the Voronoi mesh, with $\kord=2,\dots,10$ and the accuracy of the different formulations is compared using $\tau=1$ and $\mathrm{tol}=\mathrm{tol}_0$. 
A summary of the results is provided in \fig{stability_stokes_voronoi}. All the stabilized VEM and self-stabilized versions feature optimal error estimates for the velocity. A similar trend is observed for the pressure error with higher oscillations in the stabilized case.

\section{Equations with variable coefficients}\label{sec:VC}

Many practical problems involve coefficients that change spatially. An interesting example is given by modern composite laminates, where the elastic properties are a function of the position.

For most numerical methods, generic coefficients are interpolated in polynomial spaces. This approximation is adequate if either a high polynomial degree or a sufficiently fine mesh is considered.

In virtual element setting, the issue of dealing with variable coefficients has been addressed in \cite{daveiga2016virtual} for second-order elliptic problems: the authors discussed how the choice of the polynomial projection ($\Pinabla \trialfcn_i$ or $\Piokone \nabla\trialfcn_i$) may lead to suboptimal results.
Moreover, the stabilization term, which is only scaled by a measure of the coefficient, may further alter the results. This aspect is one of the motivations that justify the introduction of self-stabilized formulations.

A new strategy is proposed here to capture the influence of variable coefficients in both stabilized and self-stabilized VEM. The idea is to construct a polynomial approximation of virtual functions that combines the standard $\Piok$ and the features of the coefficient. This new strategy is built upon well-established virtual element spaces.
The acronym $\mathrm{VC}$-VEM is hereinafter employed to refer to this approach.

\subsection{Second-order elliptic problem}
\label{sec:2nd_vc}

The first example considered in this section is the second-order elliptic problem: \begin{equation}
    \begin{cases}
        \text{find } \unku \in H_0^1 \plbr{\surom} \text{ such that:}\\
        \bila \plbr{\testu,\unku} = \plbr{\testu,\force}_{\surom} \quad \forall \testu \in H_0^1 \plbr{\surom},
    \end{cases}
\label{eq:laplace_goveqs_2}
\end{equation}
where:
\begin{equation}
    \begin{aligned}
        \bila \plbr{\testu,\unku} = \int_{\surom} \nabla \testu^T \Acoeff \plbr{x,y} \nabla \unku \, \de \surom.
    \end{aligned}
\label{eq:laplace_bilabilf_2}
\end{equation}
Assume that the coefficient $\Acoeff \plbr{x,y}$ is uniformly positive definite and smooth within each element of the mesh. 

Two different discretizations of the above problem were proposed and studied in the framework of the standard virtual element space $\Vspaced$ introduced in \eq{VEMlocalspace_enhanced}, depending on the choice of the polynomial projector, namely:
\begin{equation}
\begin{aligned}
    &\bilad^\el \plbr{\trialfcn_i, \trialfcn_j } = 
        \bila^\el \plbr{ \Pinabla \trialfcn_i, \Pinabla \trialfcn_j } + \tau \KEscal \plbr{ \plbr{ I - \Pinabla } \trialfcn_i, \plbr{ I - \Pinabla } \trialfcn_j },\\
    &\bilad^\el \plbr{\trialfcn_i,\trialfcn_j } = 
        \plbr{ \Piokone \nabla\trialfcn_i, \Acoeff \plbr{x,y}\Piokone \nabla\trialfcn_j }_\el + \tau \KEscal \plbr{ \plbr{ I - \Pinabla } \trialfcn_i, \plbr{ I - \Pinabla } \trialfcn_j }.
\end{aligned}
\label{eq:general_beirao_projection}
\end{equation}
In the case of the Laplace equation, both constructions are robust, whenever $\Acoeff$ is constant. This is not anymore true in the case of variable coefficients: the discrete problem built upon $\Pinabla$ is subject to a strong deterioration in the order of convergence for $\kord\ge3$, as demonstrated in \cite{daveiga2016virtual}. Thus, the use of the $\Piokone$ was recommended for high-order schemes.

As mentioned, a new projection of $\Vspaced\plbr{\el}$ onto $\Pspace_\kord \plbr{\el}$ is proposed here, where $\Acoeff \plbr{x,y}$ is explicitly taken into account.
The stabilized virtual element method is constructed by considering, for $i=1,...,\spacedim$, the polynomial approximation $\PiA \trialfcn_i$ of $\trialfcn_i$ defined as the solution of:
\begin{equation}
    \begin{cases}
    \begin{aligned}
        &\bila^\el \plbr{ \PiA \trialfcn_i, \qpoly } = \bilVC\plbr{ \trialfcn_i, \qpoly } && \forall \qpoly \in \Pspace_\kord \plbr{\el},\\
        &P_0^\el \plbr{ \PiA \trialfcn_i, \qpoly } = P_0^\el \plbr{ \trialfcn_i, \qpoly }  && \forall \qpoly \in \Pspace_0 \plbr{\el},
    \end{aligned}    
    \end{cases}
    \label{eq:VEMH1projection_VC}
\end{equation}
where $\bilVC$ is defined as:
\begin{equation}
    \begin{aligned}
    \bilVC\plbr{ \trialfcn_i, \qpoly } \coloneqq &- \int_\el \Piok \trialfcn_i\,\diverg \plbr{\Acoeff \plbr{x,y} \nabla \qpoly}  \; \de \el\\
    &+ \sum_{e=1}^{\nedges} \int_{\bound} \trialfcn_i \cdot \Acoeff \plbr{x,y} \nabla \qpoly \normedge \; \de \bound\\
    &+ \sum_{e=1}^{\nedgescurv} \int_{\boundcurv} \trialfcn_i \cdot \Acoeff \plbr{x,y} \nabla \qpoly \normedgecurv \; \de \boundcurv.
    \end{aligned}
    \label{eq:B_matrix_integrationbyparts_VC_enhanced}
\end{equation}
The term $\bilVC\plbr{ \trialfcn_i, \qpoly }$ is a reasonable approximation of $\bila^\el\plbr{ \trialfcn_i, \qpoly }$ obtained by applying integration by parts and replacing $\trialfcn_i$ with $\Piok \trialfcn_i$ in the non-computable volume integral.
The expression above requires the partial derivatives of $\Acoeff \plbr{x,y}$ to be available.

All quantities in~\eqref{eq:VEMH1projection_VC}, which involve polynomials and the  coefficient $\Acoeff \plbr{x,y}$, can be computed by considering standard high-order quadrature rules.

Once $\PiA$ is computed, the discrete bilinear form is defined, as usual for VEM, by substituting the projection into the bilinear form and considering a suitable stabilization term (built on the standard $\Pinabla$), that is:
\begin{equation}\label{eq:VC_bilad}
    \bilad^\el \plbr{\trialfcn_i, \trialfcn_j } = 
        \bila^\el \plbr{ \PiA \trialfcn_i, \PiA \trialfcn_j } + \tau \KEscal \plbr{ \plbr{ I - \Pinabla } \trialfcn_i, \plbr{ I - \Pinabla } \trialfcn_j }.
\end{equation}
Since $\PiA$ is not an orthogonal projection, the (non-computable) terms of kind $\bila^\el \plbr{ \PiA \trialfcn_i, \plbr{ I - \PiA } \trialfcn_j }$ do not naturally vanish. Thus, they are neglected when defining $\bilad^\el$ in the stabilized VC-VEM setting.

The theoretical analysis of VC-VEM involves several technical aspects and is beyond the scope of the present work.

Concerning the self-stabilized formulations $\mathrm{V1}$, $\mathrm{V3}$, $\mathrm{V4}$ and $\mathrm{V6}$, an higher-order counterpart of $\PiA$ is computed within the framework of the discrete spaces defined in Section~\ref{subsec:Laplace_selfstabilization}.
Moreover, in lieu of the $L^2$ projection, the following operator is considered:
\begin{equation}
    \plbr{\Pi_{\kord-1}^{0,\Acoeff}\nabla\trialfcn_i,\Acoeff \plbr{x,y}\qpolyvec}_\el = \plbr{\nabla\trialfcn_i,\Acoeff \plbr{x,y}\qpolyvec}_\el \quad \forall \qpolyvec \in \sqbr{\Pspace_{\kord-1} \plbr{\el}}^2.
    \label{eq:new_l2}
\end{equation}
the $\mathrm{VC}$-VEM does not apply to formulations $\mathrm{V2}$ and $\mathrm{V5}$ as the internal degrees of freedom would have the form $(\trialfcn_i\,\diverg \plbr{\Acoeff \plbr{x,y} \nabla \qpoly})_\el$, which involves the possibly non-polynomial coefficient $\Acoeff \plbr{x,y}$.

\subsection{Linear Elasticity problem}

The VC-VEM is now extended to the linear elasticity problem defined in \eq{linearelasticity_governingeqs}. As presented in the previous example, the elasticity space in \eq{VEMspaceLinearElasticity_local} is considered. The following polynomial projection $\PiKVC\trialfcnvec_i$ of $\trialfcnvec_i$ ($i=1,...,\spacedim$) is then employed to define the polynomial contribution to the discrete bilinear form:
\begin{equation}
    \begin{cases}
    \begin{aligned}
         &\bila^\el \plbr{\PiKVC \trialfcnvec_i,\qpolyvec} = \bilVC\plbr{\trialfcnvec_i,\qpolyvec} && \forall \qpolyvec \in \sqbr{ \Pspace_{\kord} \plbr{\el}}^2,\\
         &P_0^\el \plbr{\PiKVC \trialfcnvec_i,\qpolyvec} = P_0^\el \plbr{\trialfcnvec_i,\qpolyvec} && \forall \qpolyvec \in RM\plbr{\el},
    \end{aligned}     
    \end{cases}
    \label{eq:energyproj_linearelasticityVC}
\end{equation}
where now:
\begin{equation}
    \bila^\el \plbr{\trialfcnvec_i,\qpolyvec} = \int_\el \strain\plbr{\trialfcnvec_i}^T \cost \plbr{x,y} \strain\plbr{\qpolyvec} \de \el,
    \label{eq:Cbilinearform_linearelasticity}
\end{equation}
and:
\begin{equation}
    \aligned
   \bilVC\plbr{\trialfcnvec_i,\qpolyvec} 
   \coloneqq &- \int_\el \PiZvec \trialfcnvec_i\,\bm{L} \plbr{\cost \plbr{x,y} \strain\plbr{\qpolyvec}}  \; \de \el\\
    &+ \sum_{e=1}^{\nedges} \int_{\bound} \trialfcnvec_i \cdot \hat{\bm{\sigma}}\plbr{\qpolyvec} \normedge \; \de \bound\\
    &+ \sum_{e=1}^{\nedgescurv} \int_{\boundcurv} \trialfcnvec_i \cdot \hat{\bm{\sigma}}\plbr{\qpolyvec} \normedgecurv \; \de \boundcurv.
   \endaligned
\end{equation}
The symbol $\hat{\bm{\sigma}}$ denotes the Cauchy stress tensor accounting for the spatial variability of $\cost \plbr{x,y}$.

\begin{remark}
    The projection operators $\PiA$ and $\PiKVC$ must be recomputed whenever the coefficient changes. This also applies when a mesh of uniform elements is considered, as the coefficient varies spatially. In principle, the $\mathrm{VC}$-VEM may result in a higher computational cost. However, the overhead is mitigated by the fact that the $\kord$-version of VEM is generally defined on coarse meshes.
\end{remark}

\section{Numerical Results: Equations with variable coefficients} \label{sec:NumericalResults_VC}

In this section, numerical tests are conducted on a set of problems characterized by the presence of variable coefficients. The aim of this study is to investigate the effectiveness of the procedure outlined in the previous section. In particular, the effect of coefficient variability on the accuracy of the method is analyzed.

\subsection{Second-order elliptic problem}

The problem under consideration is solved by adopting the same exact solution, mesh types and order $\kord$ used for the Laplace problem in Section \ref{subsec:Laplace_results_selfstab}.

To clarify the importance of appropriately handling variable coefficients, the simulations are run both with standard and $\mathrm{VC}$-VEM. The two bilinear forms in \eq{general_beirao_projection} are considered, together with $\KEs_3$ ($\tau=1$). Self-stabilized formulations $\mathrm{V3}$ and $\mathrm{V6}$ ($\mathrm{tol}=\mathrm{tol}_0$) are also chosen. The stabilization $\KEs_3$ is selected since it guarantees optimal accuracy with low computational burden (calculation of boundary integrals and matrix inversions are not required, as shown in Section \ref{sec:NumericalResults_SelfStab}). At the same time, self-stabilized $\mathrm{V3}$ and its counterpart $\mathrm{V6}$ are chosen to assess the robustness of $\mathrm{VC}$-VEM, for which the divergence-free 
property of the polynomial space is no longer applied due to the presence of the coefficient in the projection definition.

The error is evaluated as:
\begin{equation}
    e_{\nabla} =  \frac{ \plbr{\sum_{\el \in \globdom}  \norm{ \sqrt{\Acoeff \plbr{x,y}} \nabla \plbr{ \unku - \Piok \unkud } }_{0,\el}^2  }^{\frac{1}{2}} }{ \norm{\sqrt{\Acoeff \plbr{x,y}}\nabla \unku}_{0,\el} }.
    \label{eq:laplace_energyerror_VC}
\end{equation}
In the following, polynomial and trigonometric coefficients are considered.

 \subsubsection{Example 1: polynomial coefficients}
 The first example deals with polynomial coefficients in the form:
 \begin{equation}
    \Acoeff \plbr{x,y} = 
     \begin{bmatrix}
            x + 1 & 0\\
            0 & y + 1
     \end{bmatrix}.
     \label{eq:A_poly}
 \end{equation}

\begin{figure}[htbp]
        \centering
	\subfigure[Energy norm error. \label{fig:VC_poly_err_quad}]{
		\includegraphics[width=0.47\textwidth]{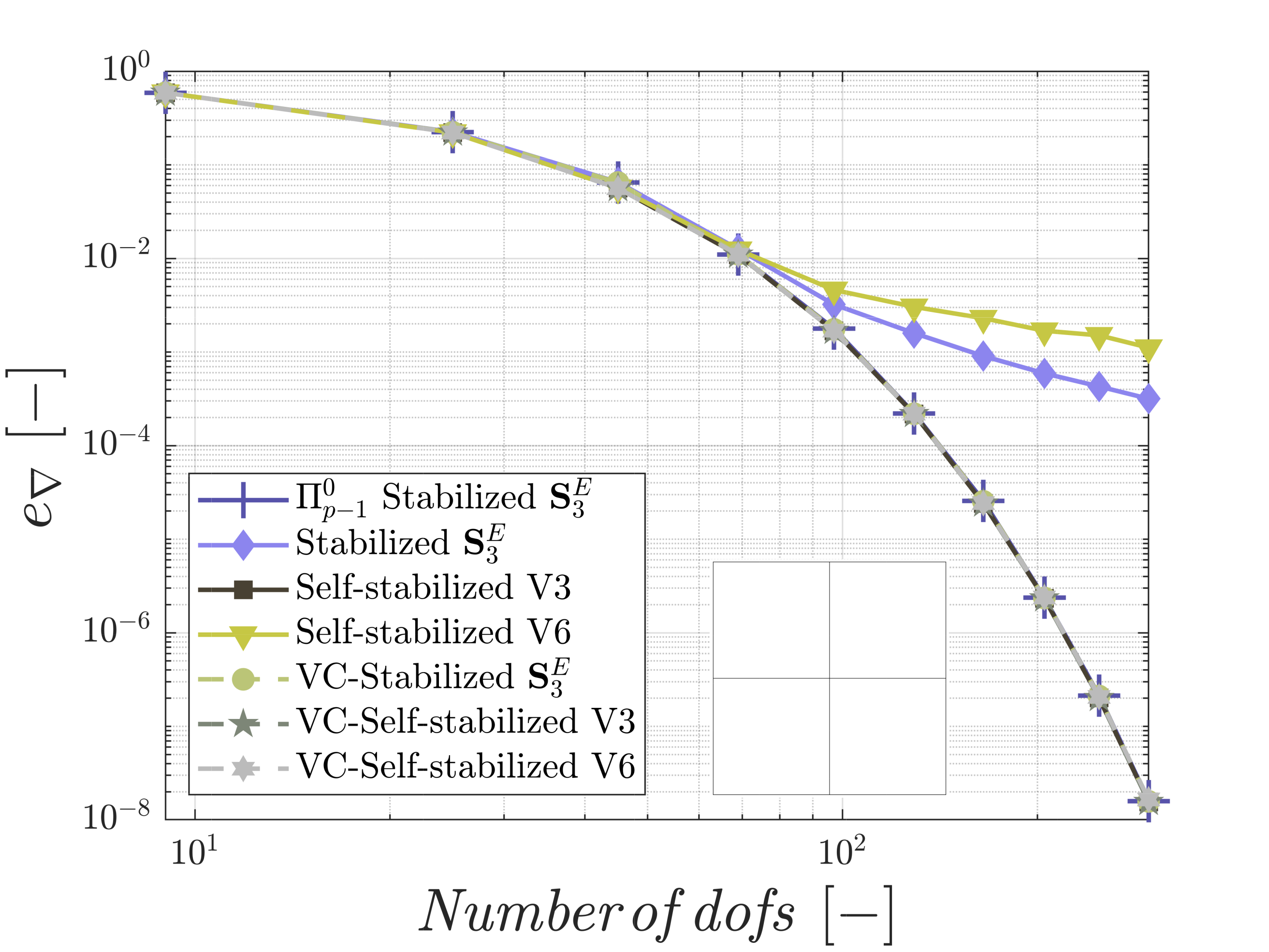}
	}
        \centering
        \hfill   
	\subfigure[Condition number. \label{fig:VC_poly_cond_quad}]{
		\includegraphics[width=0.47\textwidth]{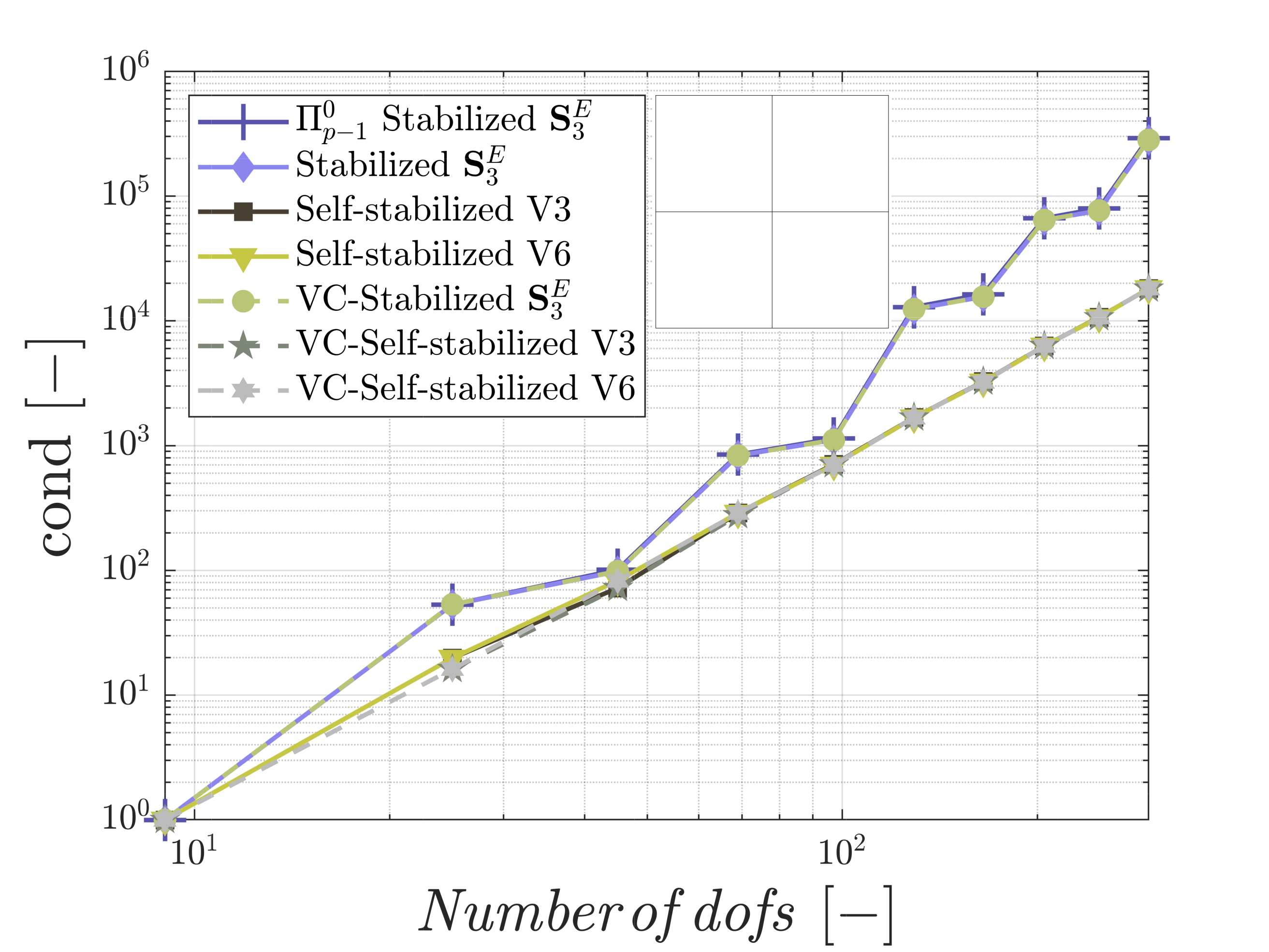}
	}     
    \caption{Second-order elliptic problem, example 1: performance comparison for quadrilateral mesh.\label{fig:VC_poly_quad}}
 \end{figure}

First, the quadrilateral mesh is considered. The results are shown in \fig{VC_poly_quad}. The accuracy of the solution, as seen from \fig{VC_poly_err_quad}, is the same for all the simulations up to the order $\kord=4$. However, when $p > 4$ some curves exhibit a diverging response. It is worth noting the excellent accuracy displayed by all the methods based on the newly introduced approach to consider coefficient-variability, even for $\kord \geq 5$.  
In this scenario, also the self-stabilized $\mathrm{V3}$ and the stabilized VEM with $\Piokone$ projection demonstrate excellent convergence. This behavior is attributed to the use of the $L^2$ projection, which tends to exhibits better accuracy in the presence of variable coefficients, as already discussed in \cite{daveiga2016virtual}. On the other hand, the formulation $\mathrm{V6}$, which does not consider variable coefficients, behaves as the stabilized VEM with $\Pinabla$. Regarding the stiffness matrix conditioning, reported in \fig{VC_poly_cond_quad}, no differences can be noted between the cases with or without variable coefficient projection.

\begin{figure}[htbp]
        \centering
	\subfigure[Energy norm error. \label{fig:VC_poly_err_voronoi}]{
		\includegraphics[width=0.47\textwidth]{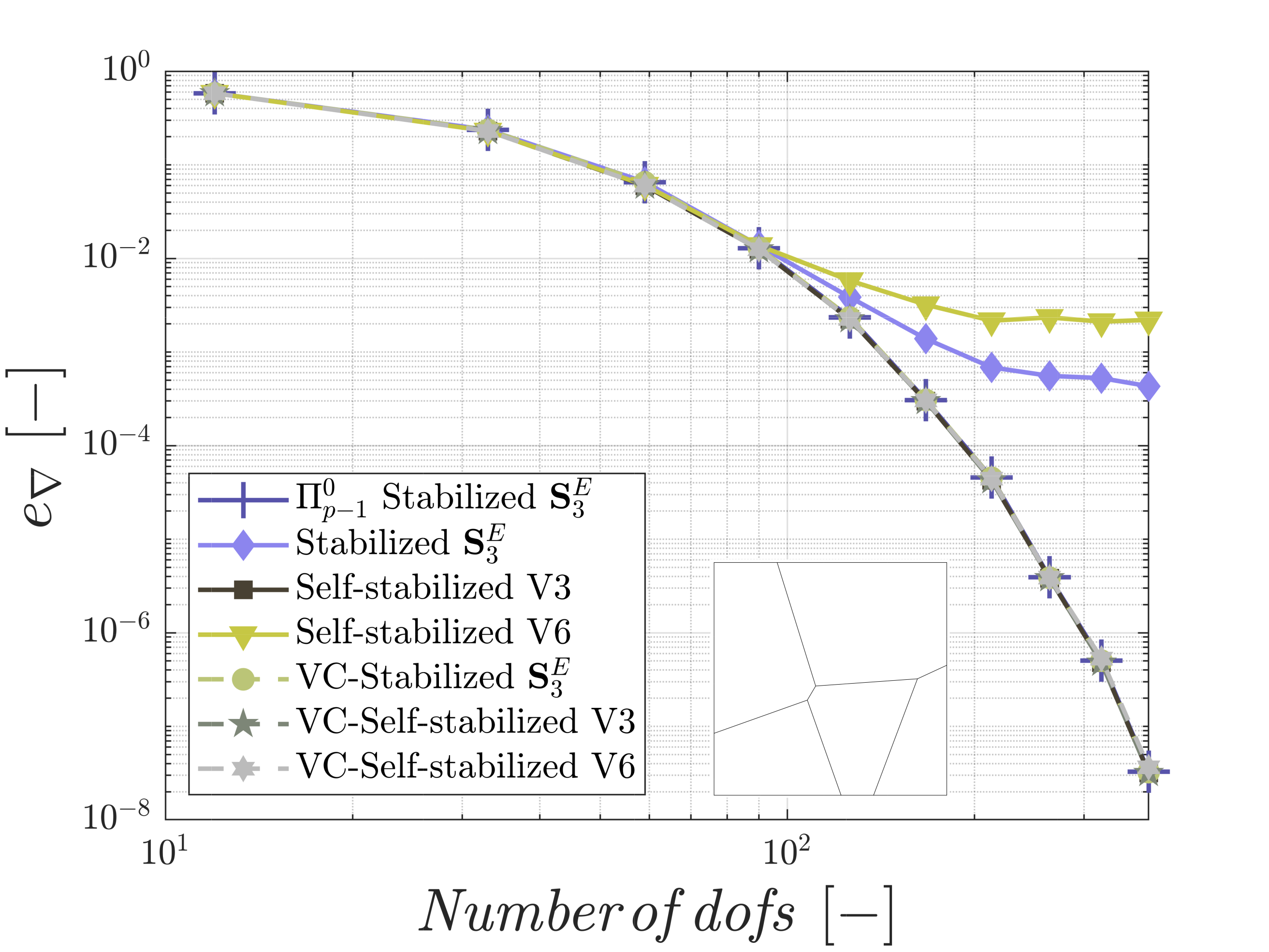}
	}
        \centering
        \hfill   
	\subfigure[Condition number. \label{fig:VC_poly_cond_voronoi}]{
		\includegraphics[width=0.47\textwidth]{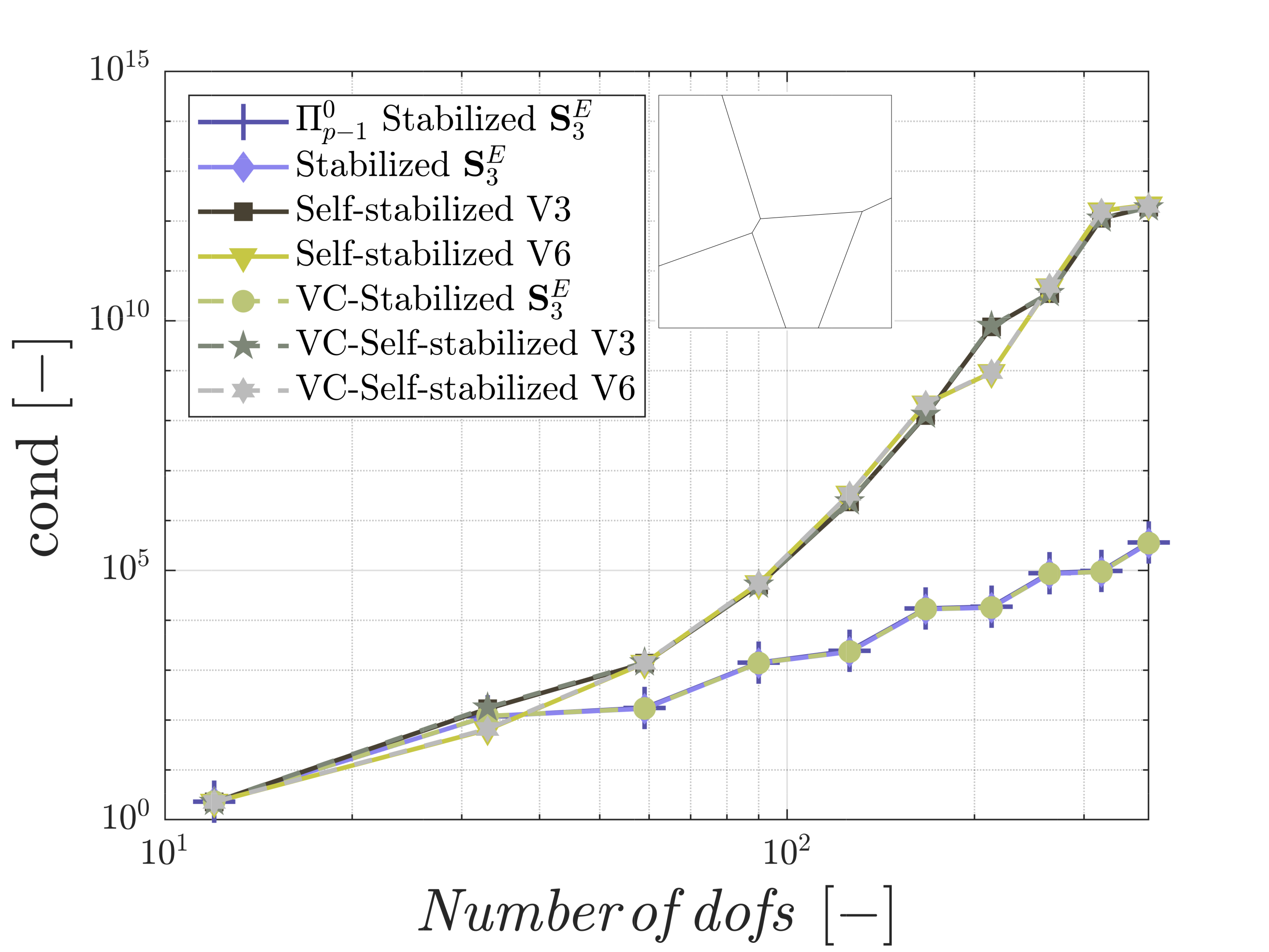}
	}     
    \caption{Second-order elliptic problem, example 1: performance comparison for Voronoi mesh.\label{fig:VC_poly_voronoi}}
 \end{figure}

 The results are shown in \fig{VC_poly_voronoi} for the Voronoi mesh. The same trend observed for the quadrilateral mesh can be noted in terms of accuracy. On the contrary, the condition number for the self-stabilized versions is significantly worse than the stabilized one, as shown in \fig{VC_poly_cond_voronoi}. Nevertheless, the presence of variable coefficients in the projection has no effect on the conditioning.

\begin{figure}[htbp]
        \centering
	\subfigure[Energy norm error. \label{fig:VC_poly_err_ottagoni}]{
		\includegraphics[width=0.47\textwidth]{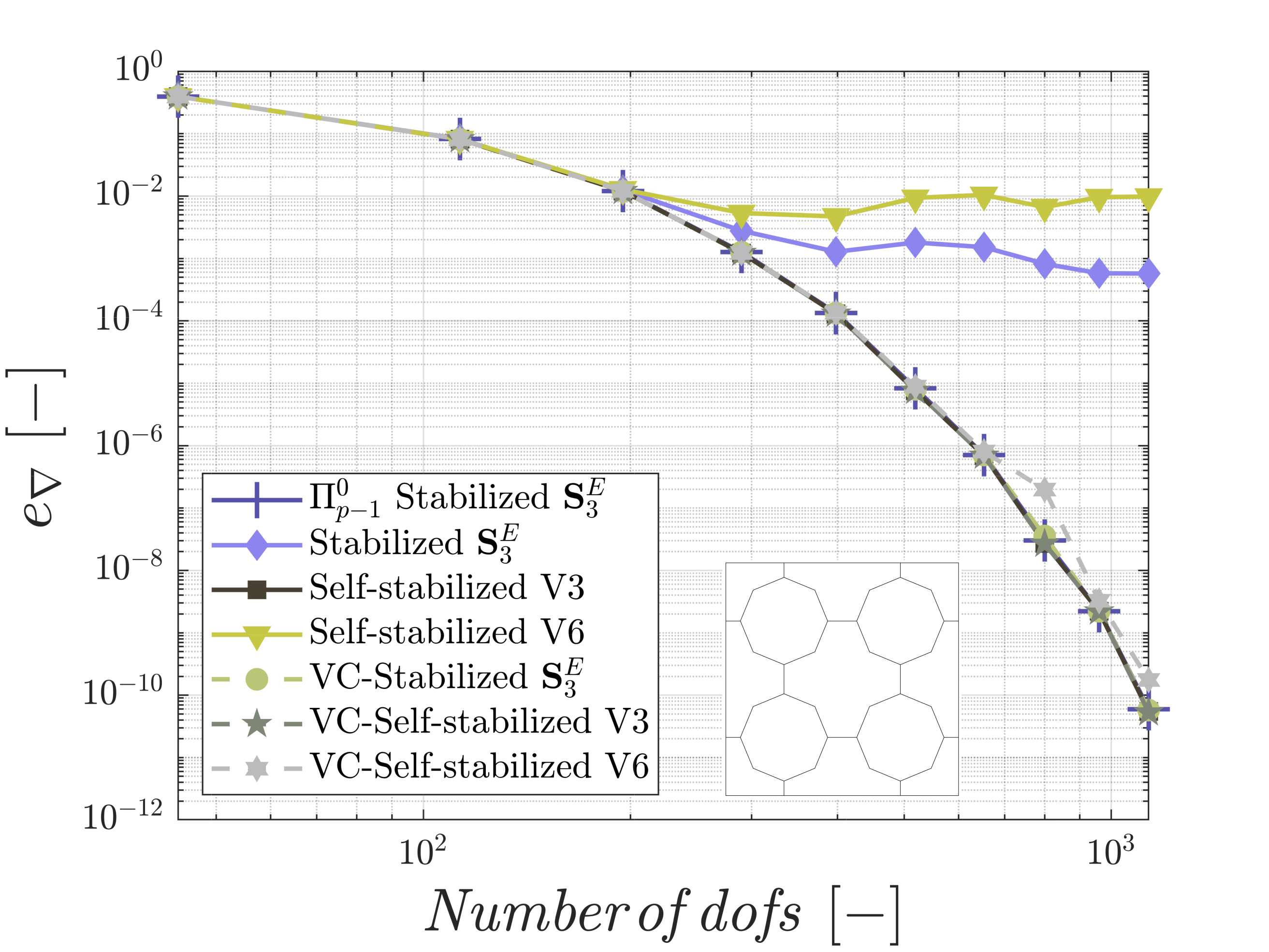}
	}
        \centering
        \hfill   
	\subfigure[Condition number. \label{fig:VC_poly_cond_ottagoni}]{
		\includegraphics[width=0.47\textwidth]{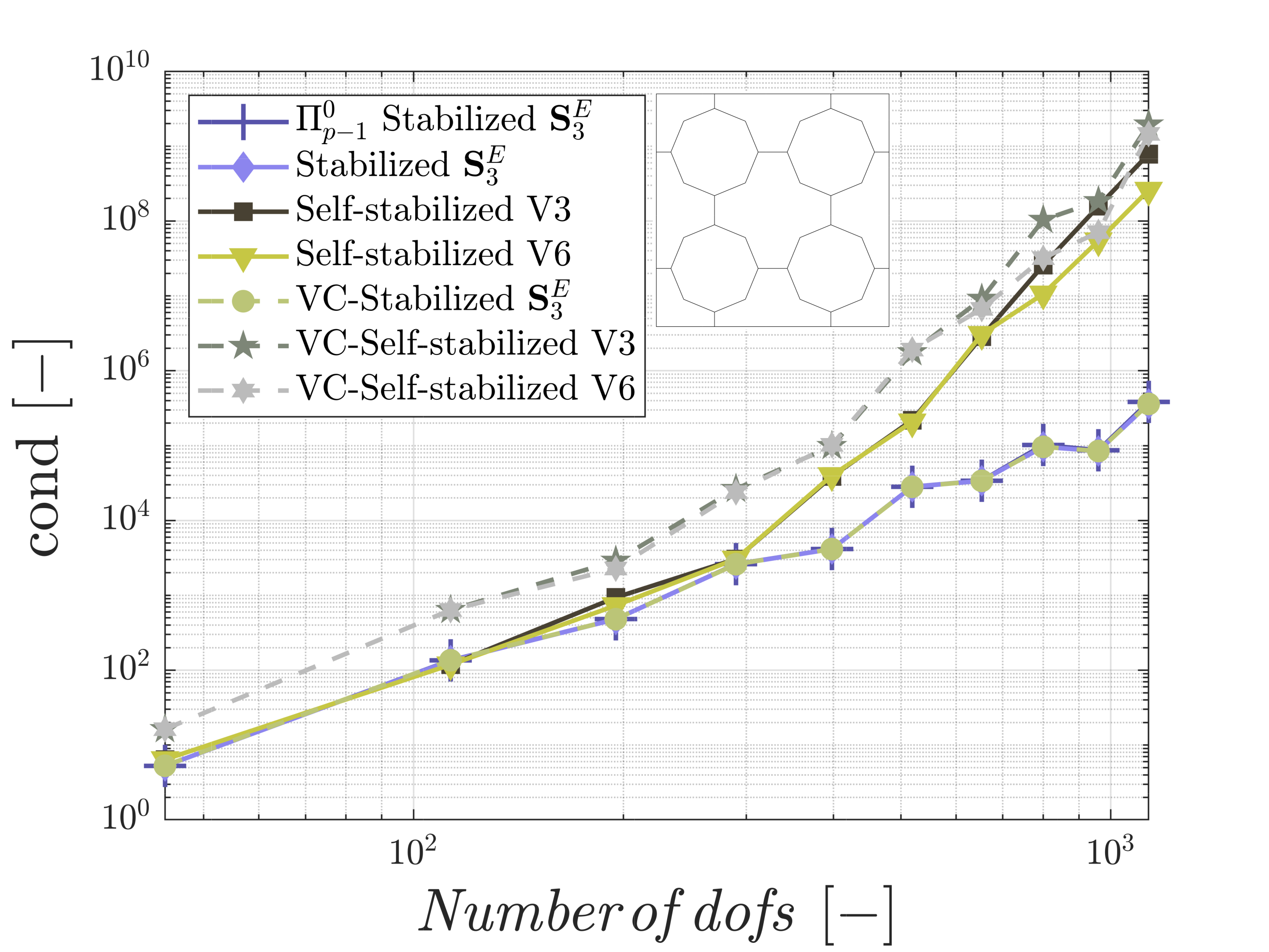}
	}     
    \caption{Second-order elliptic problem, example 1: performance comparison for octagonal mesh.\label{fig:VC_poly_ottagoni}}
 \end{figure}

 Lastly, the results for the octagonal mesh, reported in \fig{VC_poly_ottagoni}, show the same trend observed for the Voronoi mesh, both in accuracy and conditioning. 

\subsubsection{Example 2: trigonometric coefficient}

In the second example, the variation of the coefficients is assumed to be trigonometric, so:
\begin{equation}
    \Acoeff \plbr{x,y} = 
     \begin{bmatrix}
            \cos^4 \plbr{\pi x } + 1& \cos^2 \plbr{\pi x } \sin^2 \plbr{\pi y }\\
            \cos^2 \plbr{\pi x } \sin^2 \plbr{\pi y } & \sin^4 \plbr{\pi y } + 3
     \end{bmatrix}.
     \label{eq:A_trig}
 \end{equation}

\begin{figure}[htbp]
        \centering
	\subfigure[Energy norm error. \label{fig:VC_trig_err_quad}]{
		\includegraphics[width=0.47\textwidth]{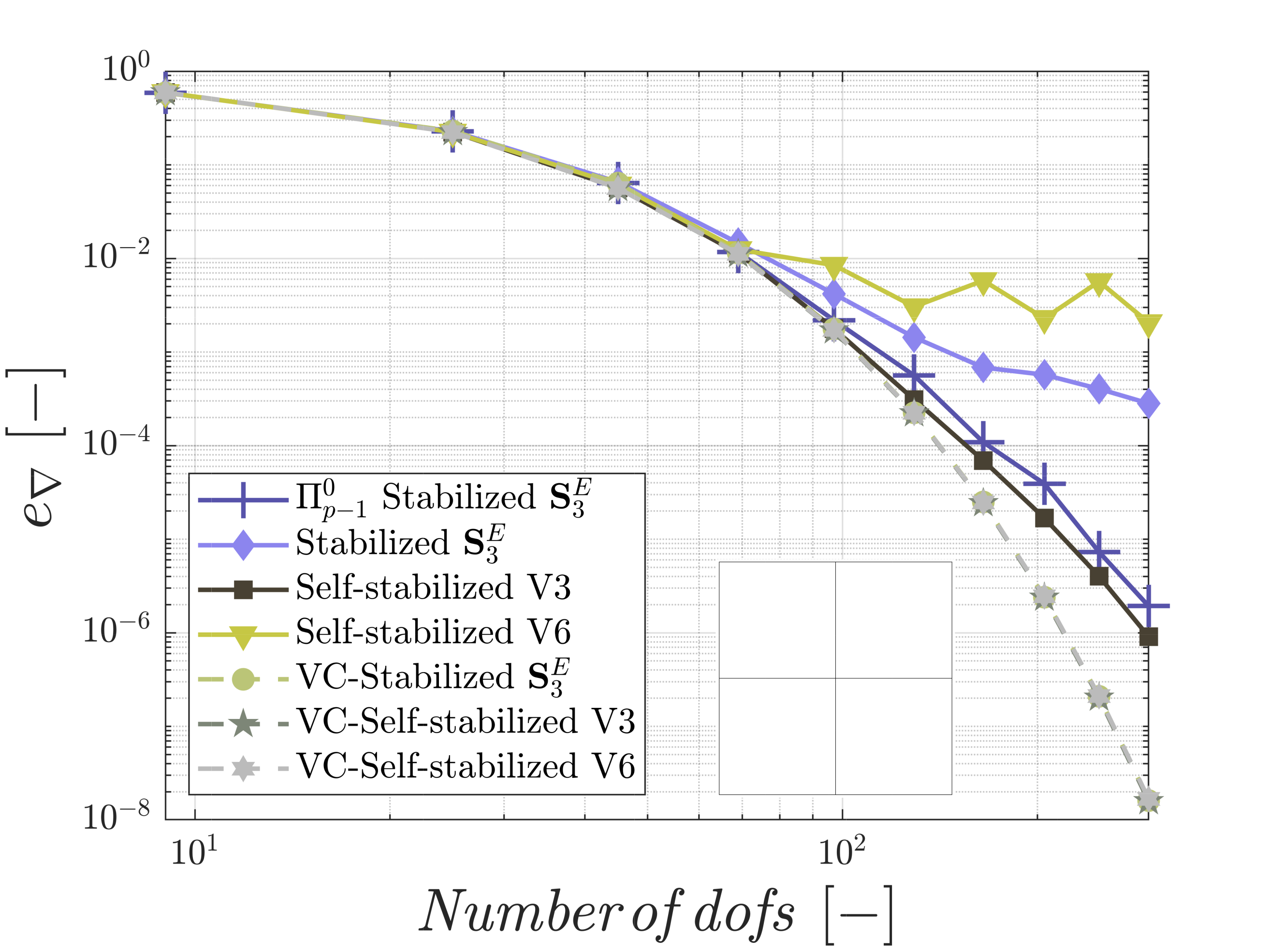}
	}
        \centering
        \hfill   
	\subfigure[Condition number. \label{fig:VC_trig_cond_quad}]{
		\includegraphics[width=0.47\textwidth]{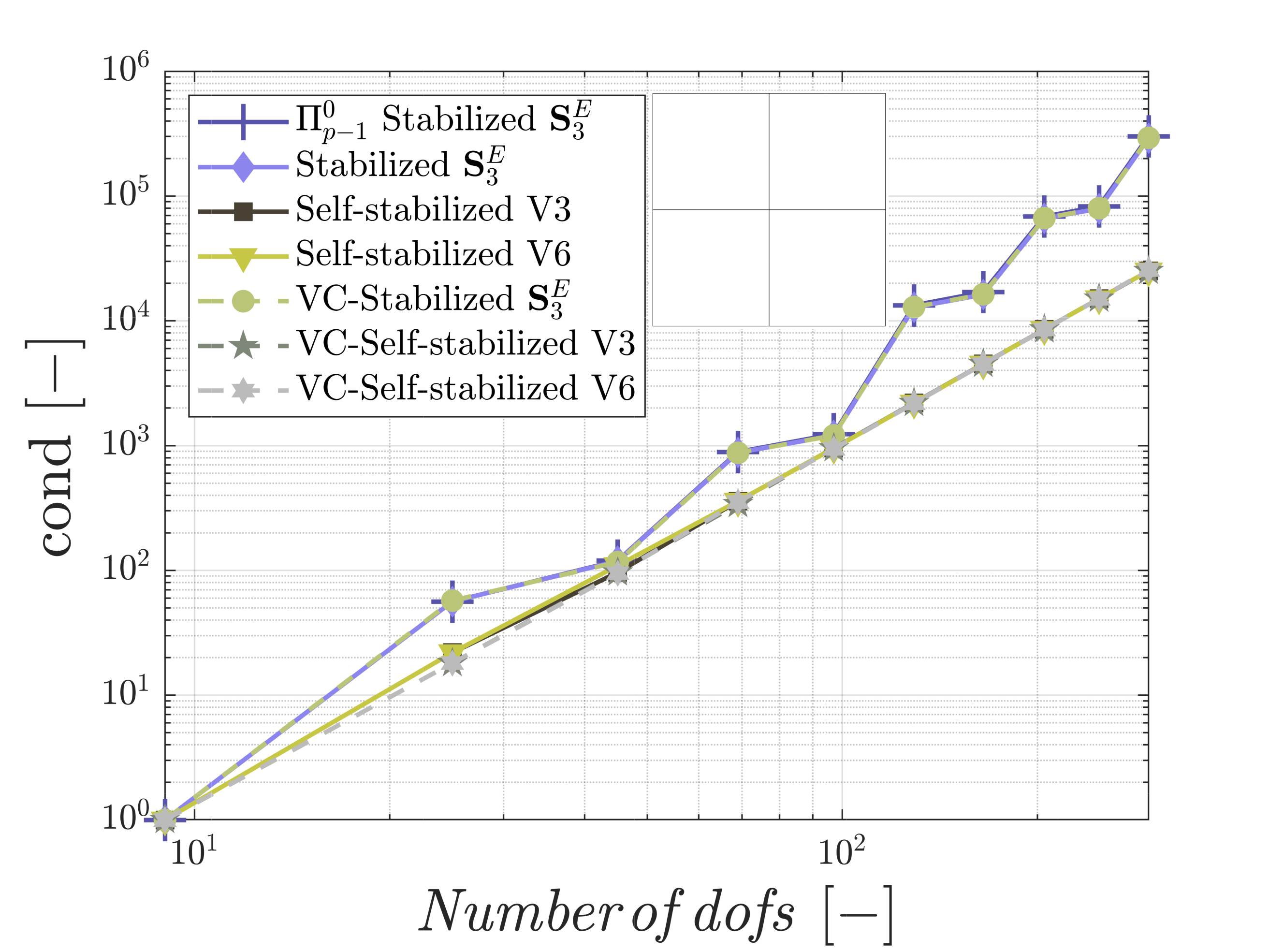}
	}     
    \caption{Second-order elliptic problem, example 2: performance comparison for quadrilateral mesh.\label{fig:VC_trig_quad}}
 \end{figure}

A summary of the results is provided in \figto{VC_trig_quad}{VC_trig_ottagoni} for the three different meshes. 
For all mesh types, the condition number features the same behavior observed in the previous test case. It is then concluded that the novel formulation with variable coefficients does not affect the conditioning. For what concerns the accuracy, all $\mathrm{VC}$-VEM simulations demonstrate optimal convergence. On the contrary, by considering this second coefficient variation, the self-stabilized $\mathrm{V3}$ and the stabilized VEM with $\Piokone$ projection no longer preserve the accuracy displayed in the previous test case. It is important to highlight that on more structured meshes, i.e. the quadrilateral and the octagonal ones, the latter formulations show better accuracy compared to the Voronoi mesh. 

\begin{figure}[htbp]
        \centering
	\subfigure[Energy norm error. \label{fig:VC_trig_err_voronoi}]{
		\includegraphics[width=0.47\textwidth]{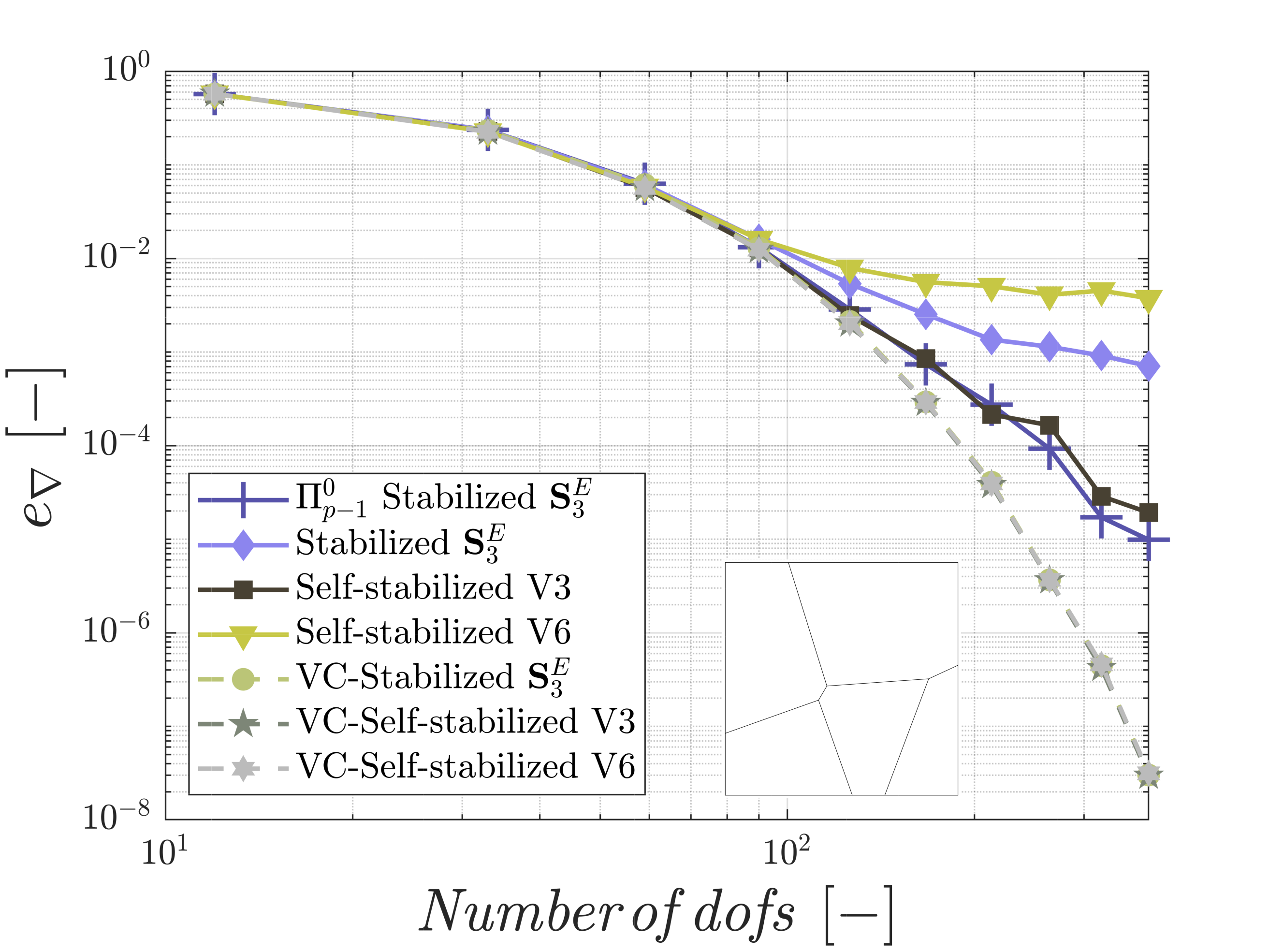}
	}
        \centering
        \hfill   
	\subfigure[Condition number. \label{fig:VC_trig_cond_voronoi}]{
		\includegraphics[width=0.47\textwidth]{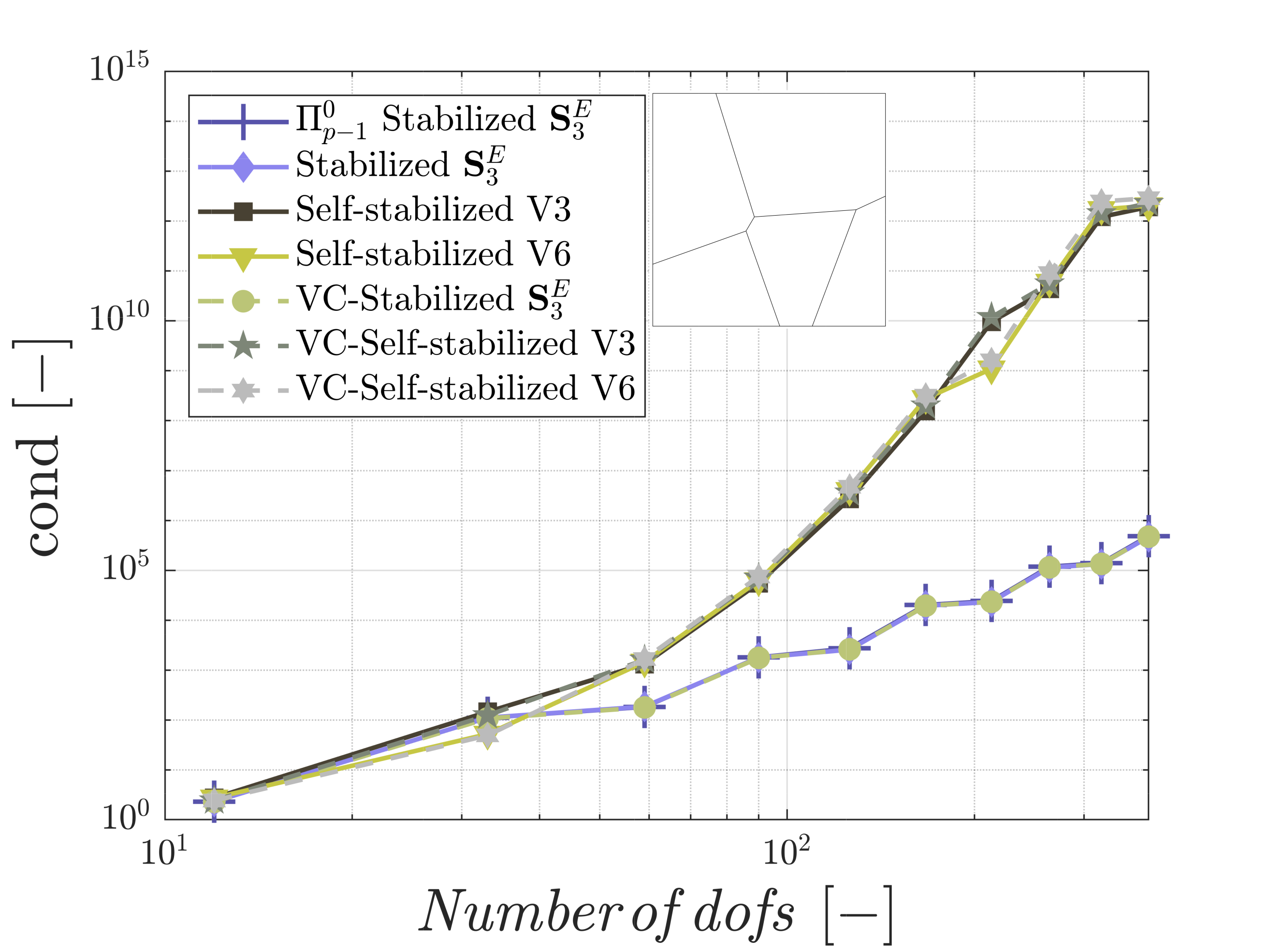}
	}     
    \caption{Second-order elliptic problem, example 2: performance comparison for Voronoi mesh.\label{fig:VC_trig_voronoi}}
 \end{figure}

 \begin{figure}[htbp]
        \centering
	\subfigure[Energy norm error. \label{fig:VC_trig_err_ottagoni}]{
		\includegraphics[width=0.47\textwidth]{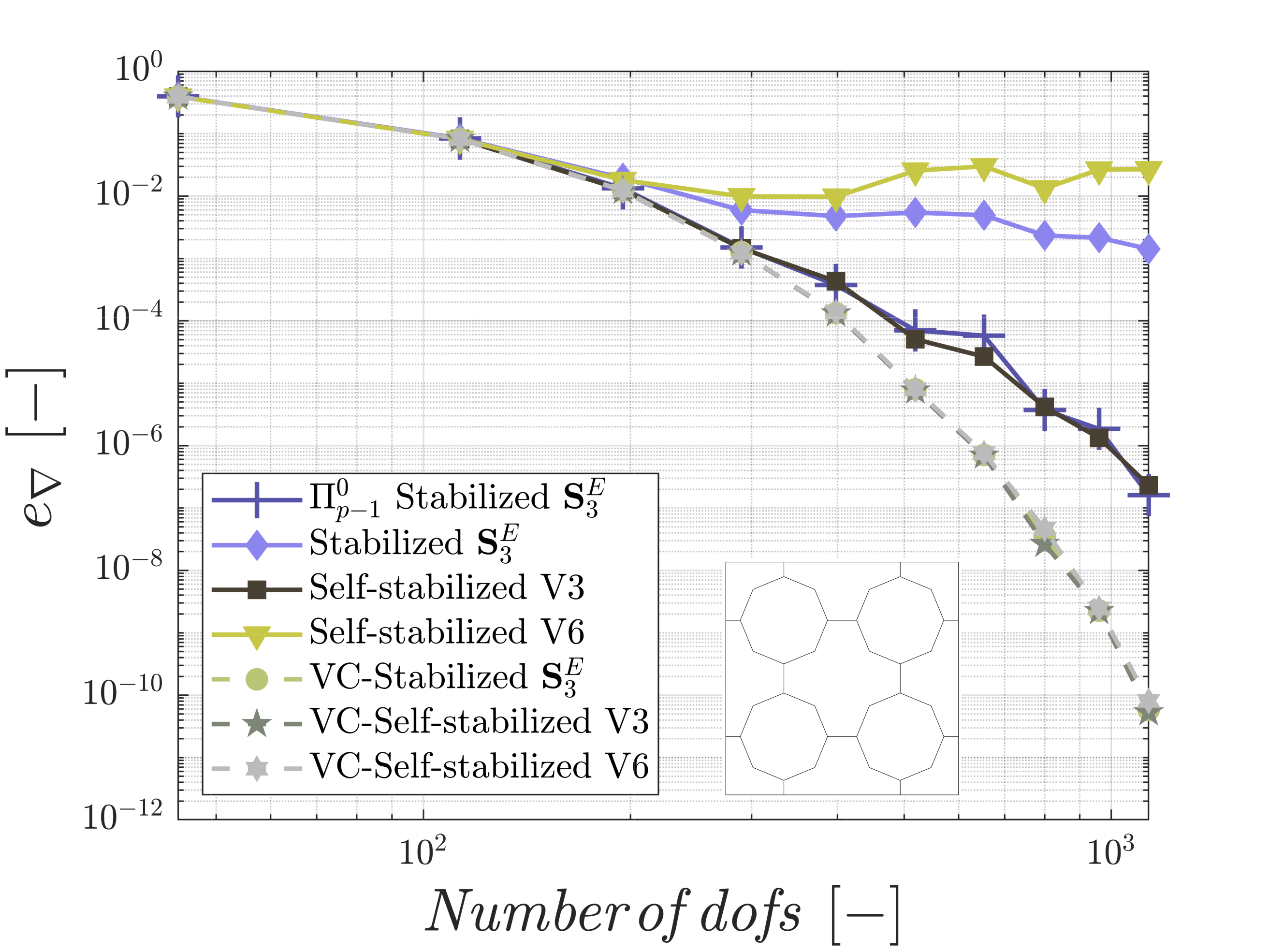}
	}
        \centering
        \hfill   
	\subfigure[Condition number. \label{fig:VC_trig_cond_ottagoni}]{
		\includegraphics[width=0.47\textwidth]{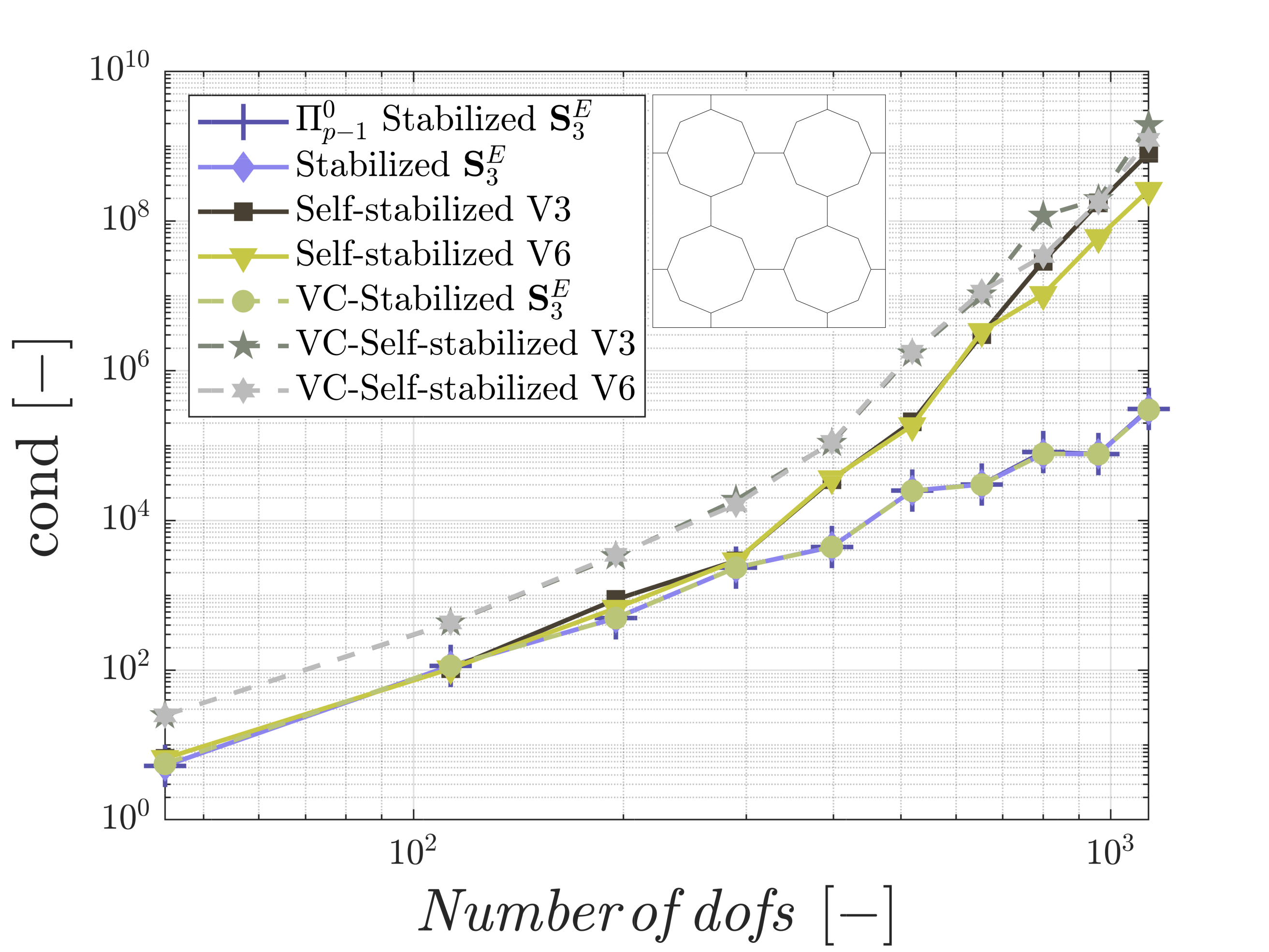}
	}     
    \caption{Second-order elliptic problem, example 2: performance comparison for octagonal mesh.\label{fig:VC_trig_ottagoni}}
 \end{figure}

\subsection{Linear elasticity problem}

Lastly, the results for the linear elasticity are presented. The same exact solution, mesh types and order $\kord$ of Section \ref{subsec:Linerelasticity_results_selfstab} are considered. By following the idea outlined in \cite{reddy2019virtual}, the projection operators for standard VEM are computed by considering a piece-wise constant approximation of the coefficient. In each element, $\cost$ is approximated by its average at the integration points. The error is measured as in \eq{linearelasticity_l2error_dispenergy} and the analysis is conducted with the stabilized VEM $\KEs_3$ ($\tau=0.5$) and the self-stabilized formulation $\mathrm{V1}$ ($\mathrm{tol}=\mathrm{tol}_0$), which has not been previously investigated.

Finally, the variable elastic properties are assumed to obey this polynomial behavior:
\begin{equation}
\begin{aligned}
& E_1 \plbr{x,y} = E_2 \plbr{x,y} = 72000 \plbr{1+100\plbr{x-x_c}^4+\plbr{y-y_c}^2},\\ 
& \nu_{12} = \nu_{21} = 0.33,\qquad
G_{12} \plbr{x,y} = \frac{E_1 \plbr{x,y}}{2\plbr{1+\nu_{12}}},
\end{aligned}    
 \label{eq:const_law_VC1}
\end{equation}
where $\plbr{x_c,y_c}=\plbr{0.5,0.5}$.

\begin{figure}[htbp]
        \centering
	\subfigure[$L^2$ error. \label{fig:VC_quad_err_disp_quad}]{
		\includegraphics[width=0.47\textwidth]{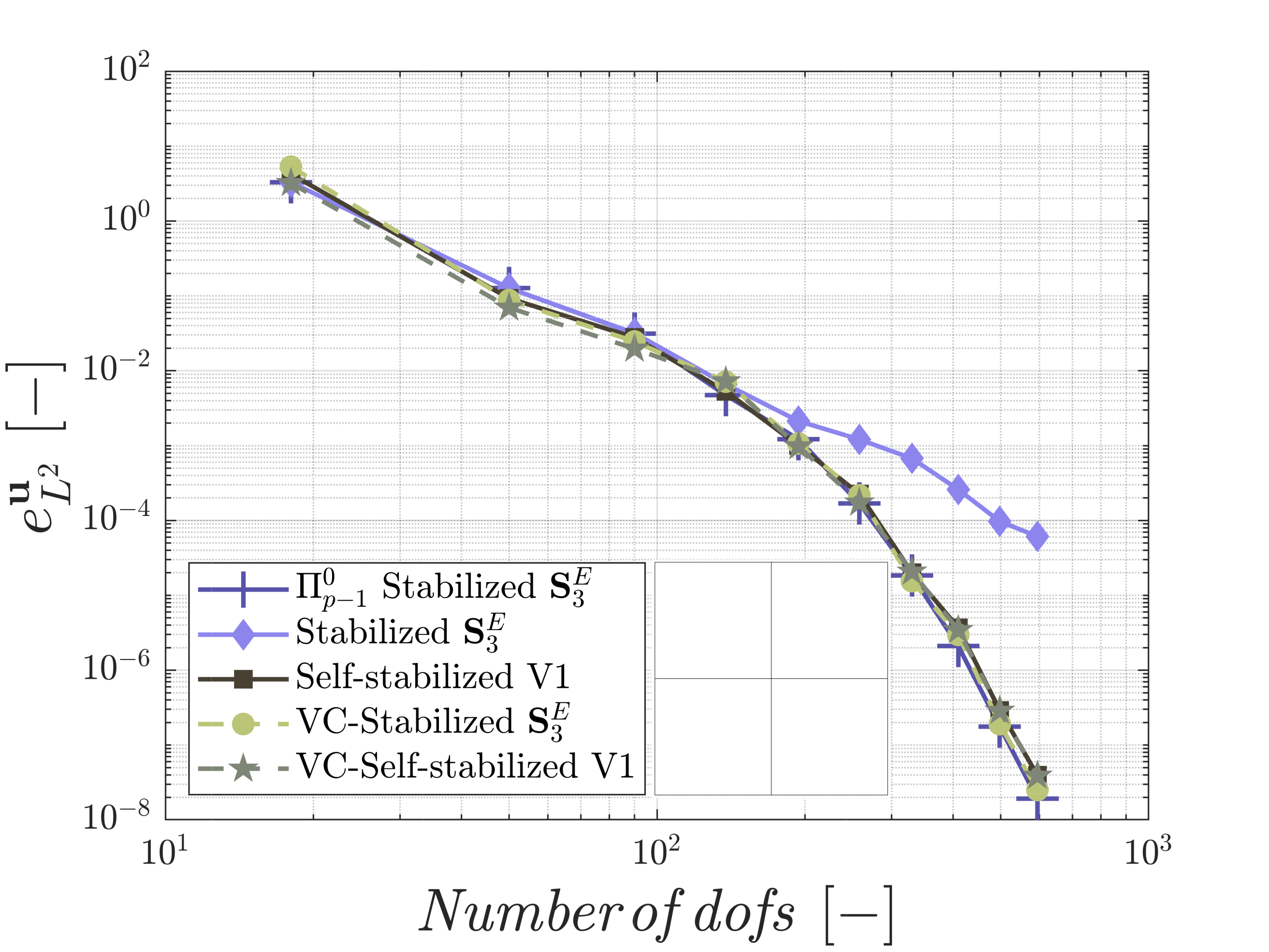}
	}
        \centering
        \hfill   
	\subfigure[Energy norm error. \label{fig:VC_quad_err_stress_quad}]{
		\includegraphics[width=0.47\textwidth]{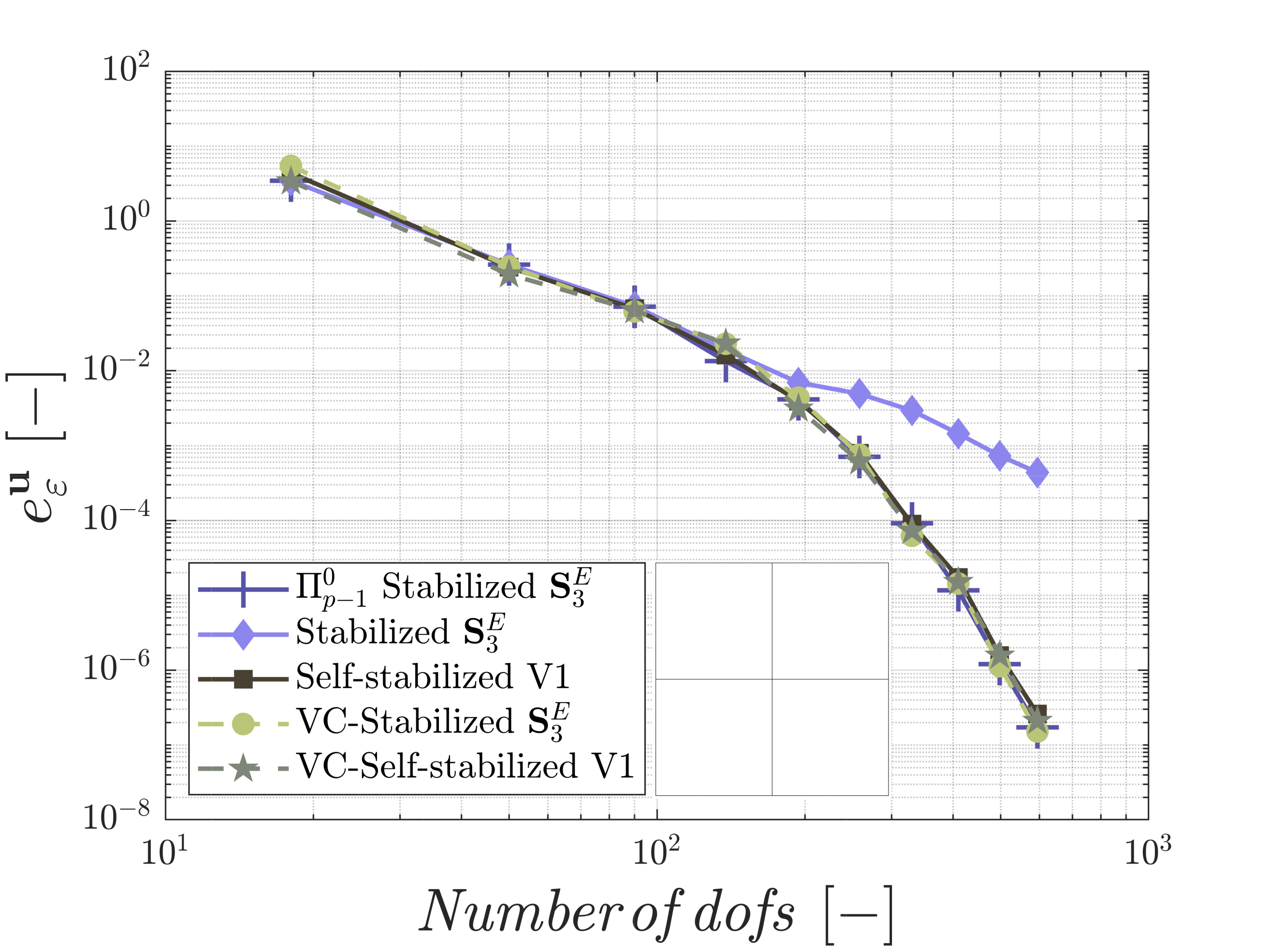}
	}     
    \caption{Linear elasticity problem: performance comparison for quadrilateral mesh.\label{fig:VC_quad_quad}}
 \end{figure}

For the quadrilateral mesh, results are shown in \fig{VC_quad_quad}. The analysis conducted with the newly developed formulation retains optimal accuracy. This is true both for the $L^2$ and for the energy error. Also, the self-stabilized $\mathrm{V1}$ and the approach presented in \cite{daveiga2016virtual}, extended to the present problem, provide an adequate description of the solution.

\begin{figure}
        \centering
	\subfigure[$L^2$ error. \label{fig:VC_quad_err_disp_curv}]{
		\includegraphics[width=0.47\textwidth]{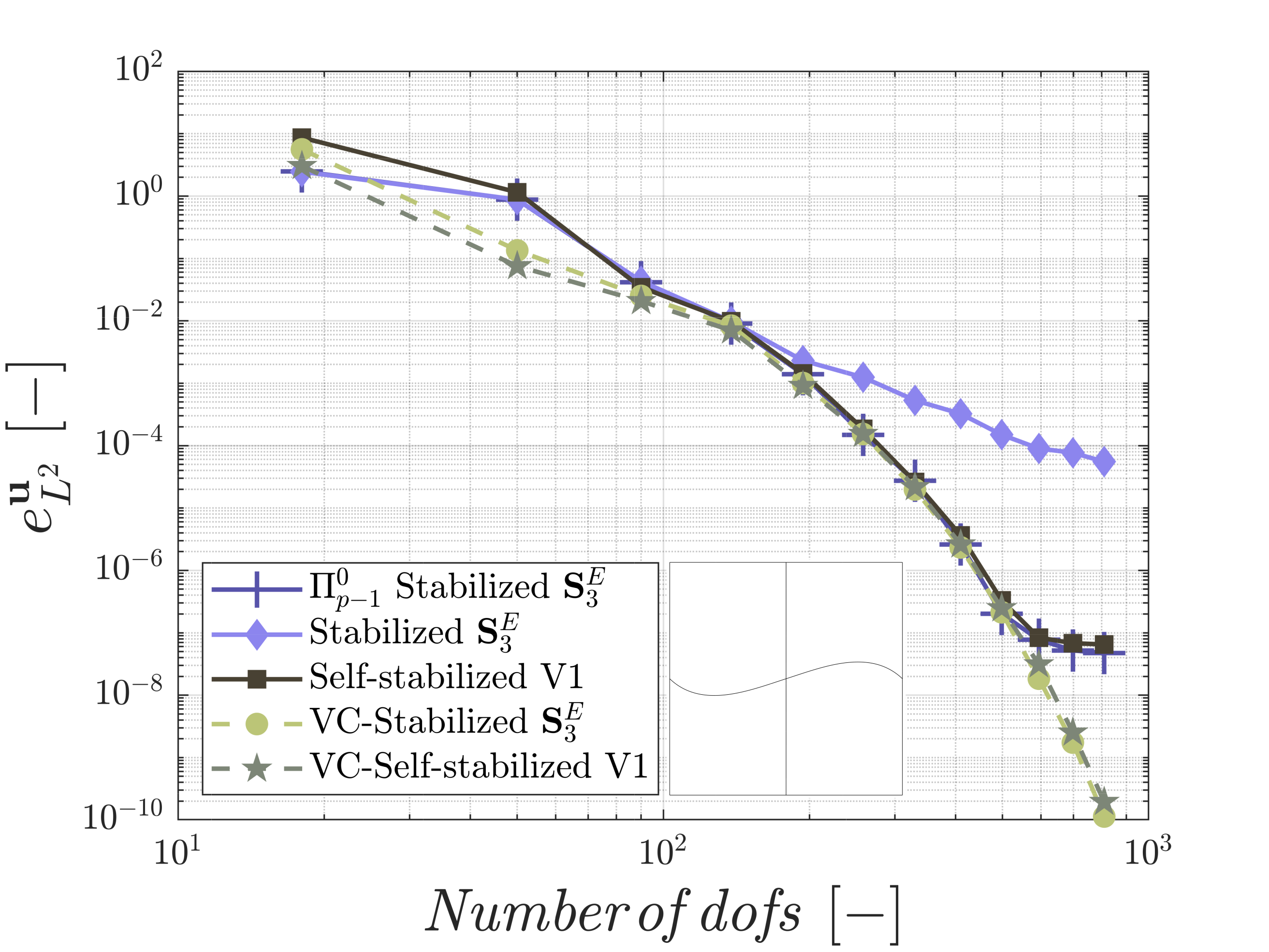}
	}
        \centering
        \hfill   
	\subfigure[Energy norm error. \label{fig:VC_quad_err_stress_curv}]{
		\includegraphics[width=0.47\textwidth]{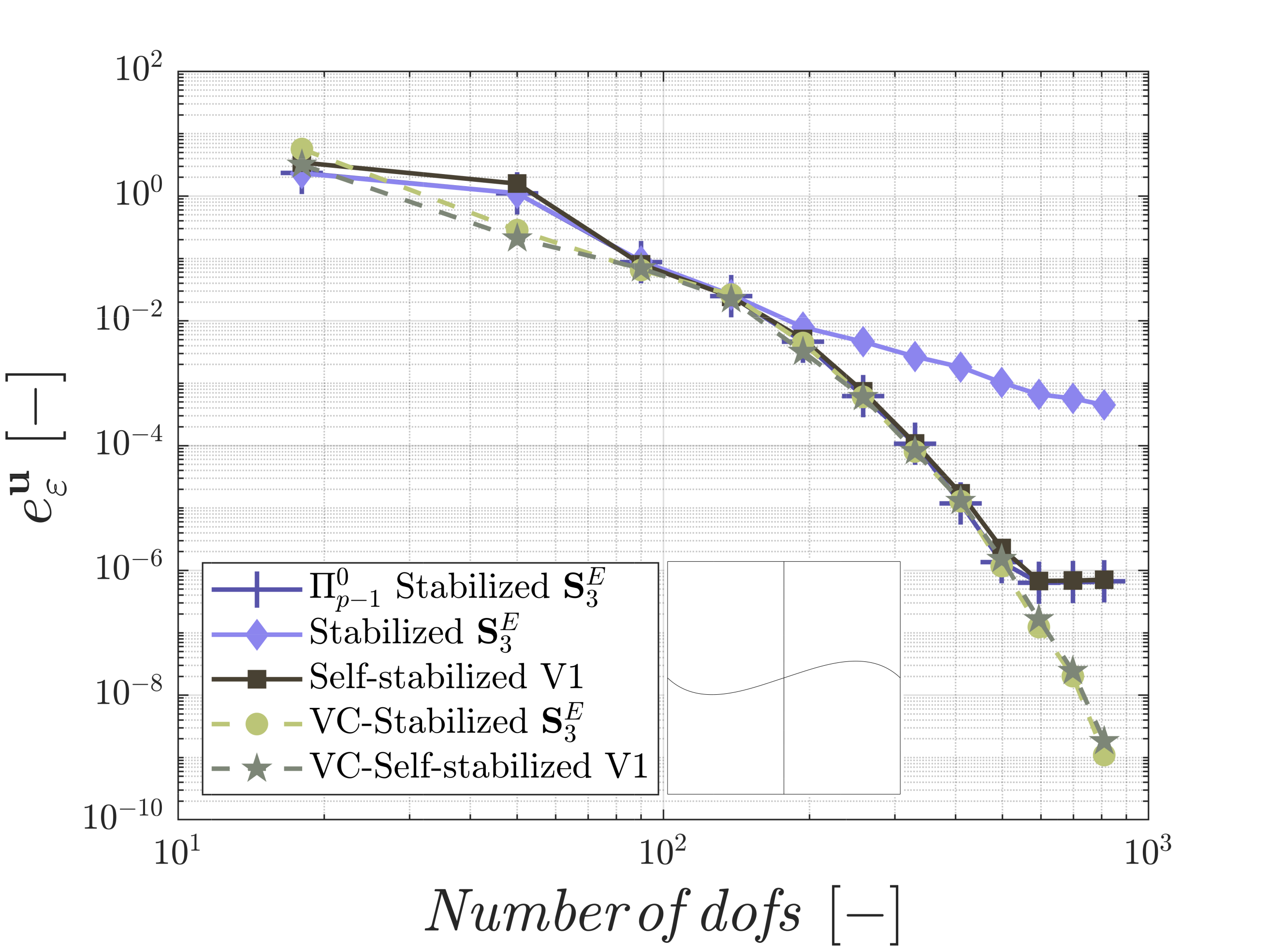}
	}     
    \caption{Linear elasticity problem: performance comparison for B\'ezier-edge mesh.\label{fig:VC_quad_curv}}
 \end{figure}

The results obtained with the B\'ezier-edge mesh are presented in \fig{VC_quad_curv} and are similar to the previous case. The approach discussed in \cite{daveiga2016virtual} features optimal accuracy up to order $\kord=9$ with a slight decrease for $\kord\geq10$. Therefore, to determine whether this accuracy loss is an isolated behavior or not, this analysis is conducted up to order $\kord=12$. It is clear that for $\kord \geq 10$ the $\Piokone$ stabilized VEM and the self-stabilized $\mathrm{V1}$ experience a flattening of the error curve for this mesh type. On the contrary, the $\mathrm{VC}$-VEM simulations maintain the accuracy both in the $L^2$ and in the energy error. As in the previous tests, the standard stabilized VEM features larger errors as the order $\kord$ increases. 

\section{Conclusions} \label{sec:Conclusions}

This paper provided an in-depth comparison between stabilized and self-stabilized formulations for the $\kord$-version of the VEM.

The investigation has been performed by keeping in mind potential applications in the field of variable stiffness plates with curvilinear stringers. In this case, the constitutive equations feature variable coefficients and their discretization might require curvilinear meshes for an accurate representation of the geometry.

In the first part of the work, the benchmarking has been carried out by considering the Laplace equation, the linear elasticity problem and the Stokes problem in the case of constant coefficients. The effect of the stabilization parameter and of the additional polynomial order required by self-stabilized formulations has been analyzed in terms of accuracy and condition number. The results can be summarized by the following highlights. 

\begin{itemize}
    \item Self-stabilized formulations feature optimal accuracy without the issue of choosing a suitable stabilization term/parameter, although at the price of worse conditioning and higher computational cost caused by higher-order polynomial projections.
    \item The choice of the tolerance for determining the additional polynomial order $\laug$ does not affect the accuracy of self-stabilized formulations. At the same time, a ``good'' tolerance should balance the conditioning and the computational cost: a ``large'' tolerance provides larger values of $\laug$ so that the stiffness matrix is far from being singular. Again, larger $\laug$ results in higher computational cost.
    \item The stabilized VEM formulations retain at least the same accuracy as the self-stabilized ones, while preserving better conditioning and lower computational cost.
\end{itemize}

The second part of the paper focused on the case of variable coefficients, with the introduction of a new virtual element approach, denoted by $\mathrm{VC}$-VEM. It consists in defining a new polynomial projection which combines the $L^2$ projector and the features of the coefficient in the framework of standard VEM spaces. The newly introduced formulation was compared with the standard approaches and provided promising results. The comparison can be summarized by the following highlights.

\begin{itemize}
    \item As already observed in \cite{daveiga2016virtual}, the stabilized VEM constructed with the operator $\Pinabla$ has a strong loss of accuracy for $\kord \geq 3$.
    \item The stabilized VEM constructed by the operator $\Piokone$ provides optimal results for polynomial coefficients, while a slight accuracy deterioration is observed at high approximation orders ($\kord \geq 6$) in the case of the trigonometric coefficient.
    \item The $\mathrm{VC}$-VEM appears to be robust for all the coefficients considered in the present investigation, with no deterioration even in more complex cases. 
    \item Since $\mathrm{VC}$-VEM takes into account the coefficients in the definition of the local space, the projection operator should be re-computed whenever the coefficient changes. The possibly higher computational cost is mitigated in a $\kord$-VEM setting, where very coarse meshes are considered.
\end{itemize}

Building on the obtained results, future developments will address the extension of the present formulation to geometrically nonlinear problems in elasticity.

\section*{Acknowledgements}
P. P. Foligno has received funding from the KAUST Visiting Student Program for developing this research. The support is gratefully acknowledged.

D. Boffi and F. Credali are members of the INdAM--GNCS research group.

\bibliographystyle{abbrv}
\bibliography{Bibliography}

\end{document}